\documentclass[12pt,a4paper]{article}
\usepackage[centertags]{amsmath}
\usepackage{amsfonts}
\usepackage{amssymb}
\usepackage{amsthm}
\usepackage{graphics}
\usepackage{newlfont}

 \theoremstyle{plain}
\newtheorem{thm}{Theorem}[section]
\newtheorem{cor}[thm]{Corollary}
\newtheorem{lem}[thm]{Lemma}
\newtheorem{prop}[thm]{Proposition}
\theoremstyle{definition}
\newtheorem{defn}[thm]{Definition}
\theoremstyle{remark}
\newtheorem{rem}[thm]{Remark}
\numberwithin{equation}{section}

\newcommand{\bq }{\begin{equation}}
\newcommand{\eq }{\end{equation}}
\newcommand{\bbb }{\begin{eqnarray}}
\newcommand{\eee }{\end{eqnarray}}
\newcommand{\bb }{\begin{eqnarray*}}
\newcommand{\ee }{\end{eqnarray*}}
\newcommand{\ed }{\end{document}}
\input{tcilatex}

\begin{document}

\title{ geometry and reducibility of 
induction
of  Hopf group co-algebra}
\author{ A.S.Hegazi$^1$ F.Ismail$^2$ and M.M. Elsofy$^3$ \\
%EndAName
{\small {$^1$ Department of Mathematics, Faculty of Science, Mansoura
University, Egypt.}}\\
{\small {Email: Hegazi@mans.edu.eg}}\\
{\small {$^2$ Department of Mathematics, Faculty of Science, Cairo
University, Egypt. \:\:\:\;\;\;\;\mbox{} }}\\
{\small {Email: fatma-ism@hotmail.com}}\\
{\small {$^3$ Department of Mathematics, Faculty of Science, Fayoum
University, Egypt. \;\;\;\mbox{} }}\\
{\small {Email: melsoffyy@yahoo.com}}}
\date{}
\maketitle

\begin{abstract}
In this work we study the induction (induced and coinduced)theory for Hopf
group coalgebra. We define a substructure B of a Hopf group coalgebra $H$,
called subHopf group coalgebra. Also, we have introduced the definition of
Hopf group suboalgebra and group coisotropic quantum subgroup of $H$. The
properties of the algebraic structure of the induced and coinduced are
given. Moreover, a framework of the geometric interperation and simplicity
theory of such representation strructure are stuided.
\end{abstract}

\section{{\protect\large {Introduction}}}

The induced representation of quantum group (quasitriangular Hopf algebra
[M] and [Mon]) is introduced by Gonzalez-Ruiz, L. A. Ibort [G-I] and Ciccoli
[C]. Recently, Hegazi, Agawany, F. Ismail and I. Saleh in [H], have
succeeded to give a new algebric structure for quantum subgroups and
subquantum groups, that is a subspase $B$ of a bialgebra $H$ is called a
sub-bi-algebra if the restriction to $B$ of the structure of $H$ turns $B$
into bialgebra. But a bi-sub-algebra of $H$ is a different thing. It is pair
$(B,\pi )$ consisting of another bialgebra and surjective homomorphism $\pi
:H\rightarrow B.$ A relation between these two notions are given if the
bialgebra allows a decomposition of the form $H=B\oplus I$ with $I$ a
bi-ideal in $H$. The projection $\pi $ into B yields a bi-sub-algebra $%
(B,\pi )$ in the above sense. Conversely if $B\subseteq H$ is both a
sub-bi-algebra and a bi-sub-algebra of $H$ such that $\pi ^{2}=\pi ,$ then $%
H=B+Ker$ $\pi $ and we have a decomposition of the form $H=B\oplus I$. These
construction made us able to introduce the induced representation of Hopf
algebra in each case. That is a $B$-comodule for a sub-Hopf algebra $B$ of a
Hopf algebra $H$ gives rise to an induced $H$-Hopf module and that a
representation of a quantum subgroup $(B,\pi )$ of a quantum group $H$
induced a Hopf representation of H. This procedure realizes a quantum group
induced representation.\newline

Recently, Turaev [T]and Virelizier [V-V1] give use a new definition for a
generalization of Hopf algebra, for Hopf algebra structure see [S] and
[Mon], Hopf group coalgebra. These generalization gives us quantum group
structure this structure is of great importance in Homotopy theory in
quantum field theory [T].\newline

In this paper, unless otherwise, every thing takes place over a field $K$,
and $K$-space means vector space over $K$. A map $f$ \ from a space $V$ into
a space $W$ always means linear map over K. The tensor product $V\otimes W$
is understood to be $V\otimes _{K}W$ , $I:V\rightarrow V$ always denotes the
identity map, and the transposition map $\tau :V\otimes W\rightarrow
W\otimes V$ is defined by $\tau (v\otimes w)=w\otimes v$ \ for $v\in V,w\in
W.$ Let $f:C\rightarrow D$ be a map. Then $f$ $^{\ast }:D^{\ast }\rightarrow
C^{\ast }$ is a map, where $f^{\ast }(\phi )(c)=\phi (f(c))$ for all $\phi
\in D^{\ast },c\in C.$\newline

Recently, Virelizier [V-V1] studied the algebraic properties of the Hopf $%
\pi -$ coalgebra. Now, let us give some basic definitions about Hopf $\pi -$%
coalgebra. For group $\pi $, a $\pi $-coalgebra is a family $C=\{C_{\alpha
}\}_{\alpha \in \pi }$ of $K$-spaces endowed with a family $K$-maps
(comultiplication) $\Delta =\{\Delta _{\alpha ,\beta }:C_{\alpha \beta
}\rightarrow C_{\alpha }\otimes C_{\beta }\}_{\alpha ,\beta \epsilon \pi }$,
$K$-map (counit) $\varepsilon :C_{1}\rightarrow k$ such that the following
diagrams are commute:
\begin{equation*}
\begin{tabular}{l}
$%
\begin{tabular}{lll}
$\quad C_{\alpha \beta \gamma }$ &
\begin{tabular}{l}
$\Delta _{\alpha ,\beta \gamma }$ \\
$\longrightarrow $%
\end{tabular}
& $C_{\alpha }\otimes C_{\beta \gamma }$ \\
$\Delta _{\alpha \beta ,\gamma }\downarrow $ &  & $\downarrow I_{\alpha
}\otimes \Delta _{\beta ,\gamma }$ \\
$\quad C_{\alpha \beta }\otimes C_{\gamma }$ &
\begin{tabular}{l}
$\longrightarrow $ \\
$\Delta _{\alpha ,\beta }\otimes I_{\gamma }$%
\end{tabular}
& $C_{\alpha }\otimes C_{\beta }\otimes C_{\gamma }$%
\end{tabular}%
\ $ \\
$\quad (\Delta _{\alpha ,\beta }\otimes I_{\gamma }$ $)\Delta _{\alpha \beta
,\gamma }=(I_{\alpha }\otimes \Delta _{\beta ,\gamma })\Delta _{\alpha
,\beta \gamma }$%
\end{tabular}%
\end{equation*}
\begin{equation*}
\begin{tabular}{l}
\begin{tabular}{lll}
& $C_{1}\otimes C_{\alpha }$ &  \\
$\epsilon \otimes I_{\alpha }\swarrow $ &  &  \\
$\quad K\otimes C_{\alpha }$ & $\quad \uparrow \Delta _{1,\alpha }$ &  \\
$\qquad \quad \sim _{\alpha }\nwarrow $ &  &  \\
& $\quad C_{\alpha }$ &
\end{tabular}
$\quad $\quad $\quad $\quad
\begin{tabular}{ll}
$C_{\alpha }\otimes C_{1}$ &  \\
& $\searrow I_{\alpha }\otimes \epsilon $ \\
$\Delta _{\alpha ,1}\uparrow $ & $C_{\alpha }\otimes K$ \\
& $\nearrow \sim _{\alpha }$ \\
$\qquad C_{\alpha }$ &  \\
&
\end{tabular}
$\quad $\quad $\quad $\quad $\quad $\quad \\
$\quad $\quad $\ \ \ \ \ \ \ \ \ \quad \ \ \ \ (\epsilon \otimes I_{\alpha
})\Delta _{1,\alpha }=$ $\sim _{\alpha }=(I_{\alpha }\otimes \epsilon
)\Delta _{\alpha ,1}$%
\end{tabular}%
\ .
\end{equation*}

A Hopf $\pi -$coalgebra is a $\pi $-coalgebra $(H=\{H_{\alpha }\}_{\alpha
\in \pi },\Delta ,\epsilon )$ with a family $S=\{S_{\alpha }:H_{\alpha
}\rightarrow H_{\alpha ^{-1}}\}_{\alpha \in \pi }$ of $K$-maps such that

\begin{enumerate}
\item $H_{\alpha }$ is an algebra with multiplication $\mu _{\alpha }$ and
unit $\eta _{\alpha }(1_{K})$ for all $\alpha \in \pi ;$

\item $\Delta _{\alpha ,\beta }$ $,$ $\epsilon $ are algebra maps for all $%
\alpha ,\beta \in \pi ,$

\item The antipode $S$ satisfy \newline
$\qquad \qquad\qquad \qquad \mu _{\alpha }(S_{\alpha ^{-1}}\otimes I_{\alpha
})\Delta _{\alpha ^{-1},\alpha }=1_{\alpha }\epsilon =\mu _{\alpha
}(I_{\alpha }\otimes S_{\alpha ^{-1}})\Delta _{\alpha ,\alpha ^{-1}}$
\end{enumerate}

A Hopf $\pi $-coalgebra $H$ is of finite type when every $H_{\alpha }$ is
finite dimensional. Note that it does not means that ${\oplus _{\alpha \in
\pi }}H_{\alpha }$ is finite dimensional (unless $H_{\alpha }=0$ for all but
a finite number of $\alpha \in \pi $). The notion of Hopf $\pi -$coalgebra
is not self dual. Note, $(H_{1},\Delta _{1,1},\epsilon )$ is a (classical)
Hopf algebra.\newline

Let $C$ be a $\pi $-coalgebra. A right $\pi $-comodule over $C$ is a family $%
M=\{M_{\alpha }\}_{\alpha \in \pi }$ of $K$-spaces endowed with a family $%
\theta =\{\theta _{\alpha ,\beta }:M_{\alpha \beta }\rightarrow M_{\alpha
}\otimes C_{\beta }\}_{\alpha ,\beta \epsilon \pi }$ of $K$-maps such that
the following diagrams are commute:
\begin{equation*}
\begin{tabular}{lll}
$%
\begin{tabular}{lll}
$\qquad M_{\alpha \beta \gamma }$ & $\quad \overset{\theta _{\alpha ,\beta
\gamma }}{\longrightarrow }$ & $\quad M_{\alpha }\otimes C_{\beta \gamma }$
\\
$\theta _{\alpha \beta ,\gamma }\downarrow $ &  & $\qquad \downarrow
I_{\alpha }\otimes \Delta _{\beta ,\gamma }$ \\
$\quad M_{\alpha \beta }\otimes C_{\gamma }$ & \quad $\overset{\theta
_{\alpha ,\beta }\otimes I_{\gamma }}{\longrightarrow }$ & $C_{\alpha
}\otimes C_{\beta }\otimes C\gamma $%
\end{tabular}%
\ $ &  &
\begin{tabular}{l}
\end{tabular}
\begin{tabular}{ll}
$\quad M_{\alpha }$ &  \\
& $\searrow \theta _{\alpha ,1}$ \\
$\quad \sim _{M_{\alpha }}\downarrow $ & $\quad M_{\alpha }\otimes C_{1}$ \\
& $\swarrow I_{\alpha }\otimes \epsilon $ \\
$M_{\alpha }\otimes K$ &
\end{tabular}
\\
$\quad (\theta _{\alpha ,\beta }\otimes I_{\gamma }$ $)\theta _{\alpha \beta
,\gamma }=(I_{\alpha }\otimes \Delta _{\beta ,\gamma })\theta _{\alpha
,\beta \gamma }$ &  & $\qquad (I_{\alpha }\otimes \epsilon )\theta _{\alpha
,1}=\sim _{M_{\alpha }}$%
\end{tabular}%
\end{equation*}

A right Hopf $\pi $-comodule over $H$ is a right $\pi $-comodule $%
M=\{M_{\alpha }\}_{\alpha \in \pi }$ by coaction $\theta =\{\theta _{\alpha
,\beta }\}_{\alpha ,\beta \epsilon \pi }$ such that

\begin{enumerate}
\item $M_{\alpha }$ is $H_{\alpha }-$module by action $\rho _{\alpha }$ for
all $\alpha \in \pi $

\item The following diagram, for all $\alpha ,\beta ,\gamma \in \pi $, is
commute:%
\begin{equation*}
\begin{tabular}{l}
\begin{tabular}{llllll}
$\ \ \ \quad M_{\alpha \beta }\otimes H_{\alpha \beta }$ & $\overset{\rho
_{\alpha \beta }}{\rightarrow }$ & $\quad M_{\alpha \beta }$ & $\overset{%
\theta _{\alpha ,\beta }}{\rightarrow }$ & $M_{\alpha }\otimes H_{\beta }$ &
\\
$\theta _{\alpha ,\beta }\otimes \Delta _{\alpha ,\beta }\downarrow $ &  &
&  & $\ \ \ \ \ \uparrow \rho _{\alpha }\otimes \mu _{\beta }$ &  \\
$\ \ \ \ (M_{\alpha }\otimes H_{\beta })\otimes (H_{\alpha }\otimes H_{\beta
})$ &  & $\overset{I_{\alpha }\otimes \tau \otimes I_{\beta }}{\rightarrow }$
&  & $(M_{\alpha }\otimes H_{\alpha })\otimes (H_{\beta }\otimes H_{\beta })$
&
\end{tabular}
\\
\begin{tabular}{l}
$\qquad i.e.,$ $\theta _{\alpha ,\beta }$ $\rho _{\alpha \beta }=(\rho
_{\alpha }\otimes \mu _{\beta })(I_{\alpha }\otimes \tau \otimes I_{\beta
})(\theta _{\alpha ,\beta }\otimes \Delta _{\alpha ,\beta })$%
\end{tabular}%
\end{tabular}%
\end{equation*}
\end{enumerate}

Let C be a $\pi -$coalgebra. A $\pi -$coideal of $C$ is a family $%
V=\{V_{\alpha }\}_{\alpha \in \pi }$ such \nolinebreak that

\begin{enumerate}
\item $V_{\alpha }$ is subspace of $C_{\alpha }$ for all $\alpha \in \pi ,$

\item $\Delta _{\alpha ,\beta }(V_{\alpha \beta })\subset V_{\alpha }\otimes
C_{\beta }+C_{\alpha }\otimes V_{\beta }$ for $\forall\alpha ,\beta \in \pi .
$ \mbox{} and $\epsilon _{1}(V_{1})=0.$
\end{enumerate}

A Hopf $\pi $-coideal of $H$ is a family of $K$-spaces $V=\{V_{\alpha
}\}_{\alpha \in \pi }$ such that $V$ is $\pi $-coideal of $H$, $V_{\alpha }$
is ideal of $H_{\alpha }$ and $S_{\alpha }(V_{\alpha })\subseteq V_{\alpha
^{-1}}$ for all $\alpha \in \pi .$ \newline

A family $A=\{A_{\alpha }\}_{\alpha \in \pi }$ is called subHopf $\pi -$%
coalgebra of $H$ $\ $if $(A,\Delta ,\epsilon )$ is $\pi -$coalgebra of $H$, $%
A_{\alpha }$ is a subalgebra of $H_{\alpha }{}$ and $S_{\alpha }(A_{\alpha
})\subseteq A_{\alpha ^{-1}}$ $\forall \alpha \in \pi .$ A family $%
A=\{A_{\alpha }\}_{\alpha \in \pi }$ of a subalgebra of Hopf $\pi $%
-coalgebra of $H$ is called isolated if there exist a family $I=\{I_{\alpha
}\}_{\alpha \in \pi }$ of Hopf $\pi $- coideal of $H$ such that $H_{\alpha
}=A_{\alpha }\oplus I_{\alpha }$ $\forall \alpha \in \pi .$\newline

We will call Hopf $\pi -$subcoalgebra of $H$ any pair $(C,\sigma )$ such
that $C=\{C_{\alpha }\}_{\alpha \in \pi }$ is a Hopf $\pi -$coalgebra and $%
\sigma =\{\sigma _{\alpha }:H_{\alpha }\rightarrow C_{\alpha }\}_{\alpha \in
\pi }$ family of algebra epimorphisms which satisfies for $\alpha ,\beta \in
\pi $

\begin{enumerate}
\item $\Delta _{\alpha ,\beta }^{C}\sigma _{\alpha \beta }=(\sigma _{\alpha
}\otimes \sigma _{\beta })\Delta _{\alpha ,\beta }^{H},$ i.e., the following
diagram is commute
\begin{equation*}
\begin{tabular}{llll}
& $H_{\alpha \beta }$ & $\overset{\sigma _{\alpha ,\beta }}{\longrightarrow }
$ & $C_{\alpha \beta }$ \\
$\Delta _{\alpha ,\beta }^{H}$ & $\downarrow $ &  & $\downarrow \Delta
_{\alpha ,\beta }^{C}$ \\
& $H_{\alpha }\otimes H_{\beta }$ &
\begin{tabular}{l}
$\longrightarrow $ \\
$\sigma _{\alpha }\otimes \sigma _{\beta }$%
\end{tabular}
& $C_{\alpha }\otimes C_{\beta }$%
\end{tabular}%
\end{equation*}

\item $\varepsilon ^{C} \sigma _{1}=\varepsilon ^{H},$ i.e., the following
diagram is commute
\begin{equation*}
\begin{tabular}{lllll}
& $H_{1}$ & $\overset{\sigma _{1}}{\longrightarrow }$ & $C_{1}$ &  \\
$\epsilon ^{H}$ & $\searrow $ &  & $\swarrow $ & $\epsilon ^{C}$ \\
&  & $K$ &  &
\end{tabular}%
\end{equation*}

\item $S_{\alpha }^{C}\sigma _{\alpha }=\sigma _{\alpha ^{-1}}S_{\alpha
}^{H} $ i.e., the following diagram is commute
\begin{equation*}
\begin{tabular}{lllll}
& $H_{\alpha }$ & $\overset{\sigma _{\alpha }}{\longrightarrow }$ & $%
C_{\alpha }$ &  \\
$S_{\alpha }^{H}$ & $\downarrow $ &  & $\downarrow $ & $S_{\alpha }^{C}$ \\
& $H_{\alpha ^{-1}}$ &
\begin{tabular}{l}
$\longrightarrow $ \\
$\sigma _{\alpha ^{-1}}$%
\end{tabular}
& $C_{\alpha ^{-1}}$ &
\end{tabular}%
\end{equation*}
\end{enumerate}

A pair $(C,\sigma )$ is called left $\pi -$coisotropic quantum subgroup of $%
H $ if

\begin{enumerate}
\item $C$ is $\pi -$coalgebra,

\item $C_{\alpha }$ is left $H_{\alpha }$-module by $\omega _{\alpha }$ for
all $\alpha \in \pi ,$

\item $\sigma =\{\sigma _{\alpha }:H_{\alpha }\rightarrow C_{\alpha }$ $%
\}_{\alpha \in \pi }$ family of surjective linear maps such that

\begin{enumerate}
\item $\sigma _{\alpha }$ is left $H_{\alpha }$-module map for all $\alpha
\in \pi ,$

\item $\Delta _{\alpha ,\beta }^{C}\sigma _{\alpha \beta }=(\sigma _{\alpha
}\otimes \sigma _{\beta })\Delta _{\alpha ,\beta }^{H}$ \mbox{} and \mbox{} $%
\varepsilon ^{C}\circ \sigma _{1}=\varepsilon ^{H}$.\newline
\end{enumerate}
\end{enumerate}

If $g=\{g_{\alpha }:C_{\alpha }\rightarrow H_{\alpha }\}_{\alpha \in \pi },$
is a family of linear maps, its convolution inverse is a family of linear
maps $g^{-1}=\{g_{\alpha }^{-1}:C_{\alpha ^{-1}}\rightarrow H_{\alpha
}\}_{\alpha \in \pi }$ such that
\begin{equation*}
\mu _{\alpha }(g_{\alpha }\otimes g_{\alpha }^{-1})\Delta _{\alpha ,\alpha
^{-1}}^{C}(c)=\epsilon ^{C}(c)1_{H_{\alpha }}=\mu _{\alpha }(g_{\alpha
}^{-1}\otimes g_{\alpha })\Delta _{\alpha ^{-1},\alpha }^{C}(c).\newline
\end{equation*}

\begin{prop}
Every isolated subHopf $\pi -$coalgebra is $\pi -$coisotropic quantum
subgroup.
\end{prop}

\begin{proof} Clear from definition.
\end{proof}
A left $\pi -$coisotropic quantum subgroup $(C,\sigma )$ of $H$ is said to
have a left section if there exist a family of linear, convolution
invertible, maps $g=\{g_{\alpha }:C_{\alpha }\rightarrow H_{\alpha
}\}_{\alpha \in \pi }$ such that

\begin{enumerate}
\item $g_{\alpha }(\sigma _{\alpha }(1))=1,$

\item for $\alpha \in \pi ,u\in H_{1},c\in C_{\alpha }$ and $v\in \sigma
_{\alpha }^{-1}(c)$ we have \newline
$(\sigma _{1}\otimes I_{\alpha }^{H})(\mu _{1}\otimes I_{\alpha
})(I_{1}\otimes \tau )(\Delta _{1,\alpha }^{H}g_{\alpha }\otimes
I_{1})(c\otimes u)\newline
=(I_{1}\otimes g_{\alpha })(\sigma _{1}\otimes \sigma _{\alpha })(\mu
_{1}\otimes I_{\alpha })(I_{1}\otimes \tau )(\Delta _{1,\alpha }^{H}\otimes
I_{1})(v\otimes u)$

\item for $\alpha \in \pi ,u\in H_{1},c\in C_{\alpha ^{-1}}$ and $v\in
\sigma _{\alpha ^{-1}}^{-1}(c)$ we have \newline
$(\sigma _{1}\otimes I_{\alpha }^{H})(\mu _{1}\otimes I_{\alpha
})(I_{1}\otimes \tau )(\Delta _{1,\alpha }^{H}g_{\alpha }^{-1}\otimes
I_{1})(c\otimes u) \newline
=(I_{1}\otimes g_{\alpha }^{-1})(\sigma _{1}\otimes \sigma _{\alpha
^{-1}})(\mu _{1}(S_{1}^{H}\otimes I_{1})\otimes I_{\alpha
^{-1}})(I_{1}\otimes \tau )(\tau \Delta _{\alpha ^{-1},1}^{H}\otimes
I_{1})(v\otimes u)$
\end{enumerate}

\begin{thm}
Let C be $\pi -$coalgebra, $V$ be a $\pi -$coideal and $\sigma =\{\sigma
_{\alpha }:C_{\alpha }\longrightarrow E_{\alpha }=C_{\alpha }/V_{\alpha
}\}_{\alpha \in \pi }$ be a family of the natural linear maps onto the
quotient vector spaces. Then $E=\{E_{\alpha }=C_{\alpha }\diagup V_{\alpha
}\}_{\alpha \in \pi }$ has a unique $\pi -$coalgebra structure such that $%
\sigma $ is a $\pi -$coalgebra maps.\newline
\end{thm}

\begin{proof} Since V is $\pi -$coideal, $\epsilon _{1}(V_{1})=0$ implies $%
V_{1}\subseteq Ker\epsilon _{1}.$ Hence there exist unique linear map $%
\overline{\epsilon }_{1}$ making the following diagram commutes%
\begin{equation*}
\begin{array}{lll}
C_{1} & \overset{\epsilon _{1}}{\longrightarrow } & K \\
\sigma _{1}\searrow &  & \nearrow \overline{\epsilon _{1}} \\
& E_{1} &
\end{array}%
\end{equation*}%
Let
\begin{equation*}
T=\{T_{\alpha ,\beta }:C_{\alpha \beta }\overset{\Delta _{\alpha ,\beta }}{%
\longrightarrow }C_{\alpha }\otimes C_{\beta }\overset{\sigma
_{\alpha }\otimes \sigma _{\beta }}{\longrightarrow }E_{\alpha
}\otimes E_{\beta }\}.
\end{equation*}%
Since $Ker\sigma _{\alpha }=V_{\alpha },$ then $Ker(\sigma
_{\alpha }\otimes \sigma _{\beta })=C_{\alpha }\otimes V_{\beta
}+V_{\alpha }\otimes C_{\beta }\supseteq \Delta _{\alpha ,\beta
}(V_{\alpha \beta }),$ so $V_{\alpha \beta }\subseteq KerT_{\alpha
,\beta }$ and since $T_{\alpha ,\beta }:C_{\alpha \beta
}\longrightarrow E_{\alpha }\otimes E_{\beta }$ then there exist
unique linear maps $\overline{\Delta }_{\alpha ,\beta }:C_{\alpha
\beta }/V_{\alpha \beta }=E_{\alpha \beta }\longrightarrow
E_{\alpha }\otimes
E_{\beta }$ such that the following diagram commutes%
\begin{equation*}
\begin{array}{lllll}
C_{\alpha \beta } &  & \overset{T_{\alpha ,\beta
}}{\longrightarrow } &  &
E_{\alpha }\otimes E_{\beta } \\
& \sigma _{\alpha \beta }\searrow &  & \nearrow \overline{\Delta
}_{\alpha
,\beta } &  \\
&  & E_{\alpha \beta } &  &
\end{array}%
\end{equation*}%
i.e., $\overline{\Delta }_{\alpha ,\beta }\sigma _{\alpha \beta
}=T_{\alpha ,\beta }=(\sigma _{\alpha }\otimes \sigma _{\beta
})\Delta _{\alpha ,\beta }$ for all $\alpha ,\beta \in \pi .$ Thus
$\overline{\Delta }_{\alpha ,\beta \gamma }\sigma _{\alpha \beta
\gamma }=$ $(\sigma _{\alpha }\otimes \sigma _{\beta \gamma
})\Delta _{\alpha ,\beta \gamma }$ for all $\alpha ,\beta ,\gamma
\in \pi .$ Now,
\begin{eqnarray*}
&&(I_{\alpha }\otimes \overline{\Delta }_{\beta ,\gamma })\overline{\Delta }%
_{\alpha ,\beta \gamma }\sigma _{\alpha \beta \gamma } \\
&=&(I_{\alpha }\otimes \overline{\Delta }_{\beta ,\gamma })(\sigma
_{\alpha
}\otimes \sigma _{\beta \gamma })\Delta _{\alpha ,\beta \gamma } \\
&=&(\sigma _{\alpha }\otimes \overline{\Delta }_{\beta ,\gamma
}\sigma
_{\beta \gamma })\Delta _{\alpha ,\beta \gamma } \\
&=&(\sigma _{\alpha }\otimes (\sigma _{\beta }\otimes \sigma
_{\gamma
})\Delta _{\beta ,\gamma })\Delta _{\alpha ,\beta \gamma } \\
&=&(\sigma _{\alpha }\otimes \sigma _{\beta }\otimes \sigma
_{\gamma })(I_{\alpha }\otimes \Delta _{\beta ,\gamma })\Delta
_{\alpha ,\beta \gamma
} \\
&=&(\sigma _{\alpha }\otimes \sigma _{\beta }\otimes \sigma
_{\gamma })(\Delta _{\alpha ,\beta }\otimes I_{\gamma })\Delta
_{\alpha \beta ,\gamma
} \\
&=&((\sigma _{\alpha }\otimes \sigma _{\beta })\Delta _{\alpha
,\beta
}\otimes \sigma _{\gamma })\Delta _{\alpha \beta ,\gamma } \\
&=&(\overline{\Delta }_{\alpha ,\beta }\sigma _{\alpha \beta
}\otimes \sigma
_{\gamma })\Delta _{\alpha \beta ,\gamma } \\
&=&(\overline{\Delta }_{\alpha ,\beta }\otimes I_{\gamma })(\sigma
_{\alpha
\beta }\otimes \sigma _{\gamma })\Delta _{\alpha \beta ,\gamma } \\
&=&(\overline{\Delta }_{\alpha ,\beta }\otimes I_{\gamma })\overline{\Delta }%
_{\alpha \beta ,\gamma }\sigma _{\alpha \beta \gamma }\text{ \ \ \
\ \ \ \ \ \ \ \ \ \ \ \ \ \ \ \ \ \ \ \ \ \ \ \ \ \ \ \ \ \ \ \ \
\ \ \ \ \ \ \ \ \ \ \ \ \ \ \ \ \ \ \ \ \ \ \ \ \ \ \ \ \ \ \ \ \
\ \ \ \ \ \ \ \ \ \ \ \ \ \ \ \ \ \ \ \ \ \ \ \ \ \ \ \ \ \ \ \ \
\ \ \ \ \ \ \ \ \ \ \ \ \ \ \ \ \ \ \ \ \ \ \ \ \ \ \ }
\end{eqnarray*}%
since $\sigma _{\alpha }$ is surjective for all $\alpha \in \pi \,\,\,$then $%
(I_{\alpha }\otimes \overline{\Delta }_{\beta ,\gamma })\overline{\Delta }%
_{\alpha ,\beta \gamma }=(\overline{\Delta }_{\alpha ,\beta
}\otimes
I_{\gamma })\overline{\Delta }_{\alpha \beta ,\gamma }.$ Also, since $%
\overline{\Delta }_{\alpha ,\beta }\sigma _{\alpha \beta }=(\sigma
_{\alpha
}\otimes \sigma _{\beta })\Delta _{\alpha ,\beta }$ implies $\overline{%
\Delta }_{\alpha ,1}\sigma _{\alpha 1}=(\sigma _{\alpha }\otimes
\sigma _{1})\Delta _{\alpha ,1}$ and so
\begin{eqnarray*}%
&&(I_{E_{\alpha }}\otimes \overline{\epsilon
}_{1})\overline{\Delta }_{\alpha
,1}\sigma _{\alpha 1}\\
&=&(I_{E_{\alpha }}\otimes \overline{\epsilon }%
_{1})(\sigma _{\alpha }\otimes \sigma _{1})\Delta _{\alpha
,1}\\
&=&(\sigma _{\alpha }\otimes I_{K})\underbrace{(I_{H_{\alpha
}}\otimes \epsilon _{1})\Delta _{\alpha ,1}}\\
&=&\sigma _{\alpha 1},\text{ \ \ \ \ \ \ \ \ \ \ \ \ \ \ \ \ \ \ \
\ \ \ \ \ \ \ \ \ \ \ \ \ \ \ \ \ \ \ \ \ \ \ \ \ \ \ \ \ \ \ \ \
\ \ \ \ \ \ \ \ \ \ \ \ \ \ \ \ \ \ \ \ \ \ \ \ \ \ \ \ \ \ \ \ \
\ \ \ \ \ \ \ \ \ \ \ \ \ \ \ \ \ \ \ \ \ \ \ \ \ \ \ \ \ \ \ \ \
\ \ \ \ \ \ \ \ \ \ \ }
\end{eqnarray*}%
$i.e.,(I_{E_{\alpha }}\otimes \overline{\epsilon }_{1})\overline{\Delta }%
_{\alpha ,1}\sigma _{\alpha 1}=\sigma _{\alpha 1}.$ Since $\sigma
_{\alpha }$ is surjective for all $\alpha \in \pi ,$ then
$(I_{E_{\alpha }}\otimes \overline{\epsilon }_{1})\overline{\Delta
}_{\alpha ,1}=I_{E_{\alpha }}.$
Similar, $(\overline{\epsilon }_{1}\otimes I_{E_{\alpha }})\overline{\Delta }%
_{1,\alpha }=I_{E_{\alpha }}.$ Thus $(E,\overline{\Delta },\overline{%
\epsilon }_{1})$ is $\pi -$coalgebra. \newline
\end{proof}

\section{{\protect\large {Induced representations of Hopf group coalgebra}}}

In this section, we study the induced representation for Hopf group
coalgebra. To reach this goal we use the definitions of subHopf group
coalgebra, Hopf group subcoalgebra and group coisotropic quantum subgroup.

\begin{thm}
Let $H$ be a Hopf $\pi -$coalgebra and $(C,\sigma )$ be a Hopf $\pi-$%
subcoalgebra of $H.$ Then $(C,\sigma )$ is a left $\pi -$coisotropic quantum
subgroup of $H.$
\end{thm}

\begin{proof}
We define $\omega =\{\omega _{\alpha }=\mu _{\alpha }^{C}(\sigma
_{\alpha }\otimes I_{\alpha }):H_{\alpha }\otimes C_{\alpha
}\rightarrow C_{\alpha }\}_{\alpha \in \pi }$ we'll prove that
$C_{\alpha }$ is a left $H_{\alpha }-$module by $\omega _{\alpha
}.$ i.e., the following diagrams are commute
\begin{equation*}
\begin{tabular}{llll}
& $H_{\alpha }\otimes H_{\alpha }\otimes C_{\alpha }$ & $%
\begin{tabular}{l}
$I_{\alpha }\otimes \omega _{\alpha }$ \\
$\ \ \rightarrow $%
\end{tabular}%
$ & $H_{\alpha }\otimes C_{\alpha }$ \\
$\mu _{\alpha }\otimes I_{\alpha }$ & $\downarrow $ &  &
$\downarrow \omega
_{\alpha }$ \\
& $H_{\alpha }\otimes C_{\alpha }$ &
\begin{tabular}{l}
$\longrightarrow $ \\
$\omega _{\alpha }$%
\end{tabular}
& $C_{\alpha }$%
\end{tabular}%
\mbox{}\hspace{.5cm}\mbox{}\\
\begin{tabular}{lll}
$K\otimes C_{\alpha }$ & $%
\begin{tabular}{l}
$\eta _{\alpha }\otimes I_{\alpha }$ \\
$\ \ \rightarrow $%
\end{tabular}%
$ & $H_{\alpha }\otimes C_{\alpha }$ \\
$\sim \searrow $ &  & $\swarrow \omega _{\alpha }$ \\
& $C_{\alpha }$ &
\end{tabular}%
\end{equation*}
\begin{eqnarray*}
&& \omega _{\alpha }(I_{\alpha }\otimes \omega _{\alpha
})(h\otimes k\otimes b)\\
&=&\omega _{\alpha }(h\otimes \sigma _{\alpha }(k)b) \\
&=&\sigma _{\alpha }(h)(\sigma _{\alpha }(k)b) \\
&=&\sigma _{\alpha }(hk)b\\
&=&\omega _{\alpha }(hk\otimes b)\\
&=&\omega _{\alpha }(\mu _{\alpha }\otimes I_{\alpha })(h\otimes
k\otimes b) \\
&=&\omega _{\alpha }(\mu _{\alpha }\otimes I_{\alpha })(h\otimes
k\otimes b)
\end{eqnarray*}
and%
\begin{eqnarray*}
\omega _{\alpha }(\eta _{\alpha }\otimes I_{\alpha })(k\otimes b)
&=&\omega _{\alpha }(\eta _{\alpha }(k)\otimes b)=\sigma _{\alpha
}(\eta _{\alpha }(k))b = kb =\sim (k\otimes b)
\end{eqnarray*}
Now, we`ll prove that $\sigma _{\alpha }$ is left $H_{\alpha
}$-module map $\forall \alpha \in \pi $ i.e., the following
diagram is commute %
\begin{equation*}
\begin{tabular}{llll}
& $H_{\alpha }$ & $%
\begin{tabular}{l}
$\sigma _{\alpha }$ \\
$\rightarrow $%
\end{tabular}%
$ & $C_{\alpha }$ \\
$\mu _{\alpha }^{H}$ & $\uparrow $ &  & $\uparrow \omega _{\alpha }$ \\
& $H_{\alpha }\otimes H_{\alpha }$ &
\begin{tabular}{l}
$\longrightarrow $ \\
$I_{\alpha }\otimes \sigma _{\alpha }$%
\end{tabular}
& $H_{\alpha }\otimes C_{\alpha }$%
\end{tabular}%
\end{equation*}
\begin{equation*}
\omega _{\alpha }(I_{\alpha }\otimes \sigma _{\alpha })=\mu
_{\alpha }^{C}(\sigma _{\alpha }\otimes I_{\alpha })(I_{\alpha
}\otimes \sigma _{\alpha })=\mu _{\alpha }^{C}(\sigma _{\alpha
}\otimes \sigma _{\alpha })=\sigma _{\alpha }\mu _{\alpha }^{H}.
\end{equation*}%
\end{proof}

\begin{rem}
If $(C,\sigma )$ is a Hopf $\pi -$subcoalgebra of $H$, then the map $%
L_{\alpha ,\beta }=(\sigma _{\alpha }\otimes I_{\beta }^{H})\Delta _{\alpha
,\beta }^{H}:H_{\alpha \beta }\rightarrow C_{\alpha }\otimes H_{\beta
}~~\forall \alpha ,\beta \in \pi $ is an algebra map.
\end{rem}

\begin{prop}
Let $H=\{H_{\alpha }\}_{\alpha \in \pi }$ be a Hopf $\pi -$coalgebra and $%
(C,\sigma )$ be a Hopf $\pi -$subcoalgebra of $H.$ Let $B=\{B_{\alpha
}\}_{\alpha \in \pi },$ where $B_{\alpha }=\{h\in H_{\alpha }:L_{1,\alpha
}(h)=1\otimes h\}, $ then $B_{\alpha }$ is subalgebra of $H_{\alpha }$ and $B
$ is right $\pi -$coideal of $H$ (this B is called $\pi -$quantum right
embeddable homogeneous space of $H$ ).
\end{prop}

\begin{proof}\mbox{}
\;\;\;\; Firstly, let $h,k\in B_{\alpha },$ then $L_{1,\alpha
}(hk)=L_{1,\alpha }(h)\cdot L_{1,\alpha }(k)=(1\otimes h)\cdot
(1\otimes k)=1\otimes hk$ and hence $hk\in B_{\alpha }.$ Secondly
let $h\in B_{\alpha \beta },$ we will prove that $\Delta _{\alpha
,\beta }^{H}(h)\in B_{\alpha }\otimes H_{\beta }.$
\begin{eqnarray*}
&&(L_{1,\alpha }\otimes I_{\beta }^{H})\Delta _{\alpha ,\beta
}^{H}(h)\\
&=&((\sigma _{1}\otimes I_{\alpha }^{H})\Delta _{1,\alpha
}^{H}\otimes I_{\beta }^{H})\Delta _{\alpha ,\beta
}^{H}(h)\\
&=&(\sigma _{1}\otimes I_{\alpha }^{H}\otimes I_{\beta
}^{H})(\Delta _{1,\alpha }^{H}\otimes I_{\beta }^{H})\Delta
_{\alpha ,\beta }^{H}(h) \\
&=&(\sigma _{1}\otimes I_{\alpha }^{H}\otimes I_{\beta
}^{H})(I_{1}^{H}\otimes \Delta _{\alpha
,\beta }^{H})\Delta _{1,\alpha \beta }^{H}(h)\\
&=&(I_{1}^{C}\otimes \Delta _{\alpha ,\beta }^{H})(\sigma
_{1}\otimes I_{\alpha \beta
}^{H})\Delta _{1,\alpha \beta }^{H}(h)\\
&=&(I_{1}^{C}\otimes \Delta%
_{\alpha ,\beta }^{H})(1\otimes h)\\
&=&1\otimes \Delta _{\alpha ,\beta%
}^{H}(h).\text{ \ \ \ \ \ \ \ \ \ \ \ \ \ \ \ \ \ \ \ \ \ \ \ \ \
\ \ \ \ \ \ \ \ \ \ \ \ \ \ \ \ \ \ \ \ \ \ \ \ \ \ \ \ \ \ \ \ \
\ \ \ \ \ \ \ \ \ \ \ \ \ \ \ \ \ \ \ \ \ \ \ \ \ \ \ \ \ \ \ \ \
\ \ \ \ \ \ \ \ \ \ \ \ \ \ \ \ \ \ \ \ \ \ \ \ \ \ \ \ \ \ \ \ \
\ \ \ \ \ }
\end{eqnarray*}
\end{proof}

\begin{lem}
Let $H=\{H_{\alpha }\}_{\alpha \in \pi }$ be a Hopf $\pi - $coalgebra and $%
(C,\sigma )$ be a left $\pi -$coisotropic quantum subgroup of $H.$ For $%
\alpha ,\beta ,\gamma \in \pi $, $a,b\in H_{\alpha \beta }$ we have

\begin{itemize}
\item $L_{\alpha ,\beta }(ab)=\Delta _{\alpha ,\beta }^{H}(a)\Theta
L_{\alpha ,\beta }(b)$

\item $(I_{\alpha }\otimes \Delta _{\beta ,\gamma }^{H})L_{\alpha ,\beta
\gamma }=(L_{\alpha ,\beta }\otimes I_{\gamma })\Delta _{\alpha \beta
,\gamma }^{H}$,
\end{itemize}

where $(m\otimes n)\Theta (u\otimes v)=\omega _{\alpha }(m\otimes u)\otimes
nv,m\in H_{\alpha },u\in C_{\alpha },n,v\in H_{\beta }.$
\end{lem}

\begin{proof}
\begin{eqnarray*}
&& L_{\alpha ,\beta }(ab) \\
&=&(\sigma _{\alpha }\otimes I_{\beta }^{H})\Delta _{\alpha ,\beta
}^{H}\mu _{\alpha \beta }(a\otimes b)\\
&=& (\sigma _{\alpha }\otimes I_{\beta }^{H})(\mu _{\alpha
}\otimes \mu _{\beta })(I\otimes \tau \otimes I)(\Delta _{\alpha
,\beta }^{H}\otimes
\Delta _{\alpha ,\beta }^{H})(a\otimes b)\\
&=&(\sigma _{\alpha }\mu _{\alpha }\otimes \mu _{\beta })(I\otimes
\tau \otimes I)(\Delta _{\alpha ,\beta }^{H}\otimes \Delta
_{\alpha ,\beta }^{H})(a\otimes b)\\
&=&\omega _{\alpha }(a_{1}^{\alpha }\otimes \sigma _{\alpha
}(b_{1}^{\alpha }))\otimes a_{2}^{\beta
}b_{2}^{\beta }\\
&=&\Delta _{\alpha ,\beta }^{H}(a)\Theta L_{\alpha ,\beta
}(b)\text{ \ \ \ \ \ \ \ \ \ \ \ \ \ \ \ \ \ \ \ \ \ \ \ \ \ \ \ \
\ \ \ \ \ \ \ \ \ \ \ \ \ \ \ \ \ \ \ \ \ \ \ \ \ \ \ \ \ \ \ \ \
\ \ \ \ \ \ \ \ \ \ \ \ \ \ \ \ \ \ \ \ \ \ \ \ \ \ \ \ \ \ \ \ \
\ \ \ \ \ \ \ \ \ \ \ \ \ \ \ \ \ \ \ \ \ \ \ \ \ \ \ \ \ \ \ \ \
\ \ }
\end{eqnarray*}%
Also,
\begin{eqnarray*}
&&(I_{\alpha }\otimes \Delta _{\beta ,\gamma }^{H})L_{\alpha
,\beta \gamma } \\
&=&(I_{\alpha }\otimes \Delta _{\beta ,\gamma }^{H})(\sigma
_{\alpha }\otimes I_{\beta \gamma }^{H})\Delta _{\alpha ,\beta
\gamma }^{H} \\
&=&(\sigma _{\alpha }\otimes I_{\beta }^{H}\otimes I_{\gamma
}^{H})(I_{\alpha }\otimes \Delta _{\beta ,\gamma }^{H})\Delta
_{\alpha
,\beta \gamma }^{H} \\
&=&(\sigma _{\alpha }\otimes I_{\beta }^{H}\otimes I_{\gamma
}^{H})(\Delta _{\alpha ,\beta }^{H}\otimes I_{\gamma })\Delta
_{\alpha \beta ,\gamma }^{H} \\
&=&((\sigma _{\alpha }\otimes I_{\beta }^{H})\Delta _{\alpha
,\beta }^{H}\otimes I_{\gamma }^{H})\Delta _{\alpha \beta ,\gamma
}^{H}\\
&=&(L_{\alpha ,\beta }\otimes I_{\gamma })\Delta _{\alpha \beta
,\gamma }^{H}\text{ \ \ \ \ \ \ \ \ \ \ \ \ \ \ \ \ \ \ \ \ \ \ \
\ \ \ \ \ \ \ \ \ \ \ \ \ \ \ \ \ \ \ \ \ \ \ \ \ \ \ \ \ \ \ \ \
\ \ \ \ \ \ \ \ \ \ \ \ \ \ \ \ \ \ \ \ \ \ \ \ \ \ \ \ \ \ \ \ \
\ \ \ \ \ \ \ \ \ \ \ \ \ \ \ \ \ \ \ \ \ \ \ \ \ \ \ \ \ \ \ \ \
\ \ \ \ \ \ \ }
\end{eqnarray*}
\end{proof}

\begin{prop}
Let $H=\{H_{\alpha }\}_{\alpha \in \pi }$ be a Hopf $\pi -$coalgebra and let
$(C,\sigma )$ be a left $\pi -$coisotropic quantum subgroup of $H.$ Let $%
G=\{G_{\alpha }\}_{\alpha \in \pi },$ where $G_{\alpha }=\{h\in H_{\alpha
}:L_{1,\alpha }(h)=\sigma _{1}(1)\otimes h\}$, then $G_{\alpha }$ is
subalgebra of $H_{\alpha }$ and $G$ is right $\pi -$ coideal of $H$.
\end{prop}

\begin{proof}
For $h,g\in G_{\alpha },$ we have $L_{1,\alpha }(h)=\sigma
_{1}(1)\otimes h$, $L_{1,\alpha }(g)=\sigma _{1}(1)\otimes g$ and
$\omega _{1}(I_{1}\otimes \sigma _{1})=\sigma _{1}\mu _{1}$.
\\
We`ll prove that $hg\in G_{\alpha },i.e.,L_{1,\alpha }(hg)=\sigma
_{1}(1)\otimes hg$ \\
\begin{eqnarray*}
&& L_{1,\alpha }(hg) \\
&=& \Delta_{1,\alpha}(h)\Theta
L_{1,\alpha}(g)\\
&=&\Delta_{1,\alpha}(h)\Theta (\sigma_{1}(1)\otimes g)\\
&=&(\omega _{1}\otimes \mu _{\alpha })(h_{1}^{1}\otimes \sigma
_{1}(1)\otimes h_{2}^{\alpha }\otimes g) \\
&=&\omega _{1}(I_{1}\otimes \sigma _{1})(h_{1}^{1}\otimes
1)\otimes h_{2}^{\alpha }g \\
&=&\sigma _{1}\mu _{1}(h_{1}^{1}\otimes 1)\otimes h_{2}^{\alpha
}g\\
&=&\sigma _{1}(h_{1}^{1})\otimes
h_{2}^{\alpha }g \\
&=&(\sigma _{1}\otimes \mu _{\alpha })(\Delta _{1,\alpha
}^{H}\otimes I_{\alpha })(h\otimes g)\\
&=&(I_{1}\otimes \mu _{\alpha })((\sigma _{1}\otimes I_{1})\Delta
_{1,\alpha }^{H}(h)\otimes g)
\\
&=&(I_{1}\otimes \mu _{\alpha })(\sigma _{1}(1)\otimes h\otimes
g)\\
&=&\sigma _{1}(1)\otimes hg \text{ \ \ \ \ \ \ \ \ \ \ \ \ \ \ \ \
\ \ \ \ \ \ \ \ \ \ \ \ \ \ \ \ \ \ \ \ \ \ \ \ \ \ \ \ \ \ \ \ \
\ \ \ \ \ \ \ \ \ \ \ \ \ \ \ \ \ \ \ \ \ \ \ \ \ \ \ \ \ \ \ \ \
\ \ \ \ \ \ \ \ \ \ \ \ \ \ \ \ \ \ \ \ \ \ \ \ \ \ \ \ \ \ \ \ \
\ \ \ \ \ \ \ \ \ \ \ \ \ \ }
\end{eqnarray*}%
We`ll prove that $G$ is right $\pi $-coideal of $H$, i.e., for
$h\in G_{\alpha \beta },\Delta _{\alpha ,\beta }^{H}(h)\in
G_{\alpha }\otimes H_{\beta }$
\begin{eqnarray*}
&&(L_{1,\alpha }\otimes I_{\beta }^{H})\Delta _{\alpha ,\beta
}^{H}(h) \\
&=&((\sigma _{1}\otimes I_{\alpha }^{H})\Delta _{1,\alpha
}^{H}\otimes I_{\beta }^{H})\Delta _{\alpha ,\beta
}^{H}(h)\\
&=&(\sigma _{1}\otimes I_{\alpha }^{H}\otimes I_{\beta
}^{H})(\Delta
_{1,\alpha }^{H}\otimes I_{\beta }^{H})\Delta _{\alpha ,\beta }^{H}(h) \\
&=&(\sigma _{1}\otimes I_{\alpha }^{H}\otimes I_{\beta
}^{H})(I_{1}^{H}\otimes \Delta _{\alpha ,\beta }^{H})\Delta
_{1,\alpha \beta }^{H}(h)\\
&=&(I_{1}^{C}\otimes \Delta _{\alpha ,\beta }^{H})(\sigma
_{1}\otimes
I_{\alpha \beta }^{H})\Delta _{1,\alpha \beta }^{H}(h) \\
&=&(I_{1}^{C}\otimes \Delta _{\alpha ,\beta }^{H})(\sigma
_{1}(1)\otimes h)\\
&=&\sigma _{1}(1)\otimes \Delta _{1,\alpha \beta }^{H}(h)\text{ \
\ \ \ \ \ \ \ \ \ \ \ \ \ \ \ \ \ \ \ \ \ \ \ \ \ \ \ \ \ \ \ \ \
\ \ \ \ \ \ \ \ \ \ \ \ \ \ \ \ \ \ \ \ \ \ \ \ \ \ \ \ \ \ \ \ \
\ \ \ \ \ \ \ \ \ \ \ \ \ \ \ \ \ \ \ \ \ \ \ \ \ \ \ \ \ \ \ \ \
\ \ \ \ \ \ \ \ \ \ \ \ \ \ \ \ \ \ \ \ \ \ \ \ \ \ \ \ \ }
\end{eqnarray*}%
\end{proof}

\begin{thm}
Suppose $V=\{V_{\alpha }\}_{\alpha \in \pi },$ be a right $\pi -$comodule
over the left $\pi -$coisotropic quantum subgroup $C$ of Hopf $\pi -$%
coalgebra $H$ by $\rho =\{\rho _{\alpha ,\beta }:V_{\alpha \beta
}\rightarrow V_{\alpha }\otimes C_{\beta }\}_{\alpha ,\beta \in \pi }.$ Then
$Ind(\rho )=\{Ind(\rho )_{\alpha }\}_{\alpha \in \pi }$ where $Ind(\rho
)_{\alpha }=\{x\in V_{1}\otimes H_{\alpha }:(I_{1}\otimes L_{1,\alpha
})x=(\rho _{1,1}\otimes I_{\alpha })x\} $ is right $\pi -$comodule over $H$
by $(I\otimes \Delta )=\{(I\otimes \Delta )_{\alpha ,\beta }=I_{1}\otimes
\Delta _{\alpha ,\beta }\}_{\alpha ,\beta \in \pi }$
\end{thm}

\begin{proof}
We will prove that, for $\alpha ,\beta  \in \pi ,$ $(I\otimes
\Delta )_{\alpha ,\beta }(Ind(\rho )_{\alpha \beta })\subseteq
Ind(\rho )_{\alpha }\otimes H_{\beta }.$ Since, for $v\otimes h\in
Ind(\rho )_{\alpha \beta },$ we have
\begin{eqnarray*}
&&(I_{1}\otimes L_{1,\alpha }\otimes I_{\beta })(I_{1}\otimes
\Delta _{\alpha ,\beta })(v\otimes h) \\
&=&(I_{1}\otimes (L_{1,\alpha }\otimes I_{\beta })\Delta _{\alpha
,\beta })(v\otimes h)\\
&=&(I_{1}\otimes (I_{1}\otimes \Delta _{\alpha ,\beta
})L_{1,\alpha \beta
})(V\otimes h)\;\;\;\;\;\;\;\;\;\;\;\;\;\;\;\;\;\;\;\;\;\;\; \textrm{by Lemma 2.4  }\\
&=&(I_{1}\otimes I_{1}\otimes \Delta _{\alpha ,\beta
})(I_{1}\otimes L_{1,\alpha \beta })(v\otimes h) \\
&=&(I_{1}\otimes I_{1}\otimes \Delta _{\alpha ,\beta })(\rho
_{1,1}\otimes
I_{\alpha \beta })(v\otimes h) \\
&=&(\rho _{1,1}\otimes I_{\alpha }\otimes I_{\beta })(I_{1}\otimes
\Delta _{\alpha ,\beta })(v\otimes h)\text{ \ \ \ \ \ \ \ \ \ \ \
\ \ \ \ \ \ \ \ \ \ \ \ \ \ \ \ \ \ \ \ \ \ \ \ \ \ \ \ \ \ \ \ \
\ \ \ \ \ \ \ \ \ \ \ \ \ \ \ \ \ \ \ \ \ \ \ \ \ \ \ \ \ \ \ \ \
\ \ \ \ \ \ \ \ \ \ \ \ \ \ \ \ \ \ \ \ \ \ \ \ \ \ \ \ \ \ \ \ \
\ \ \ \ \ \ \ \ \ \ \ \ \ \ \ \ \ \ \ }
\end{eqnarray*}

Now, we`ll prove that the following diagrams are commute
\begin{equation*}
\begin{tabular}{llll}
& $Ind(\rho )_{\alpha \beta \gamma }$ & $\overset{(I\otimes \Delta
)_{\alpha \beta ,\gamma }}{\longrightarrow }$ & $Ind(\rho
)_{\alpha \beta }\otimes
H_{\gamma }$ \\
$(I\otimes \Delta )_{\alpha ,\beta \gamma }$ & $\downarrow $ &  & $%
\downarrow (I\otimes \Delta )_{\alpha ,\beta }\otimes I_{\gamma }$ \\
& $Ind(\rho )_{\alpha }\otimes H_{\beta \gamma }$ &
\begin{tabular}{l}
$\longrightarrow $ \\
$I_{Ind(\rho )_{\alpha }}\otimes \Delta _{\beta ,\gamma }$%
\end{tabular}
& $Ind(\rho )_{\alpha }\otimes H_{\beta }\otimes H_{\gamma }$%
\end{tabular}%
\end{equation*}%
\mbox{and}\\
\begin{equation*}
\begin{tabular}{llll}
$Ind(\rho )_{\alpha }$ &  & $%
\begin{tabular}{l}
$(I\otimes \Delta )_{\alpha ,1}$ \\
$\rightarrow $%
\end{tabular}%
$ & $Ind(\rho )_{\alpha }\otimes H_{1}$ \\
& $\sim _{Ind(\rho )_{\alpha }}\searrow $ &  & $\swarrow
I_{Ind(\rho
)_{\alpha }}\otimes \epsilon ^{H}$ \\
&  & $Ind(\rho )_{\alpha }\otimes K$ &
\end{tabular}%
\end{equation*}
\begin{eqnarray*}
&&((I\otimes \Delta )_{\alpha ,\beta }\otimes I_{\gamma
})(I\otimes \Delta )_{\alpha \beta ,\gamma } \\
&=&(I_{1}\otimes \Delta _{\alpha ,\beta }\otimes I_{\gamma
})(I_{1}\otimes \Delta _{\alpha \beta ,\gamma })\\
&=&(I_{1}\otimes (\Delta _{\alpha ,\beta }\otimes
I_{\gamma })\Delta _{\alpha \beta ,\gamma }) \\
&=&(I_{1}\otimes (I_{\alpha }\otimes \Delta _{\beta ,\gamma
})\Delta _{\alpha ,\beta \gamma }) \\
&=&(I_{1}\otimes I_\alpha\otimes\Delta_{\beta,\gamma})(I_1\otimes
\Delta_{\alpha,\beta\gamma})\\
&=& (I_{Ind(\rho )_{\alpha }}\otimes
\Delta_{\beta,\gamma})(I\otimes\Delta)_{\alpha,\beta\gamma}\text{
\ \ \ \ \ \ \ \ \ \ \ \ \ \ \ \ \ \ \ \ \ \ \ \ \ \ \ \ \ \ \ \ \
\ \ \ \ \ \ \ \ \ \ \ \ \ \ \ \ \ \ \ \ \ \ \ \ \ \ \ \ \ \ \ \ \
\ \ \ \ \ \ \ \ \ \ \ \ \ \ \ \ \ \ \ \ \ \ \ \ \ \ \ \ \ \ \ \ \
\ \ \ \ \ \ \ \ \ \ \ \ \ \ \ \ \ \ \ \ \ \ \ \ \ \ \ \ \ \ }
\end{eqnarray*}
and
\begin{eqnarray*}
&&(I_{Ind(\rho )_{\alpha }}\otimes \epsilon ^{H})(I\otimes \Delta
)_{\alpha ,1}\\
&=&(I_{1}\otimes I_{\alpha }\otimes \epsilon
^{H})(I_{1}\otimes \Delta _{\alpha ,1}) \\
&=&(I_{1}\otimes (I_{\alpha }\otimes \epsilon )\Delta _{\alpha
,1})\\
&=&(I_{1}\otimes \sim _{H_{\alpha }})~\\
&=&\sim _{Ind(\rho )_{\alpha }}\text{ \ \ \ \ \ \ \ \ \ \ \ \ \ \
\ \ \ \ \ \ \ \ \ \ \ \ \ \ \ \ \ \ \ \ \ \ \ \ \ \ \ \ \ \ \ \ \
\ \ \ \ \ \ \ \ \ \ \ \ \ \ \ \ \ \ \ \ \ \ \ \ \ \ \ \ \ \ \ \ \
\ \ \ \ \ \ \ \ \ \ \ \ \ \ \ \ \ \ \ \ \ \ \ \ \ \ \ \ \ \ \ \ \
\ \ \ \ \ \ \ \ \ \ \ \ \ \ \ \ }
\end{eqnarray*}
\end{proof}

\begin{rem}
Given a right corepresentation $\rho $ of the $\pi -$coisotropic quantum
subgroup of $(C,\sigma )$ the corresponding corepresentation $(I\otimes
\Delta)$ on $Ind(\rho )$ of $H$ is called induced representation from $\rho$
on $H$.
\end{rem}

\section{{\protect\large {\ Geometric realization for induced representation}%
}}

In this section we given the geometric interperation of the induced
representation. Throughout this section $H$ is Hopf $\pi -$ coalgebra, $%
(C,\sigma )$ is Hopf $\pi -$subcoalgebra of $H$ and $V=\{V_{\alpha
}\}_{\alpha \in \pi }$ be a right $\pi -$comodule over $C$. The purpose of
this section is to show that the induced representation $Ind(\rho )$ from
Hopf group subcoalgebra $H$ is isomorphic as module to the tensor product of
$\pi -$quantum embeddable homogeneous space $B$ (in Proposition 2.3) and the
given comodule $V$. In case $(C,\sigma )$ is left $\pi -$coisotropic quantum
subgroup, then $H$ is isomorphic to $C\otimes G$ as a vector space where $%
G=\{G_{\alpha }\}_{\alpha \in \pi },$ where
\begin{equation*}
G_{\alpha }=\{h\in H_{\alpha }:L_{1,\alpha }(h)=(\sigma _{1}\otimes
I_{\alpha }^{H})\Delta _{1,\alpha }^{H}(h)=\sigma _{1}(1)\otimes h\}.
\end{equation*}

\begin{lem}
$Ind(\rho )_{\alpha }$ is a right (left) $B_{\alpha }$-module for
all $\alpha \in \pi .$
\end{lem}

\begin{proof}
Let $v_{1}\otimes h_\alpha \in Ind(\rho)_\alpha, b_{\alpha} \in
B_{\alpha},$ we define the right action as follows $\lambda
_{\alpha }(v_{1}\otimes h_{\alpha }\otimes b_{\alpha
})=v_{1}\otimes h_{\alpha }b_{\alpha }.$ We need only to prove
that $\lambda _{\alpha }(v_{1}\otimes h_{\alpha }\otimes b_{\alpha
})=v_{1}\otimes b_{\alpha }h_{\alpha }\in Ind(\rho )_{\alpha }.$
We have $ (I_{1}\otimes L_{1,\alpha })(v_{1}\otimes h_{\alpha
})=(\rho _{1,1}\otimes I\smallskip
_{\alpha })(v_{1}\otimes h_{\alpha })$ imply that%
$((I_{1}\otimes L_{1,\alpha })(v_{1}\otimes h_{\alpha }))(1\otimes
1\otimes b_{\alpha }) =((\rho _{1,1}\otimes I_{\alpha
})(v_{1}\otimes h_{\alpha }))(1\otimes 1\otimes b_{\alpha })$
imply that%
\begin{eqnarray*}
v_{1}\otimes \sigma _{1}((h_{\alpha })_{1})\otimes (h_{\alpha
})_{2}b_{\alpha } &=&\rho _{1,1}(v_{1})\otimes h_{\alpha
}b_{\alpha }.
\;\;\;\;\;\;\;\;\;\;\;\;\;\;\;\;\;\;\;\;\;\;\;\;\;\;\;\;\;\;\;\;\;\;\;\;\;\;\;\;\;\;\;\;
(*)
\end{eqnarray*}
Also, we have
\begin{eqnarray*}
&& L_{1,\alpha }(h_{\alpha }b_{\alpha })\\
&=&\Delta _{1,\alpha
}^{H}(h_{\alpha })\Theta L_{1,\alpha }(b_{\alpha })\\
&=&((h_{\alpha })_{1}\otimes (h_{\alpha })_{2})\Theta
(1_{C_{1}}\otimes b_{\alpha })
\;\;\;\;\;\;\;\;\;\;\;\;\;\;\;\;\;\;\;\;\;\;\;\;\;\;\;\;\;\;\;\;\;\;\;\
\;\;\;\;\;\;\;\;\;\;\;\;\;\;\;
\text{by lemma 2.5}\\
&=&\omega _{\alpha }((h_{\alpha })_{1}\otimes \sigma
_{1}(1_{H_{1}}))\otimes (h_{\alpha })_{2}b_{\alpha }\\
&=&\sigma _{1}(h_{\alpha _{1}})\otimes (h_{\alpha })_{2}b_{\alpha
}.
\;\;\;\;\;\;\;\;\;\;\;\;\;\;\;\;\;\;\;\;\;\;\;\;\;\;\;\;\;\;\;\;\;\;\;\;\;\;\;\;\;\;\;\;
\;\;\;\;\;\;\;\;\;\;\;\;\;\;\;\;\;\;\;\;\;\;\;\;\;\;\;\;\;\;\;\;\;\;\;\;
(**)
\end{eqnarray*}%
Now,
\begin{eqnarray*}
&& (I_{1}\otimes L_{1,\alpha })\lambda _{\alpha }(v_{1}\otimes
h_{\alpha }\otimes b_{\alpha })\\
&=&(I_{1}\otimes L_{1,\alpha })(v_{1}\otimes h_{\alpha }b_{\alpha })\\
&=& v_{1}\otimes L_{1,\alpha }(h_{\alpha }b_{\alpha }) \\
&=&v_{1}\otimes \sigma _{1}((h_{\alpha })_{1})\otimes (h_{\alpha
})_{2}b_{\alpha } \;\;\;\;\;\;\;\;\;\;\;\;\;\;\;\;\;\;\;\;\;\;\;\;\;\;\;\;\;\;\;\;\;\;\;\;\;\;\;\;\;\;\;\;\;\;\;\;\;\;\;\;\;\;\;\;\;\;\;\;\;\;\;\;\;\;\; \text{by (**) }\\
&=&\rho _{1,1}(v_{1})\otimes h_{\alpha }b_{\alpha }\;\;\;\;\;\;\;\;\;\;\;\;\;\;\;\;\;\;\;\;\;\;\;\;\;\;\;\;\;\;\;\;\;\;\;\;\;\;\;\;\;\;\;\;\;\;\;\;\;\;\;\;\;\;\;\;\;\;\;\;\;\;\;\;\;\;\;\;\;\;\;\;\;\;\;\;\;\;\;\;\;\;\text{by (*)} \\
&=&(\rho _{1,1}\otimes I_{\alpha })(v_{1}\otimes h_{\alpha
}b_{\alpha })\\
&=& (\rho _{1,1}\otimes I_{\alpha }) \lambda _{\alpha
}(v_{1}\otimes h_{\alpha }\otimes b_{\alpha }).
\end{eqnarray*}
The left action is similar.
\end{proof}

\begin{defn}
Let $(C,\sigma )$ be a Hopf $\pi -$subcoalgebra of $H.$ If $g=\{g_{\alpha
}:C_{\alpha }\rightarrow H_{\alpha }\}_{\alpha \in \pi },$ is a family of
linear maps, its convolution inverse, if it exists, is a family of linear
maps $g^{-1}=\{g_{\alpha }^{-1}:C_{\alpha ^{-1}}\rightarrow H_{\alpha
}\}_{\alpha \in \pi }$ such that
\begin{equation*}
\mu _{\alpha }(g_{\alpha }\otimes g_{\alpha }^{-1})\Delta _{\alpha ,\alpha
^{-1}}^{C}(c)=\epsilon ^{C}(c)1_{H_{\alpha }}=\mu _{\alpha }(g_{\alpha
}^{-1}\otimes g_{\alpha })\Delta _{\alpha ^{-1},\alpha }^{C}(c)
\end{equation*}
\end{defn}

\begin{defn}
A Hopf $\pi -$subcoalgebra $(C,\sigma )$ of $H$ is said to have a left
section if there exists a family of linear, convolution invertible, maps $%
g=\{g_{\alpha }:C_{\alpha }\rightarrow H_{\alpha }\}_{\alpha \in \pi }$ such
that for all $\alpha\in \pi$, $g_{\alpha }(1)=1$ and $L_{1,\alpha }g_{\alpha
}=(I_{1}\otimes g_{\alpha })\Delta _{1,\alpha }^{C}$
\end{defn}

\begin{lem}
For any Hopf $\pi -$coalgebra $H$ we have

\begin{enumerate}
\item $(\Delta _{1,\alpha }^{H}\otimes \Delta _{\alpha ^{-1},1}^{H})\Delta
_{\alpha ,\alpha ^{-1}}^{H}=(I_{1}\otimes \Delta _{\alpha ,\alpha
^{-1}}^{H}\otimes I_{1})(\Delta _{1,1}^{H}\otimes I_{1})\Delta _{1,1}^{H}$

\item $(\Delta _{\alpha ^{-1},1}^{H}\otimes \Delta _{1,\alpha }^{H})\Delta
_{\alpha ^{-1},\alpha }^{H}=(I_{\alpha ^{-1}}\otimes \Delta
_{1,1}^{H}\otimes I_{\alpha })$($\Delta _{\alpha ^{-1},1}^{H}\otimes
I_{\alpha })\Delta _{\alpha ^{-1},\alpha }^{H} $
\end{enumerate}
\end{lem}

\begin{proof}
\begin{eqnarray*}
&& (\Delta _{1,\alpha }^{H}\otimes \Delta _{\alpha
^{-1},1}^{H})\Delta _{\alpha ,\alpha ^{-1}}^{H}\\
&=&(\Delta _{1,\alpha }^{H}\otimes I_{\alpha ^{-1}}\otimes
I_{1})(I_{\alpha
}\otimes \Delta _{\alpha ^{-1},1}^{H})\Delta _{\alpha ,\alpha ^{-1}}^{H}%
\\
&=&(\Delta _{1,\alpha }^{H}\otimes I_{\alpha ^{-1}}\otimes
I_{1})(\Delta _{\alpha ,\alpha ^{-1}}^{H}\otimes I_{1})\Delta
_{1,1}^{H} \\
&=&((\Delta _{1,\alpha }^{H}\otimes I_{\alpha ^{-1}})\Delta
_{\alpha ,\alpha
^{-1}}^{H}\otimes I_{1})\Delta _{1,1}^{H} \\
&=&((I_{1}\otimes \Delta _{\alpha ,\alpha ^{-1}}^{H})\Delta
_{1,1}^{H}\otimes I_{1})\Delta _{1,1}^{H}\\
&=&(I_{1}\otimes \Delta _{\alpha ,\alpha ^{-1}}^{H}\otimes
I_{1})(\Delta _{1,1}^{H}\otimes I_{1})\Delta _{1,1}^{H}\text{ \ \
\ \ \ \ \ \ \ \ \ \ \ \ \ \ \ \ \ \ \ \ \ \ \ \ \ \ \ \ \ \ \ \ \
\ \ \ \ \ \ \ \ \ \ \ \ \ \ \ \ \ \ \ \ \ \ \ \ \ \ \ \ \ \ \ \ \
\ \ \ \ \ \ \ \ \ \ \ \ \ \ \ \ \ \ \ \ \ \ \ \ \ \ \ \ \ \ \ \ \
\ \ \ \ \ \ \ \ \ \ \ \ \ \ \ \ \ \ \ \ \ \ \ \ \ \ \ \ }
\end{eqnarray*}%
\begin{equation*}
i.e.,(\Delta _{1,\alpha }^{H}\otimes \Delta _{\alpha
^{-1},1}^{H})\Delta _{\alpha ,\alpha
^{-1}}^{H}(h)=h_{11}^{1}\otimes h_{121}^{\alpha }\otimes
h_{122}^{\alpha ^{-1}}\otimes h_{2}^{1}
\end{equation*}%
The second statement is similar.
\end{proof}
Now, from lemma 3.5 up to theorem 3.12, $H$ is Hopf $\pi
-$coalgebra and $(C,\sigma )$ is Hopf $\pi -$subcoalgebra. Also,
$(C,\sigma
)$ have a section $g=\{g_{\alpha }\}_{\alpha \in \pi }$ and antipode $%
S^C=\{S^{C}_{\alpha}: C_\alpha \to
C_{\alpha^{-1}}\}_{\alpha\in\pi}$

\begin{lem}
For $\alpha \in \pi ,L_{1,\alpha }g_{\alpha }^{-1}=(S_{1}^{C}\otimes
g_{\alpha }^{-1})\tau \Delta _{\alpha ^{-1},1}^{C}$
\end{lem}

\begin{proof}
We will prove that for $\alpha \in \pi $, $L_{1,\alpha }g_{\alpha
}^{-1}$ and $(S_{1}^{C}\otimes g_{\alpha }^{-1})\tau \Delta
_{\alpha ^{-1},1}^{C}$ are inverse to the same element
$L_{1,\alpha }g_{\alpha }$ in the convolution algebra
$Conv(C,C_{1}\otimes H_{\alpha }).$ For $c\in C_{1}$,
\begin{eqnarray*}
&&(L_{1,\alpha }g_{\alpha }\ast L_{1,\alpha }g_{\alpha
}^{-1})(c)\\
&=&\{(\sigma _{1}\otimes I_{\alpha }^{H})\Delta _{1,\alpha
}^{H}g_{\alpha }\ast (\sigma _{1}\otimes I_{\alpha }^{H})\Delta
_{1,\alpha }^{H}g_{\alpha
}^{-1}\}(c) \\
&=&\mu _{C_{1}\otimes H_{\alpha }}\{(\sigma _{1}\otimes I_{\alpha
}^{H})\Delta _{1,\alpha }^{H}g_{\alpha }\otimes (\sigma
_{1}\otimes I_{\alpha }^{H})\Delta _{1,\alpha }^{H}g_{\alpha
}^{-1}\}\Delta _{\alpha
,\alpha ^{-1}}^{C}(c) \\
&=&(\mu _{C_{1}}\otimes \mu _{H_{\alpha }})(I\otimes \tau \otimes
I)\{(\sigma _{1}\otimes I_{\alpha }^{H})\Delta _{1,\alpha
}^{H}g_{\alpha }\otimes (\sigma _{1}\otimes I_{\alpha }^{H})\Delta
_{1,\alpha
}^{H}g_{\alpha }^{-1}\}\Delta _{\alpha ,\alpha ^{-1}}^{C}(c) \\
&=&(\mu _{C_{1}}\otimes \mu _{H_{\alpha }})(\sigma _{1}\otimes
\sigma _{1}\otimes I_{\alpha }^{H}\otimes I_{\alpha
}^{H})(I\otimes \tau \otimes I)(\Delta _{1,\alpha }^{H}\otimes
\Delta _{1,\alpha }^{H})(g_{\alpha }\otimes g_{\alpha
}^{-1})\Delta _{\alpha ,\alpha
^{-1}}^{C}(c) \\
&=&(\mu _{C_{1}}(\sigma _{1}\otimes \sigma _{1})\otimes \mu
_{H_{\alpha }})(I\otimes \tau \otimes I)(\Delta _{1,\alpha
}^{H}\otimes \Delta _{1,\alpha }^{H})(g_{\alpha }\otimes g_{\alpha
}^{-1})\Delta _{\alpha
,\alpha ^{-1}}^{C}(c) \\
&=& (\sigma_1 \mu_{H_1}\otimes  \mu _{H_{\alpha }})(I\otimes \tau
\otimes I)(\Delta _{1,\alpha }^{H}\otimes \Delta _{1,\alpha
}^{H})(g_{\alpha }\otimes g_{\alpha }^{-1})\Delta _{\alpha
,\alpha ^{-1}}^{C}(c) \\
&=&(\sigma _{1}\otimes I_{\alpha }^{H})(\mu _{H_{1}}\otimes \mu
_{H_{\alpha }})(I\otimes \tau \otimes I)(\Delta _{1,\alpha
}^{H}\otimes \Delta _{1,\alpha }^{H})(g_{\alpha }\otimes g_{\alpha
}^{-1})\Delta _{\alpha
,\alpha ^{-1}}^{C}(c) \\
&=&(\sigma _{1}\otimes I_{\alpha }^{H})\Delta _{1,\alpha }^{H}\mu
_{H_{\alpha }}(g_{\alpha }\otimes g_{\alpha }^{-1})\Delta _{\alpha
,\alpha ^{-1}}^{C}(c)\\
&=&(\sigma _{1}\otimes I_{\alpha }^{H})\Delta _{1,\alpha
}^{H}(\epsilon
^{C}(c)1_{H_{\alpha }}) \\
&=&\epsilon ^{C}(c)(\sigma _{1}\otimes I_{\alpha }^{H})\Delta
_{1,\alpha }^{H}(1_{H_{\alpha }})\\
&=&\epsilon^{C}(c)(1_{C_{1}}\otimes 1_{H_{\alpha }}).
\end{eqnarray*}%
and%
\begin{eqnarray*}
&& \{L_{1,\alpha }g_{\alpha }\ast (S_{1}^{C}\otimes g_{\alpha
}^{-1})\tau \Delta _{\alpha ^{-1},1}^{C}\}(c) \\
&=&\{(\sigma _{1}\otimes I_{\alpha }^{H})\Delta _{1,\alpha
}^{H}g_{\alpha }\ast (S_{1}^{C}\otimes g_{\alpha }^{-1})\tau
\Delta _{\alpha
^{-1},1}^{C}\}(c) \\
&=&\{(I_{1}\otimes g_{\alpha })\Delta _{1,\alpha }^{C}\ast
(S_{1}^{C}\otimes
g_{\alpha }^{-1})\tau \Delta _{\alpha ^{-1},1}^{C}\}(c) \\
&=&\mu_{C_1\otimes H_{\alpha}}\{(I_{1}\otimes g_{\alpha })\Delta
_{1,\alpha }^{C}\otimes(S_{1}^{C}\otimes g_{\alpha }^{-1})\tau
\Delta _{\alpha ^{-1},1}^{C}\}\Delta
_{\alpha ,\alpha ^{-1}}^{C}(c) \\
&=&(\mu _{C_{1}}\otimes \mu _{H_{\alpha }})(I\otimes \tau \otimes
I)\{(I_{1}\otimes g_{\alpha })\Delta _{1,\alpha
}^{C}\otimes(S_{1}^{C}\otimes g_{\alpha }^{-1})\tau \Delta
_{\alpha ^{-1},1}^{C}\}\Delta
_{\alpha ,\alpha ^{-1}}^{C}(c) \\
&=&(\mu _{C_{1}}\otimes \mu _{H_{\alpha }})(I\otimes \tau \otimes
I)(I_{1}\otimes g_{\alpha }\otimes S_{1}^{C}\otimes g_{\alpha
}^{-1})(I_{1}\otimes I_{\alpha }\otimes \tau )(\Delta _{1,\alpha
}^{C}\otimes
\Delta _{\alpha ^{-1},1}^{C})\Delta _{\alpha ,\alpha ^{-1}}^{C}(c) \\
&=&(\mu _{C_{1}}\otimes \mu _{H_{\alpha }})(I\otimes \tau \otimes
I)(I_{1}\otimes g_{\alpha }\otimes S_{1}^{C}\otimes g_{\alpha
}^{-1})(I_{1}\otimes I_{\alpha }\otimes \tau )(c_{11}^{1}\otimes
c_{121}^{\alpha }\otimes c_{122}^{\alpha ^{-1}}\otimes
c_{2}^{1}) \\
&=&c_{11}^{1}S_{1}^{C}(c_{2}^{1})\otimes g_{\alpha
}(c_{121}^{\alpha })g_{\alpha }^{-1}(c_{122}^{\alpha
^{-1}})\\
&=&c_{11}^{1}S_{1}^{C}(c_{2}^{1})\otimes \epsilon
^{C}(c_{12}^{\alpha })1_{H_{\alpha }} \\
&=&\epsilon ^{C}(c_{12}^{\alpha
})c_{11}^{1}S_{1}^{C}(c_{2}^{1})\otimes 1_{H_{\alpha
}}\\
&=&c_{1}^{1}S_{1}^{C}(c_{2}^{1})\otimes 1_{H_{\alpha }}\\
&=&\epsilon ^{C}(c)(1_{C_{1}}\otimes 1_{H_{\alpha }})
\end{eqnarray*}
\end{proof}

\begin{lem}
For $h\in H_{\alpha },$ we have $h=g_{\alpha }\sigma _{\alpha
}(h_{1}^{\alpha })g_{\alpha }^{-1}\sigma _{\alpha ^{-1}}(h_{21}^{\alpha
^{-1}})h_{22}^{\alpha }.$
\end{lem}

\begin{proof}
\begin{eqnarray*}
h &=&\epsilon ^{H}(h_{1}^{1})h_{2}^{\alpha } \\
&=&\varepsilon ^{C}\sigma _{1}(h_{1}^{1})h_{2}^{\alpha }\\
&=&(\mu _{\alpha }(g_{\alpha }\otimes g_{\alpha }^{-1})\Delta
_{\alpha
,\alpha ^{-1}}^{C}(\sigma _{1}(h_{1}^{1})))h_{2}^{\alpha } \\
&=&(\mu _{\alpha }(g_{\alpha }\otimes g_{\alpha }^{-1})(\sigma
_{\alpha }\otimes \sigma _{\alpha ^{-1}})\Delta _{\alpha ,\alpha
^{-1}}^{H}(h_{1}^{1}))h_{2}^{\alpha } \\
&=&\mu _{\alpha }(\mu _{\alpha }(g_{\alpha }\otimes g_{\alpha
}^{-1})(\sigma _{\alpha }\otimes \sigma _{\alpha ^{-1}})\otimes
I_{\alpha })(\Delta
_{\alpha ,\alpha ^{-1}}^{H}\otimes I_{\alpha })\Delta _{1,\alpha }^{H}(h) \\
&=&\mu _{\alpha }(\mu _{\alpha }(g_{\alpha }\otimes g_{\alpha
}^{-1})(\sigma _{\alpha }\otimes \sigma _{\alpha ^{-1}})\otimes
I_{\alpha })(I_{\alpha
}\otimes \Delta _{\alpha ^{-1},\alpha }^{H})\Delta _{\alpha ,1}^{H}(h) \\
&=&g_{\alpha }\sigma _{\alpha }(h_{1}^{\alpha })g_{\alpha
}^{-1}\sigma _{\alpha ^{-1}}(h_{21}^{\alpha ^{-1}})h_{22}^{\alpha
}.
\end{eqnarray*}
\end{proof}

\begin{lem}
For $h\in H_{1}$ and $\alpha \in \pi ,$ we have $\mu _{\alpha }(g_{\alpha
}^{-1}\sigma _{\alpha ^{-1}}\otimes I_{\alpha })\Delta _{\alpha ,\alpha
^{-1}}^{H}(h)\in B_{\alpha }$.
\end{lem}

\begin{proof}
\begin{eqnarray*}
&& L_{1,\alpha }\mu _{\alpha }(g_{\alpha }^{-1}\sigma _{\alpha
^{-1}}\otimes I_{\alpha })\Delta _{\alpha ^{-1},\alpha }^{H}(h)\\
&=&(\sigma _{1}\otimes I_{\alpha }^{H})\Delta _{1,\alpha }^{H}\mu
_{\alpha }(g_{\alpha }^{-1}\sigma _{\alpha ^{-1}}\otimes I_{\alpha
})\Delta _{\alpha
^{-1},\alpha }^{H}(h) \\
&=&\mu _{C_{1}\otimes H_{\alpha }}((\sigma _{1}\otimes I_{\alpha
}^{H})\Delta _{1,\alpha }^{H}g_{\alpha }^{-1}\sigma _{\alpha
^{-1}}\otimes (\sigma _{1}\otimes I_{\alpha }^{H})\Delta
_{1,\alpha }^{H})\Delta _{\alpha
^{-1},\alpha }^{H}(h) \\
&=&\mu _{C_{1}\otimes H_{\alpha }}((S_{1}^{C}\otimes g_{\alpha
}^{-1})\tau \Delta _{\alpha ^{-1},1}^{C}\sigma _{\alpha
^{-1}}\otimes (\sigma _{1}\otimes I_{\alpha }^{H})\Delta
_{1,\alpha }^{H})\Delta _{\alpha
^{-1},\alpha }^{H}(h) \\
&=&\mu _{C_{1}\otimes H_{\alpha }}((S_{1}^{C}\otimes g_{\alpha
}^{-1})\tau (\sigma _{\alpha ^{-1}}\otimes \sigma _{1})\Delta
_{\alpha ^{-1},1}^{H}\otimes(\sigma _{1}\otimes I_{\alpha
}^{H})\Delta _{1,\alpha }^{H})\Delta _{\alpha
^{-1},\alpha }^{H}(h) \\
&=&\mu _{C_{1}\otimes H_{\alpha }}((S_{1}^{C}\otimes g_{\alpha
}^{-1})\tau \otimes I_{1}\otimes I_{\alpha }^{H})(\sigma _{\alpha
^{-1}}\otimes \sigma _{1}\otimes \sigma _{1}\otimes I_{\alpha
}^{H})(\Delta _{\alpha ^{-1},1}^{H}\otimes \Delta _{1,\alpha
}^{H})\Delta _{\alpha
^{-1},\alpha }^{H}(h) \\
&=&\mu _{C_{1}\otimes H_{\alpha }}((S_{1}^{C}\otimes g_{\alpha
}^{-1})\tau \otimes I_{1}\otimes I_{\alpha }^{H})(\sigma _{\alpha
^{-1}}(h_{11}^{\alpha ^{-1}})\otimes \sigma
_{1}(h_{121}^{1})\otimes \sigma _{1}(h_{122}^{1})\otimes h_{2}^{\alpha }) \\
&=&S_{1}^{C}(\sigma _{1}(h_{121}^{1}))\sigma
_{1}(h_{122}^{1})\otimes g_{\alpha }^{-1}\sigma _{\alpha
^{-1}}(h_{11}^{\alpha ^{-1}})h_{2}^{\alpha } \\
&=&\epsilon ^{C}(\sigma _{1}(h_{12}^{1})1\otimes g_{\alpha
}^{-1}\sigma
_{\alpha ^{-1}}(h_{11}^{\alpha ^{-1}})h_{2}^{\alpha } \\
&=&\epsilon ^{H}(h_{12}^{1})1\otimes g_{\alpha }^{-1}\sigma
_{\alpha ^{-1}}(h_{11}^{\alpha ^{-1}})h_{2}^{\alpha } \\
&=&1\otimes g_{\alpha }^{-1}\sigma _{\alpha ^{-1}}(\epsilon
^{H}(h_{12}^{1})h_{11}^{\alpha ^{-1}})h_{2}^{\alpha } \\
&=&1\otimes g_{\alpha }^{-1}\sigma _{\alpha ^{-1}}(h_{1}^{\alpha
^{-1}})h_{2}^{\alpha } \\
&=&1\otimes \mu _{\alpha }(g_{\alpha }^{-1}\sigma _{\alpha
^{-1}}\otimes I_{\alpha })\Delta _{\alpha ^{-1},\alpha }^{H}(h).
\end{eqnarray*}
\end{proof}

\begin{lem}
For $\alpha \in \pi ,$%
\begin{equation*}
(\sigma _{\alpha }\otimes \sigma _{\alpha ^{-1}}\otimes I_{\alpha
}^{H})(\Delta _{\alpha ,\alpha ^{-1}}^{H}\otimes I_{\alpha }^{H})\Delta
_{1,\alpha }^{H}g_{\alpha } =(I_{\alpha }\otimes I_{\alpha ^{-1}}\otimes
g_{\alpha })(I_{\alpha }\otimes \Delta _{\alpha ^{-1},\alpha }^{C})\Delta
_{\alpha ,1}^{C}
\end{equation*}
\end{lem}

\begin{proof}
\begin{eqnarray*}
&& (\sigma _{\alpha }\otimes \sigma _{\alpha ^{-1}}\otimes
I_{\alpha }^{H})(\Delta _{\alpha ,\alpha ^{-1}}^{H}\otimes
I_{\alpha }^{H})\Delta _{1,\alpha }^{H}g_{\alpha } \\
&=&((\sigma _{\alpha }\otimes \sigma _{\alpha ^{-1}})\Delta
_{\alpha ,\alpha ^{-1}}^{H}\otimes I_{\alpha }^{H})\Delta
_{1,\alpha }^{H}g_{\alpha
} \\
&=&(\Delta _{\alpha ,\alpha ^{-1}}^{C}\sigma _{1}\otimes I_{\alpha
}^{H})\Delta _{1,\alpha }^{H}g_{\alpha }\\
&=&(\Delta _{\alpha ,\alpha ^{-1}}^{C}\otimes I_{\alpha
}^{H})(\sigma _{1}\otimes I_{\alpha }^{H})\Delta _{1,\alpha
}^{H}g_{\alpha } \\
&=&(\Delta _{\alpha ,\alpha ^{-1}}^{C}\otimes I_{\alpha
}^{H})(I_{1}\otimes g_{\alpha })\Delta
_{1,\alpha }^{C}\\
&=&(I_{\alpha }\otimes I_{\alpha ^{-1}}\otimes g_{\alpha })(\Delta
_{\alpha ,\alpha ^{-1}}^{C}\otimes I_{\alpha })\Delta _{1,\alpha
}^{C} \\
&=&(I_{\alpha }\otimes I_{\alpha ^{-1}}\otimes g_{\alpha
})(I_{\alpha }\otimes \Delta _{\alpha ^{-1},\alpha }^{C})\Delta
_{\alpha ,1}^{C}\text{ \ \ \ \ \ \ \ \ \ \ \ \ \ \ \ \ \ \ \ \ \ \
\ \ \ \ \ \ \ \ \ \ \ \ \ \ \ \ \ \ \ \ \ \ \ \ \ \ \ \ \ \ \ \ \
\ \ \ \ \ \ \ \ \ \ \ \ \ \ \ \ \ \ \ \ \ \ \ \ \ \ \ \ \ \ \ \ \
\ \ \ \ \ \ \ \ \ \ \ \ \ \ \ \ \ \ \ \ \ \ \ \ \ \ \ \ \ \ \ \ \
\ \ \ \ \ \ \ \ }
\end{eqnarray*}
\end{proof}

\begin{lem}
For $\alpha \in \pi $ and $b\in B_{\alpha }$
\begin{equation*}
(\sigma _{\alpha }\otimes \sigma _{\alpha ^{-1}}\otimes I_{\alpha
}^{H})(\Delta _{\alpha ,\alpha ^{-1}}^{H}\otimes I_{\alpha }^{H})\Delta
_{1,\alpha }^{H}(b)=1\otimes 1\otimes b
\end{equation*}
\end{lem}

\begin{proof}
\begin{eqnarray*}
&&(\sigma _{\alpha }\otimes \sigma _{\alpha ^{-1}}\otimes
I_{\alpha }^{H})(\Delta _{\alpha ,\alpha ^{-1}}^{H}\otimes
I_{\alpha }^{H})\Delta _{1,\alpha }^{H}(b) \\
&=&((\sigma _{\alpha }\otimes \sigma _{\alpha ^{-1}})\Delta
_{\alpha ,\alpha ^{-1}}^{H}\otimes I_{\alpha }^{H})\Delta
_{1,\alpha
}^{H}(b)\\
&=&(\Delta _{\alpha ,\alpha ^{-1}}^{C}\sigma _{1}\otimes I_{\alpha
}^{H})\Delta _{1,\alpha }^{H}(b) \\
&=&(\Delta _{\alpha ,\alpha ^{-1}}^{C}\otimes I_{\alpha
}^{H})(\sigma _{1}\otimes I_{\alpha }^{H})\Delta _{1,\alpha
}^{H}(b)\\
&=&(\Delta _{\alpha ,\alpha ^{-1}}^{C}\otimes I_{\alpha
}^{H})(1\otimes b)\\
&=&1\otimes 1\otimes b \text{ \ \ \ \ \ \ \ \ \ \ \ \ \ \ \ \ \ \
\ \ \ \ \ \ \ \ \ \ \ \ \ \ \ \ \ \ \ \ \ \ \ \ \ \ \ \ \ \ \ \ \
\ \ \ \ \ \ \ \ \ \ \ \ \ \ \ \ \ \ \ \ \ \ \ \ \ \ \ \ \ \ \ \ \
\ \ \ \ \ \ \ \ \ \ \ \ \ \ \ \ \ \ \ \ \ \ \ \ \ \ \ \ \ \ \ \ \
\ \ \ \ \ \ \ \ \ \ \ \ }
\end{eqnarray*}
\end{proof}

\begin{thm}
$H$ is isomorphic to $C\otimes B$ as vector space.
\end{thm}

\begin{proof}
We define $A=\{A_{\alpha }:C_{\alpha }\otimes B_{\alpha
}\rightarrow H_{\alpha }\}_{\alpha \in \pi }$ as follow
\begin{equation*}
A_{\alpha }(c\otimes b)=\mu _{\alpha }(g_{\alpha }\otimes
I_{\alpha })(c\otimes b)=g_{\alpha }(c)b.
\end{equation*}%
Clear $A_{\alpha }$ is linear and by lemma 3.6 and lemma 3.7 that
$A_{\alpha }$ is surjective for all $\alpha \in \pi .$ We define
$A_{\alpha }^{-1}=(\sigma _{\alpha }\otimes \mu _{\alpha
}(g_{\alpha }^{-1}\sigma _{\alpha ^{-1}}\otimes I_{\alpha
}))(\Delta _{\alpha ,\alpha ^{-1}}^{H}\otimes I_{\alpha })\Delta
_{1,\alpha }^{H}.$ We`ll prove that for $\alpha \in \pi ,A_{\alpha
}A_{\alpha }^{-1}=I_{H_{\alpha }}$ and $A_{\alpha }^{-1}A_{\alpha
}=I_{C_{\alpha }\otimes B_{\alpha }}.$ Firstly,  let $h\in
H_{\alpha }$
and $c\otimes b\in C_{\alpha }\otimes B_{\alpha }$%
\begin{eqnarray*}
&&A_{\alpha }A_{\alpha }^{-1}(h) \\
&=&A_{\alpha }((\sigma _{\alpha }\otimes \mu _{\alpha }(g_{\alpha
}^{-1}\sigma _{\alpha ^{-1}}\otimes I_{\alpha }))(\Delta _{\alpha
,\alpha
^{-1}}^{H}\otimes I_{\alpha })\Delta _{1,\alpha }^{H}(h)) \\
&=&A_{\alpha }((\sigma _{\alpha }\otimes \mu _{\alpha }(g_{\alpha
}^{-1}\sigma _{\alpha ^{-1}}\otimes I_{\alpha }))(I_{\alpha
}\otimes \Delta
_{\alpha ^{-1},\alpha }^{H})\Delta _{\alpha ,1}^{H}(h)) \\
&=&A_{\alpha }(\sigma _{\alpha }(h_{1}^{\alpha })\otimes g_{\alpha
}^{-1}\sigma _{\alpha ^{-1}}(h_{21}^{\alpha ^{-1}})h_{22}^{\alpha
})\\
&=&g_{\alpha }\sigma _{\alpha }(h_{1}^{\alpha })g_{\alpha
}^{-1}\sigma _{\alpha ^{-1}}(h_{21}^{\alpha ^{-1}})h_{22}^{\alpha
} \\
&=&h~\
\;\;\;\;\;\;\;\;\;\;\;\;\;\;\;\;\;\;\;\;\;\;\;\;\;\;\;\;\;\;\;\;\;\;\;\;\;\;\;\;\;\;\;\;\;\;\;\;\;\;\;\;\;\;\;\;\;\;\;\;\;\;\;\;\;\;\;\;\;\;\;\;
\text{by lemma 3.6}\text{ \ \ \ \ \ \ \ \ \ \ \ \ \ \ \ \ \ \ \ \
\ \ \ \ \ \ \ \ \ \ \ \ \ \ \ \ \ \ \ \ \ \ \ \ \ \ \ \ \ \ \ \ \
\ \ \ \ \ \ \ \ \ \ \ \ \ \ \ \ \ \ \ \ \ \ \ \ \ \ \ \ \ \ \ \ \
\ \ \ \ \ \ \ \ \ \ \ \ \ \ \ \ \ \ \ \ \ \ \ \ \ \ \ \ \ \ \ \ \
\ \ \ \ \ \ \ \ \ \ }
\end{eqnarray*}%
Secondly, since, for $\alpha \in \pi $ we have $((\sigma _{\alpha
}\otimes \sigma _{\alpha ^{-1}})\Delta _{\alpha ,\alpha
^{-1}}^{H}\otimes I_{\alpha }^{H})\Delta _{1,\alpha }^{H}$ is an
algebra map, then
\begin{eqnarray*}
&&A_{\alpha }^{-1}A_{\alpha }(c\otimes b)\\
&=&(\sigma _{\alpha }\otimes \mu _{\alpha }(g_{\alpha }^{-1}\sigma
_{\alpha ^{-1}}\otimes I_{\alpha }))(\Delta _{\alpha ,\alpha
^{-1}}^{H}\otimes I_{\alpha })\Delta _{1,\alpha }^{H}\mu _{\alpha
}(g_{\alpha }\otimes
I_{\alpha })(c\otimes b) \\
&=&(I_{\alpha }\otimes \mu _{\alpha })(I_{\alpha }\otimes
g_{\alpha }^{-1}\otimes I_{\alpha })(\sigma _{\alpha }\otimes
\sigma _{\alpha ^{-1}}\otimes I_{\alpha }^{H})(\Delta _{\alpha
,\alpha ^{-1}}^{H}\otimes I_{\alpha }^{H})\Delta _{1,\alpha
}^{H}\mu _{\alpha }(g_{\alpha }\otimes
I_{\alpha })(c\otimes b) \\
&=&(I_{\alpha }\otimes \mu _{\alpha })(I_{\alpha }\otimes
g_{\alpha }^{-1}\otimes I_{\alpha })((\sigma _{\alpha }\otimes
\sigma _{\alpha ^{-1}})\Delta _{\alpha ,\alpha ^{-1}}^{H}\otimes
I_{\alpha }^{H})\Delta _{1,\alpha }^{H}\mu _{\alpha }(g_{\alpha
}\otimes I_{\alpha
})(c\otimes b) \\
&=&(I_{\alpha }\otimes \mu _{\alpha })(I_{\alpha }\otimes
g_{\alpha }^{-1}\otimes I_{\alpha })\mu _{H_{\alpha }\otimes
H_{\alpha ^{-1}}\otimes
H_{\alpha }} \\
&&(((\sigma _{\alpha }\otimes \sigma _{\alpha ^{-1}})\Delta
_{\alpha ,\alpha ^{-1}}^{H}\otimes I_{\alpha }^{H})\Delta
_{1,\alpha }^{H}g_{\alpha }\otimes ((\sigma _{\alpha }\otimes
\sigma _{\alpha ^{-1}})\Delta _{\alpha ,\alpha
^{-1}}^{H}\otimes I_{\alpha }^{H})\Delta _{1,\alpha }^{H})(c\otimes b) \\
&=&(I_{\alpha }\otimes \mu _{\alpha })(I_{\alpha }\otimes
g_{\alpha
}^{-1}\otimes I_{\alpha })\mu _{H_{\alpha }\otimes H_{\alpha ^{-1}}} \\
&&((I_{\alpha }\otimes I_{\alpha ^{-1}}\otimes g_{\alpha
})(I_{\alpha }\otimes \Delta _{\alpha ^{-1},\alpha }^{C})\Delta
_{\alpha
,1}^{C}(c)\otimes 1\otimes 1\otimes b) \;\;\;\;\;\;\;\;\; \textrm{by lemmas 3.8, 3.9}\\
&=&c_{1}^{\alpha }\otimes g_{\alpha }^{-1}(c_{21}^{\alpha
^{-1}})g_{\alpha }(c_{22}^{\alpha })b\\
&=&c_{1}^{\alpha }\otimes
\epsilon (c_{2}^{1})b\\
&=&\epsilon (c_{2}^{1})c_{1}^{\alpha }\otimes
b\\
&=&c\otimes b.
\end{eqnarray*}%
\end{proof}

\begin{rem}
In any Hopf $\pi $ -coalgebra $H$, every $H_{\alpha }$ is left $H_{1}$
-comodule by $\Delta _{1,\alpha }$. i.e., the following diagrams are commute
\begin{equation*}
\begin{tabular}{lllllllll}
& $H_{\alpha }$ & $%
\begin{tabular}{l}
$\Delta _{1,\alpha }$ \\
$\rightarrow $%
\end{tabular}%
$ & $H_{1}\otimes H_{\alpha }$ &  & $H_{\alpha }$ & $%
\begin{tabular}{l}
$\Delta _{1,\alpha }$ \\
$\rightarrow $%
\end{tabular}%
$ & $H_{1}\otimes H_{\alpha }$ &  \\
$\Delta _{1,\alpha }$ & $\downarrow $ &  & $\downarrow \Delta _{1,1}\otimes
I_{\alpha }$ &  &  & $\sim \searrow $ & $\downarrow \epsilon ^{H}\otimes
I_{\alpha }$ &  \\
& $H_{1}\otimes H_{\alpha }$ &
\begin{tabular}{l}
$\longrightarrow $ \\
$I_{1}\otimes \Delta _{1,\alpha }$%
\end{tabular}
& $H_{1}\otimes H_{1}\otimes H_{\alpha }$ &  &  &  & $K\otimes H_{\alpha }$
&
\end{tabular}%
\end{equation*}
\end{rem}

A Hopf $\pi -$coalgebra $H$ is called have a left (right) cosection if there
exist a family of algebra maps $\eta =\{\eta _{\alpha }:H_{1}\rightarrow
H_{\alpha }\}_{\alpha \in \pi }$ such that $\eta _{\alpha }\,$is left
(right) $H_{1}$-comodule map $\forall\alpha \in \pi $ i.e., the following
diagram is commute
\begin{equation*}
\begin{tabular}{llll}
& $H_{1}$ & $\overset{\eta _{\alpha }}{\longrightarrow }$ & $H_{\alpha }$ \\
$\Delta _{1,1}$ & $\downarrow $ &  & $\downarrow \Delta _{1,\alpha }$ \\
& $H_{1}\otimes H_{1}$ &
\begin{tabular}{l}
$\longrightarrow $ \\
$I_{1}\otimes \eta _{\alpha }$%
\end{tabular}
& $H_{1}\otimes H_{\alpha }$%
\end{tabular}%
\end{equation*}

\begin{thm}
If $H$ have a left cosection, then $Ind(\rho )$ is isomorphic to $V\otimes B=
$ $\{(V\otimes B)_{\alpha }=V_{1}\otimes B_{\alpha }\}_{\alpha \in \pi }$ as
right $B$-module.
\end{thm}

\begin{proof}
Firstly, we`ll prove that $L_{1,\alpha }\eta _{\alpha
}g_{1}=(I_{1}\otimes \eta _{\alpha }g_{1})\Delta _{1,1}^{C}.$
\begin{eqnarray*}
&&L_{1,\alpha }\eta _{\alpha }g_{1} \\
&=&(\sigma _{1}\otimes I_{\alpha }^{H})\Delta _{1,\alpha }^{H}\eta
_{\alpha }g_{1}\\
&=&(\sigma _{1}\otimes I_{\alpha }^{H})(I_{1}\otimes \eta _{\alpha
})\Delta _{1,1}^{H}g_{1} \\
&=&(I_{1}\otimes \eta _{\alpha })(\sigma _{1}\otimes
I_{1}^{H})\Delta
_{1,1}^{H}g_{1} \\
&=&(I_{1}\otimes \eta _{\alpha })L_{1,1}g_{1} \\
&=&(I_{1}\otimes \eta _{\alpha })(I_{1}\otimes g_{1})\Delta
_{1,1}^{C} \\
&=&(I_{1}\otimes \eta _{\alpha }g_{1})\Delta _{1,1}^{C}.\text{ \ \
\ \ \ \ \ \ \ \ \ \ \ \ \ \ \ \ \ \ \ \ \ \ \ \ \ \ \ \ \ \ \ \ \
\ \ \ \ \ \ \ \ \ \ \ \ \ \ \ \ \ \ \ \ \ \ \ \ \ \ \ \ \ \ \ \ \
\ \ \ \ \ \ \ \ \ \ \ \ \ \ \ \ \ \ \ \ \ \ \ \ \ \ \ \ \ \ \ \ \
\ \ \ \ \ \ \ \ \ \ \ \ \ \ \ \ \ \ \ \ \ \ \ \ \ \ \ \ }
\end{eqnarray*}
We define $T_{\alpha }=(I_{1}\otimes \eta _{\alpha }g_{1})\rho
_{1,1}:V_{1}\rightarrow Ind(\rho )_{\alpha }.$ We`ll prove
$T_{\alpha }(v_{1})\in Ind(\rho )_{\alpha }.$
\begin{eqnarray*}
&&(I_{1}\otimes L_{1,\alpha })T_{\alpha }(v_{1})\\
&=&(I_{1}\otimes L_{1,\alpha })(I_{1}\otimes \eta _{\alpha
}g_{1})\rho
_{1,1}(v_{1})\\
&=&(I_{1}\otimes L_{1,\alpha }\eta _{\alpha }g_{1})\rho _{1,1}(v_{1}) \\
&=&(I_{1}\otimes (I_{1}\otimes \eta _{\alpha }g_{1})\Delta
_{1,1}^{C})\rho _{1,1}(v_{1})\\
&=&(I_{1}\otimes I_{1}\otimes \eta _{\alpha }g_{1})(I_{1}\otimes
\Delta _{1,1}^{C})\rho _{1,1}(v_{1}) \\
&=&(I_{1}\otimes I_{1}\otimes \eta _{\alpha }g_{1})(\rho
_{1,1}\otimes I_{1})\rho _{1,1}(v_{1}) \\
&=&(\rho _{1,1}\otimes I_{\alpha })(I_{1}\otimes \eta _{\alpha
}g_{1})\rho _{1,1}(v_{1}) \\
&=&(\rho _{1,1}\otimes I_{\alpha })T_{\alpha }(v_{1}).\text{ \ \ \
\ \ \ \ \ \ \ \ \ \ \ \ \ \ \ \ \ \ \ \ \ \ \ \ \ \ \ \ \ \ \ \ \
\ \ \ \ \ \ \ \ \ \ \ \ \ \ \ \ \ \ \ \ \ \ \ \ \ \ \ \ \ \ \ \ \
\ \ \ \ \ \ \ \ \ \ \ \ \ \ \ \ \ \ \ \ \ \ \ \ \ \ \ \ \ \ \ \ \
\ \ \ \ \ \ \ \ \ \ \ \ \ \ \ \ \ \ \ \ \ \ \ \ \ \ \ }
\end{eqnarray*}
For $\alpha \in \pi $ we define $q_{\alpha }:V_{1}\otimes
B_{\alpha }\rightarrow Ind(\rho )_{\alpha }$ where
\begin{equation*}
q_{\alpha }(v\otimes b)=\lambda _{\alpha }(T_{\alpha }(v)\otimes
b)=v_{1}\otimes \eta _{\alpha }g_{1}(v_{2})b.
\end{equation*}%
Clear $q_{\alpha }(v\otimes b)\in Ind(\rho )_{\alpha }$ and
$q_{\alpha }$ is linear. We define $q_{\alpha }^{-1}:Ind(\rho
)_{\alpha }\rightarrow V_{1}\otimes B_{\alpha }$;
\begin{eqnarray*}
q_{\alpha }^{-1} &=&(I_{1}\otimes \mu _{\alpha }(\eta _{\alpha
}g_{1}^{-1}\otimes I_{\alpha }))(\rho _{1,1}\otimes I_{\alpha }) \;\;\;\;\;\;\;\;\;\;\;\;\;\;\;\;\;\;\;\;\;\;\;\;\;\;\;\;\;\;\;\;\;\;\;\;\;\;\; \text{(i)}\\
&=&(I_{1}\otimes \mu _{\alpha }(\eta _{\alpha }g_{1}^{-1}\otimes
I_{\alpha }))(I_{1}\otimes L_{1,\alpha })
\;\;\;\;\;\;\;\;\;\;\;\;\;\;\;\;\;\;\;\;\;\;\;\;\;\;\;\;\;\;\;\;\;\;\;\;\;\;\;
\text{(ii)}
\end{eqnarray*}%
\begin{equation*}
i.e.,q_{\alpha }^{-1}(v\otimes h)=v_{1}\otimes \eta _{\alpha
}g_{1}^{-1}(v_{2})h=v\otimes \eta _{\alpha }g_{1}^{-1}\sigma
_{1}(h_{1}^{1})h_{2}^{\alpha }.
\end{equation*}
We`ll prove that $q_{\alpha }^{-1}(v\otimes h)\in V_{1}\otimes
B_{\alpha }.$
\begin{eqnarray*}
&&(I_{1}\otimes L_{1,\alpha })q_{\alpha }^{-1}(v\otimes h)\\
&=&(I_{1}\otimes L_{1,\alpha })(I_{1}\otimes \mu _{\alpha }(\eta
_{\alpha
}g_{1}^{-1}\otimes I_{\alpha }))(I_{1}\otimes L_{1,\alpha })(v\otimes h)\:\;\;\;\;\;\;\;\;\;\;\;\;\;\text{by equation (ii)} \\
&=&(I_{1}\otimes \mu _{C_{1}\otimes H_{\alpha }}(L_{1,\alpha }\eta
_{\alpha
}g_{1}^{-1}\otimes L_{1,\alpha }))(I_{1}\otimes L_{1,\alpha })(v\otimes h) \\
&=&v\otimes \mu _{C_{1}\otimes H_{\alpha }}((I_{1}\otimes \eta
_{\alpha
})L_{1,1}g_{1}^{-1}\otimes L_{1,\alpha })L_{1,\alpha }(h) \\
&=&v\otimes \mu _{C_{1}\otimes H_{\alpha }}((I_{1}\otimes \eta
_{\alpha })(S_{1}^{C}\otimes g_{1}^{-1})\tau \Delta
_{1,1}^{H}\otimes L_{1,\alpha
})L_{1,\alpha }(h) \:\:\:\:\:\:\:\:\:\:\:\:\textrm{ by lemma 3.5}\\
&=&v\otimes \mu _{C_{1}\otimes H_{\alpha }}((S_{1}^{C}\otimes \eta
_{\alpha }g_{1}^{-1})\tau \Delta _{1,1}^{H}\otimes (\sigma
_{1}\otimes I_{\alpha
}^{H})\Delta _{1,\alpha }^{H}) (\sigma _{1}\otimes I_{\alpha }^{H})\Delta _{1,\alpha }^{H}(h) \\
&=&v\otimes \mu _{C_{1}\otimes H_{\alpha }}((S_{1}^{C}\otimes \eta
_{\alpha }g_{1}^{-1})\tau \otimes \sigma _{1}\otimes I_{\alpha
}^{H})(\Delta _{1,1}^{H}\otimes \Delta _{1,\alpha }^{H})(\sigma
_{1}\otimes I_{\alpha }^{H})\Delta _{1,\alpha }^{H}(h)
\\
&=&v\otimes \mu _{C_{1}\otimes H_{\alpha }}((S_{1}^{C}\otimes \eta
_{\alpha }g_{1}^{-1})\tau \otimes \sigma _{1}\otimes I_{\alpha
}^{H})(\sigma _{1}\otimes \sigma _{1}\otimes I_{1}^{H}\otimes
I_{\alpha }^{H})(\Delta _{1,1}^{H}\otimes \Delta _{1,\alpha
}^{H})\Delta _{1,\alpha
}^{H}(h) \\
&=&v\otimes S_{1}^{C}\sigma _{1}(h_{121}^{1})\sigma
_{1}(h_{122}^{1})\otimes \eta _{\alpha }g_{1}^{-1}\sigma
_{1}(h_{11}^{1}))h_{2}^{\alpha }\\
&=&v\otimes \epsilon ^{C}\sigma _{1}(h_{12}^{1})1\otimes \eta
_{\alpha }g_{\alpha }^{-1}\sigma
_{1}(h_{11}^{1}))h_{2}^{\alpha } \\
&=&v\otimes \epsilon ^{H}(h_{12}^{1})1\otimes \eta _{\alpha
}g_{1}^{-1}\sigma _{1}(h_{11}^{1}))h_{2}^{\alpha }\\
&=&v\otimes 1\otimes \eta _{\alpha }g_{1}^{-1}\sigma _{1}(\epsilon
^{H}(h_{12}^{1})h_{11}^{1}))h_{2}^{\alpha } \\
&=&v\otimes 1\otimes \eta _{\alpha }g_{1}^{-1}\sigma
_{1}(h_{1}^{1})h_{2}^{\alpha }\\
&=&v\otimes 1\otimes \mu _{\alpha }(\eta _{\alpha
}g_{1}^{-1}\otimes I_{\alpha }))L_{1,\alpha }(h)
\end{eqnarray*}%
Now, we prove that $q_{\alpha }q_{\alpha }^{-1}=I$ and $q_{\alpha
}^{-1}q_{\alpha }=I$.
\begin{eqnarray*}
&&q_{\alpha }q_{\alpha }^{-1}(v\otimes h) \\
&=&(I_{1}\otimes \mu _{\alpha })(I_{1}\otimes \eta _{\alpha
}g_{1}\otimes I_{\alpha })(\rho _{1,1}\otimes I_{\alpha
})\\
&&(I_{1}\otimes \mu _{\alpha }(\eta _{\alpha }g_{1}^{-1}\otimes
I_{\alpha
}))(\rho _{1,1}\otimes I_{\alpha })(v\otimes h) \ \ \ \ \ \ \ \ \ \ \ \ \ \ \ \ \ \ \ \ \ \ \ \ \ \text{by equation (ii)}\\
&=&(I_{1}\otimes \mu _{\alpha })(I_{1}\otimes \eta _{\alpha
}g_{1}\otimes I_{\alpha })(I_{1}\otimes I_{1}\otimes \mu _{\alpha
}(\eta _{\alpha }g_{1}^{-1}\otimes I_{\alpha }))(\rho
_{1,1}\otimes I_{1}\otimes I_{\alpha })(\rho _{1,1}\otimes
I_{\alpha })(v\otimes h)\\
&=&(I_{1}\otimes \mu _{\alpha })(I_{1}\otimes \eta _{\alpha
}g_{1}\otimes I_{\alpha })(I_{1}\otimes I_{1}\otimes \mu _{\alpha
}(\eta _{\alpha }g_{1}^{-1}\otimes I_{\alpha }))((\rho
_{1,1}\otimes I_{1})\rho _{1,1}\otimes I_{\alpha })(v\otimes h)\\
&=&(I_{1}\otimes \mu _{\alpha })(I_{1}\otimes \eta _{\alpha
}g_{1}\otimes I_{\alpha })(I_{1}\otimes I_{1}\otimes \mu _{\alpha
}(\eta _{\alpha }g_{1}^{-1}\otimes I_{\alpha }))((I_{1}\otimes
\Delta _{1,1}^{C}\otimes I_{\alpha })(\rho _{1,1}\otimes I_{\alpha
})(v\otimes h) \\
&=&v_{1}^{1}\otimes \eta _{\alpha }g_{1}(v_{21}^{1})\eta _{\alpha
}g_{1}^{-1}(v_{22}^{1})h \\
&=&v_{1}^{1}\otimes \eta _{\alpha
}(g_{1}(v_{21}^{1})g_{1}^{-1}(v_{22}^{1}))h \\
&=&v_{1}^{1}\otimes \eta _{\alpha }(\epsilon
(v_{2}^{1})1)h\\
&=&\epsilon (v_{2}^{1})v_{1}^{1}\otimes \eta _{\alpha }(1)h
\\
&=&v\otimes h
\end{eqnarray*}%
\begin{eqnarray*}
&&q_{\alpha }^{-1}q_{\alpha }(v\otimes b) \\
&=&(I_{1}\otimes \mu _{\alpha }(\eta _{\alpha }g_{1}^{-1}\otimes
I_{\alpha }))(I_{1}\otimes L_{1,\alpha })(I_{1}\otimes \mu
_{\alpha })(I_{1}\otimes \eta _{\alpha }g_{1}\otimes I_{\alpha
})(\rho _{1,1}\otimes
I_{\alpha })(v\otimes b)\text{by equation (i)}\\
&=&(I_{1}\otimes \mu _{\alpha }(\eta _{\alpha }g_{1}^{-1}\otimes
I_{\alpha }))(I_{1}\otimes L_{1,\alpha })(v_{1}\otimes \eta
_{\alpha
}g_{1}(v_{2}^{1})b) \\
&=&(I_{1}\otimes \mu _{\alpha }(\eta _{\alpha }g_{1}^{-1}\otimes
I_{\alpha }))(v_{1}\otimes L_{1,\alpha }(\eta _{\alpha
}g_{1}(v_{2}^{1}))L_{1,\alpha
}(b)) \\
&=&(I_{1}\otimes \mu _{\alpha }(\eta _{\alpha }g_{1}^{-1}\otimes
I_{\alpha }))(v_{1}\otimes (I_{1}\otimes \eta _{\alpha
})L_{1,1}g_{1}(v_{2}^{1})\cdot
(1\otimes b)) \\
&=&(I_{1}\otimes \mu _{\alpha }(\eta _{\alpha }g_{1}^{-1}\otimes
I_{\alpha }))(v_{1}\otimes (I_{1}\otimes \eta _{\alpha
})(I_{1}\otimes g_{1})\Delta
_{1,1}^{C}(v_{2}^{1})\cdot (1\otimes b)) \\
&=&(I_{1}\otimes \mu _{\alpha }(\eta _{\alpha }g_{1}^{-1}\otimes
I_{\alpha }))(v_{1}\otimes v_{21}^{1}\otimes \eta _{\alpha
}g_{1}(v_{22}^{1})b) \\
&=&v_{1}\otimes \eta _{\alpha }g_{1}^{-1}(v_{21}^{1})\eta _{\alpha
}g_{1}(v_{22}^{1})b \\
&=&v_{1}\otimes \eta _{\alpha
}(g_{1}^{-1}(v_{21}^{1})g_{1}(v_{22}^{1}))b\\
&=&v_{1}\otimes \eta _{\alpha }(\epsilon ^{C}(v_{2}^{1})1)b\\
&=&\epsilon ^{C}(v_{2}^{1})v_{1}\otimes \eta _{\alpha
}(1)b\\
&=&v\otimes b.
\end{eqnarray*}
Now we prove that $q_{\alpha }$ is module map for all $\alpha \in \pi ,$%
\begin{eqnarray*}
&&q_{\alpha }(I_{1}\otimes \mu _{\alpha })(v\otimes b\otimes k) \\
&=&q_{\alpha }(v\otimes bk)\\
&=&\lambda _{\alpha }(T_{\alpha }(v)\otimes bk)\\
&=&\lambda _{\alpha }((I_{1}\otimes \eta _{\alpha }g_{1})\rho
_{1,1}(v)\otimes bk)\\
&=&\lambda _{\alpha }(v_{1}\otimes
\eta _{\alpha }g_{1}(v_{2})\otimes bk)\\
&=&v_{1}\otimes \eta _{\alpha }g_{1}(v_{2})bk\\
&=&v_{1}\otimes (\eta _{\alpha }g_{1}(v_{2})b)k\\
&=&\lambda _{\alpha }(v_{1}\otimes (\eta _{\alpha
}g_{1}(v_{2})b)\otimes k)\\
&=&\lambda _{\alpha }(q_{\alpha }\otimes I_{\alpha })(v\otimes
b\otimes k).\text{ \ \ \ \ \ \ \ \ \ \ \ \ \ \ \ \ \ \ \ \ \ \ \ \
\ \ \ \ \ \ \ \ \ \ \ \ \ \ \ \ \ \ \ \ \ \ \ \ \ \ \ \ \ \ \ \ \
\ \ \ \ \ \ \ \ \ \ \ \ \ \ \ \ \ \ \ \ \ \ \ \ \ \ \ \ \ \ \ \ \
\ \ \ \ \ \ \ \ \ \ \ \ \ \ \ \ \ \ \ \ \ \ \ \ \ \ \ \ \ \ \ \ \
\ \ \ \ \ \ }
\end{eqnarray*}
\end{proof}
Throughout this last part, we consider $(C,\sigma )$ as a left $\pi -$
coisotropic quantum subgroup of $H$, $V=\{V_{\alpha }\}_{\alpha \in \pi }$
as a right $\pi -$comodule over $C$ and $G=\{G_{\alpha }\}_{\alpha \in \pi },
$where $G_{\alpha }=\{h\in H_{\alpha }:L_{1,\alpha }(h)=(\sigma _{1}\otimes
I_{\alpha }^{H})\Delta _{1,\alpha }^{H}(h)=\sigma _{1}(1)\otimes h\}.$

\begin{lem}
$Ind(\rho )$ is right $G$-module
\end{lem}

\begin{proof}
Similar to lemma 3.1.
\end{proof}

\begin{thm}
$H$ is isomorphic to $C\otimes G$ as vector space.
\end{thm}

\begin{proof}

\begin{enumerate}
\item For $h\in H_{\alpha },$ as in lemma 3.6, $h=g_{\alpha
}\sigma _{\alpha }(h_{1}^{\alpha })g_{\alpha }^{-1}\sigma _{\alpha
^{-1}}(h_{21}^{\alpha ^{-1}})h_{22}^{\alpha }.$ \item For $h\in
H_{1}$ and $\alpha \in \pi ,$ we'll prove that $\mu _{\alpha
}(g_{\alpha }^{-1}\sigma _{\alpha ^{-1}}\otimes I_{\alpha })\Delta
_{\alpha ^{-1},\alpha }^{H}(h)\in G_{\alpha }.$
\begin{eqnarray*}
&&L_{1,\alpha }\mu _{\alpha }(g_{\alpha }^{-1}\sigma _{\alpha
^{-1}}\otimes I_{\alpha })\Delta _{\alpha ^{-1},\alpha }^{H}(h)\\
&=&(\sigma _{1}\otimes I_{\alpha }^{H})\Delta _{1,\alpha }^{H}\mu
_{\alpha }(g_{\alpha }^{-1}\sigma _{\alpha ^{-1}}\otimes I_{\alpha
})\Delta _{\alpha
^{-1},\alpha }^{H}(h) \\
&=&(\sigma _{1}\otimes I_{\alpha }^{H})\mu _{H_{1}\otimes
H_{\alpha }}(\Delta _{1,\alpha }^{H}\otimes \Delta _{1,\alpha
}^{H})(g_{\alpha }^{-1}\sigma _{\alpha ^{-1}}\otimes I_{\alpha
}^{H})\Delta _{\alpha
^{-1},\alpha }^{H}(h) \\
&=&(\sigma _{1}\mu _{H_{1}}\otimes \mu _{H_{\alpha }})(I\otimes
\tau \otimes I)(\Delta _{1,\alpha }^{H}\otimes \Delta _{1,\alpha
}^{H})(g_{\alpha }^{-1}\sigma _{\alpha ^{-1}}\otimes I_{\alpha
}^{H})\Delta _{\alpha ^{-1},\alpha }^{H}(h)\\
&=&(I_{1}\otimes \mu _{H_{\alpha }})((\sigma _{1}\mu
_{H_{1}}\otimes I_{\alpha })(I\otimes \tau )\otimes I_{\alpha
}^{H})\\
&&(\Delta _{1,\alpha }^{H}g_{\alpha }^{-1}\sigma _{\alpha
^{-1}}\otimes I_{1}\otimes I_{\alpha }^{H}) (I_{\alpha
^{-1}}\otimes \Delta _{1,\alpha }^{H})\Delta _{\alpha ^{-1},\alpha
}^{H}(h) \\
&=&(I_{1}\otimes \mu _{H_{\alpha }})((\sigma _{1}\mu
_{H_{1}}\otimes I_{\alpha })(I\otimes \tau )(\Delta _{1,\alpha
}^{H}g_{\alpha }^{-1}\sigma _{\alpha ^{-1}}\otimes I_{1})\otimes
I_{\alpha }^{H})(I_{\alpha ^{-1}}\otimes \Delta _{1,\alpha
}^{H})\Delta _{\alpha
^{-1},\alpha }^{H}(h) \\
&=&(I_{1}\otimes \mu _{H_{\alpha }})((I_{1}\otimes g_{\alpha
}^{-1})(\sigma _{1}\otimes \sigma _{\alpha ^{-1}})(\mu
_{1}(S_{1}^{H}\otimes I_{1})\otimes
I_{\alpha ^{-1}}) \\
&&(I_{1}\otimes \tau )(\tau \Delta _{\alpha ^{-1},1}^{H}\otimes
I_{1})\otimes I_{\alpha }^{H})(I_{\alpha ^{-1}}\otimes \Delta
_{1,\alpha
}^{H})\Delta _{\alpha ^{-1},\alpha }^{H}(h) \\
&=&(I_{1}\otimes \mu _{H_{\alpha }})((I_{1}\otimes g_{\alpha
}^{-1})(\sigma _{1}\otimes \sigma _{\alpha ^{-1}})(\mu
_{1}(S_{1}^{H}\otimes I_{1})\otimes
I_{\alpha ^{-1}}) \\
&&(I_{1}\otimes \tau )(\tau \otimes I_{1})\otimes I_{\alpha
}^{H})(\Delta _{\alpha ^{-1},1}^{H}\otimes \Delta _{1,\alpha
}^{H})\Delta _{\alpha
^{-1},\alpha }^{H}(h) \\
&=&(I_{1}\otimes \mu _{H_{\alpha }})[(\sigma _{1}\otimes g_{\alpha
}^{-1}\sigma _{\alpha ^{-1}})(\mu _{1}(S_{1}^{H}\otimes
I_{1})\otimes
I_{\alpha ^{-1}}) \\
&&(I_{1}\otimes \tau )(\tau \otimes I_{1})\otimes I_{\alpha
}^{H}](h_{11}^{\alpha ^{-1}}\otimes h_{121}^{1}\otimes
h_{122}^{1}\otimes
h_{2}^{\alpha }) \\
&=&\sigma _{1}(S_{1}^{H}(h_{121}^{1})h_{122}^{1})\otimes g_{\alpha
}^{-1}(\sigma _{\alpha ^{-1}}(h_{11}^{\alpha ^{-1}}))h_{2}^{\alpha
} \\
&=&\sigma _{1}(\epsilon ^{H}(h_{12}^{1})1)\otimes g_{\alpha
}^{-1}(\sigma
_{\alpha ^{-1}}(h_{11}^{\alpha ^{-1}}))h_{2}^{\alpha } \\
&=&\sigma _{1}(1)\otimes g_{\alpha }^{-1}(\sigma _{\alpha
^{-1}}(\epsilon ^{H}(h_{12}^{1})h_{11}^{\alpha
^{-1}}))h_{2}^{\alpha }\\
&=&\sigma _{1}(1)\otimes g_{\alpha }^{-1}(\sigma _{\alpha
^{-1}}(h_{1}^{\alpha ^{-1}}))h_{2}^{\alpha } \\
&=&\sigma _{1}(1)\otimes \mu _{\alpha }(g_{\alpha }^{-1}\sigma
_{\alpha ^{-1}}\otimes I_{\alpha })\Delta _{\alpha ^{-1},\alpha
}^{H}(h)\text{ \ \ \ \ \ \ \ \ \ \ \ \ \ \ \ \ \ \ \ \ \ \ \ \ \ \
\ \ \ \ \ \ \ \ \ \ \ \ \ \ \ \ \ \ \ \ \ \ \ \ \ \ \ \ \ \ \ \ \
\ \ \ \ \ \ \ \ \ \ \ \ \ \ \ \ \ \ \ \ \ \ \ \ \ \ \ \ \ \ \ \ \
\ \ \ \ \ \ \ \ \ \ \ \ \ \ \ \ \ \ \ \ \ \ \ \ \ \ \ \ \ \ \ \ \
\ \ \ \ }
\end{eqnarray*}

\item We define $A_{\alpha }:C_{\alpha }\otimes G_{\alpha
}\rightarrow
H_{\alpha }$ as follow%
\begin{equation*}
A_{\alpha }(c\otimes b)=\mu _{\alpha }(g_{\alpha }\otimes
I_{\alpha })(c\otimes b)=g_{\alpha }(c)b.
\end{equation*}
\end{enumerate}

Clear $A_{\alpha }$ is linear. We define $A_{\alpha
}^{-1}:H_{\alpha }\rightarrow C_{\alpha }\otimes G_{\alpha }$ as
follow
\begin{equation*}
A_{\alpha }^{-1}(h)=(\sigma _{\alpha }\otimes \mu _{\alpha
}(g_{\alpha }^{-1}\sigma _{\alpha ^{-1}}\otimes I_{\alpha
}))(\Delta _{\alpha ,\alpha ^{-1}}^{H}\otimes I_{\alpha })\Delta
_{1,\alpha }^{H}(h).
\end{equation*}%
By step 2 and property of $\Delta$ , we have $A_{\alpha }^{-1}(h)
\in C_{\alpha }\otimes G_{\alpha }$. We will prove $A_{\alpha
}A_{\alpha }^{-1}=I_{H_{\alpha }}$ and $A_{\alpha }^{-1}A_{\alpha
}=I_{C_{\alpha }\otimes B_{\alpha }}.$ Let $h\in H_{\alpha
},c\otimes b\in C_{\alpha }\otimes G_{\alpha },$
\begin{eqnarray*}
&&A_{\alpha }A_{\alpha }^{-1}(h)\\
&=&A_{\alpha }((\sigma _{\alpha }\otimes \mu _{\alpha }(g_{\alpha
}^{-1}\sigma _{\alpha ^{-1}}\otimes I_{\alpha }))(\Delta _{\alpha
,\alpha
^{-1}}^{H}\otimes I_{\alpha })\Delta _{1,\alpha }^{H}(h)) \\
&=&A_{\alpha }((\sigma _{\alpha }\otimes \mu _{\alpha }(g_{\alpha
}^{-1}\sigma _{\alpha ^{-1}}\otimes I_{\alpha }))(I_{\alpha
}\otimes \Delta
_{\alpha ^{-1},\alpha }^{H})\Delta _{\alpha ,1}^{H}(h)) \\
&=&A_{\alpha }(\sigma _{\alpha }(h_{1}^{\alpha })\otimes g_{\alpha
}^{-1}\sigma _{\alpha ^{-1}}(h_{21}^{\alpha ^{-1}})h_{22}^{\alpha
})\\
&=&g_{\alpha }\sigma _{\alpha }(h_{1}^{\alpha })g_{\alpha
}^{-1}\sigma _{\alpha ^{-1}}(h_{21}^{\alpha ^{-1}})h_{22}^{\alpha
}\\
&=&h.\text{ \ \ \ \ \ \ \ \ \ \ \ \ \ \ \ \ \ \ \ \ \ \ \ \ \ \ \
\ \ \ \ \ \ \ \ \ \ \ \ \ \ \ \ \ \ \ \ \ \ \ \ \ \ \ \ \ \ \ \ \
\ \ \ \ \ \ \ \ \ \ \ \ \ \ \ \ \ \ \ \ \ \ \ \ \ \ \ \ \ \ \ \ \
\ \ \ \ \ \ \ \ \ \ \ \ \ \ \ \ \ \ \ \ \ \ \ \ \ \ \ \ \ \ \ \ \
\ \ \ }
\end{eqnarray*}
and
\begin{eqnarray*}
&&A_{\alpha }^{-1}A_{\alpha }(c\otimes k)\\
&=&(\sigma _{\alpha }\otimes \mu _{\alpha }(g_{\alpha }^{-1}\sigma
_{\alpha ^{-1}}\otimes I_{\alpha }))(\Delta _{\alpha ,\alpha
^{-1}}^{H}\otimes I_{\alpha })\Delta _{1,\alpha }^{H}\mu _{\alpha
}(g_{\alpha }\otimes
I_{\alpha })(c\otimes k) \\
&=&(I_{\alpha }\otimes \mu _{\alpha })(I_{\alpha }\otimes
g_{\alpha }^{-1}\otimes I_{\alpha })(\sigma _{\alpha }\otimes
\sigma _{\alpha ^{-1}}\otimes I_{\alpha }^{H}) (\Delta _{\alpha
,\alpha ^{-1}}^{H}\otimes I_{\alpha }^{H})\Delta
_{1,\alpha }^{H}\mu _{\alpha }(g_{\alpha }\otimes I_{\alpha })(c\otimes k) \\
&=&(I_{\alpha }\otimes \mu _{\alpha })(I_{\alpha }\otimes
g_{\alpha }^{-1}\otimes I_{\alpha })((\sigma _{\alpha }\otimes
\sigma _{\alpha ^{-1}})\Delta _{\alpha ,\alpha ^{-1}}^{H}\otimes
I_{\alpha }^{H})\Delta
_{1,\alpha }^{H}\mu _{\alpha } (g_{\alpha }\otimes I_{\alpha })(c\otimes k) \\
&=&(I_{\alpha }\otimes \mu _{\alpha })(I_{\alpha }\otimes
g_{\alpha }^{-1}\otimes I_{\alpha })(\Delta _{\alpha ,\alpha
^{-1}}^{C}\sigma _{1}\otimes I_{\alpha }^{H})\mu _{H_{\alpha
}\otimes H_{\alpha ^{-1}}} (\Delta _{1,\alpha }^{H}g_{\alpha
}\otimes \Delta _{1,\alpha
}^{H})(c\otimes k) \\
&=&(I_{\alpha }\otimes \mu _{\alpha })(I_{\alpha }\otimes
g_{\alpha }^{-1}\otimes I_{\alpha })(\Delta _{\alpha ,\alpha
^{-1}}^{C}\otimes I_{\alpha }^{H})(\sigma _{1}\mu _{1}\otimes \mu
_{\alpha })(I\otimes \tau \otimes I)(\Delta _{1,\alpha
}^{H}g_{\alpha }\otimes \Delta
_{1,\alpha }^{H})(c\otimes k) \\
&=&(I_{\alpha }\otimes \mu _{\alpha })(I_{\alpha }\otimes
g_{\alpha }^{-1}\otimes I_{\alpha })(\Delta _{\alpha ,\alpha
^{-1}}^{C}\otimes \mu _{\alpha })[(\sigma _{1}\mu _{1}\otimes
I_{\alpha }^{H})(I\otimes \tau ) (\Delta _{1,\alpha }^{H}g_{\alpha
}(c)\otimes k_{1}^{1}\otimes
k_{2}^{\alpha })] \\
&=&(I_{\alpha }\otimes \mu _{\alpha })(I_{\alpha }\otimes
g_{\alpha }^{-1}\otimes I_{\alpha })(\Delta _{\alpha ,\alpha
^{-1}}^{C}\otimes \mu
_{\alpha })\\
&& [(I_{1}\otimes g_{\alpha })(\sigma _{1}\otimes \sigma _{\alpha
})(\mu _{1}\otimes I_{\alpha })(I_{1}\otimes \tau )(\Delta
_{1,\alpha }^{C}\otimes
I_{1})(v\otimes k_{1}^{1})\otimes k_{2}^{\alpha }] \\
&=&(I_{\alpha }\otimes \mu _{\alpha })(I_{\alpha }\otimes
g_{\alpha }^{-1}\otimes I_{\alpha })(\Delta _{\alpha ,\alpha
^{-1}}^{C}\otimes \mu _{\alpha })[\sigma
_{1}(v_{1}^{1}k_{1}^{1})\otimes g_{\alpha }\sigma
_{\alpha }(v_{2}^{\alpha })\otimes k_{2}^{\alpha }] \\
&=&(I_{\alpha }\otimes \mu _{\alpha })(I_{\alpha }\otimes
g_{\alpha }^{-1}\otimes I_{\alpha })(\Delta _{\alpha ,\alpha
^{-1}}^{C}\otimes \mu _{\alpha }) [\omega _{1}(v_{1}^{1}\otimes
\sigma _{1}(k_{1}^{1}))\otimes g_{\alpha
}\sigma _{\alpha }(v_{2}^{\alpha })\otimes k_{2}^{\alpha }] \\
&=&(I_{\alpha }\otimes \mu _{\alpha })(I_{\alpha }\otimes
g_{\alpha }^{-1}\otimes I_{\alpha })(\Delta _{\alpha ,\alpha
^{-1}}^{C}\otimes \mu _{\alpha }) [\omega _{1}(v_{1}^{1}\otimes
\sigma _{1}(1))\otimes g_{\alpha }\sigma _{\alpha }(v_{2}^{\alpha
})\otimes k]\\
&=&(I_{\alpha }\otimes \mu _{\alpha })(I_{\alpha }\otimes
g_{\alpha }^{-1}\otimes I_{\alpha })(\Delta _{\alpha ,\alpha
^{-1}}^{C}\otimes \mu _{\alpha })[\sigma _{1}(v_{1}^{1})\otimes
g_{\alpha }\sigma _{\alpha
}(v_{2}^{\alpha })\otimes k] \\
&=&(I_{\alpha }\otimes \mu _{\alpha })[(I_{\alpha }\otimes
g_{\alpha }^{-1})\Delta _{\alpha ,\alpha ^{-1}}^{C}\sigma
_{1}(v_{1}^{1})\otimes g_{\alpha }\sigma _{\alpha }(v_{2}^{\alpha
})k]\\
&=&(I_{\alpha }\otimes \mu _{\alpha })[(I_{\alpha }\otimes
g_{\alpha }^{-1})(\sigma _{\alpha }\otimes \sigma _{\alpha
^{-1}})\Delta _{\alpha ,\alpha ^{-1}}^{H}(v_{1}^{1})\otimes
g_{\alpha }\sigma _{\alpha }(v_{2}^{\alpha })k]\\
&=&(I_{\alpha }\otimes \mu _{\alpha })[(I_{\alpha }\otimes
g_{\alpha }^{-1})(\sigma _{\alpha }(v_{1}^{\alpha })\otimes \sigma
_{\alpha ^{-1}}(v_{2}^{\alpha ^{-1}})\otimes g_{\alpha }\sigma
_{\alpha }(v_{3}^{\alpha })k] \\
&=&\sigma _{\alpha }(v_{1}^{\alpha })\otimes g_{\alpha
}^{-1}\sigma _{\alpha ^{-1}}(v_{2}^{\alpha ^{-1}})g_{\alpha
}\sigma _{\alpha }(v_{3}^{\alpha })k \\
&=&\sigma _{\alpha
}(v_{1}^{\alpha })\otimes \epsilon ^{C}\sigma _{1}(v_{2}^{1})k\\
&=&\sigma _{\alpha }(\epsilon ^{H}(v_{2}^{1})v_{1}^{\alpha
})\otimes k \\
&=&\sigma _{\alpha }(v)\otimes k \\
&=&c\otimes k.
\end{eqnarray*}
\end{proof}

\section{{\protect\large {\ Coinduced representations of Hopf group coalgebra%
}}}

In this section, we study coinduced representation from left $\pi -$%
coisotropic quantum subgroup. We restrict our attention to finite
dimensional case of Hopf $\pi -$coalgebra, i.e., dim $H_{\alpha }=n^{\alpha
} $ $\prec \infty $ for all $\alpha \in \pi .$ Let us start now from left $%
\pi -$corepresentation $\rho =\{\rho _{\alpha ,\beta }\}_{\alpha ,\beta \pi
} $ from left $\pi -$coisotropic quantum subgroup $(C,\sigma )$ on $%
V=\{V_{\alpha }\}_{\alpha \in \pi }$. We define $W=\{W_{\alpha }\}_{\alpha
\in \pi },$ where $W_{\alpha }=\{F_{\alpha }\in Hom(V_{1},H_{\alpha }):$ $\
L_{1,\alpha }F_{\alpha }=(I_{C_{1}}\otimes F_{\alpha })\rho _{1,1}\}$ and $%
L_{\alpha ,\beta }=(\sigma _{\alpha }\otimes I_{\beta })\Delta _{\alpha
,\beta }$

\begin{lem}
For $\alpha ,\beta ,\gamma \in \pi ,$ $(L_{\alpha ,\beta }\otimes I_{\gamma
})\Delta _{\alpha \beta ,\gamma }=(I_{\alpha }\otimes \Delta _{\beta ,\gamma
})L_{\alpha ,\beta \gamma } $
\end{lem}

\begin{proof}
Similar to lemma 2.4 .
\end{proof}

\begin{lem}
Suppose that $\{e_{i}^{\alpha }\}_{i=1}^{n^{\alpha }}$ and $\{g_{i}^{\alpha
}\}_{i=1}^{n^{\alpha }}$ are dual bases of $H_{\alpha }$ and $H_{\alpha
}^{\ast }$ for all $\alpha \in \pi $. For $\alpha ,\beta ,\gamma \in \pi $,
fix $F_{\alpha \beta \gamma }\in W_{\alpha \beta \gamma }$, if we define the
two maps $\xi _{1}:V_{1}\rightarrow H_{\alpha }\otimes H_{\beta }\otimes
H_{\gamma }$ and $\xi _{2}:V_{1}\rightarrow H_{\alpha }\otimes H_{\beta
}\otimes H_{\gamma }$ such that
\begin{eqnarray*}
\xi _{1}(v_{1}) &=&\sum_{i=1}^{n^{\beta \gamma }}\sim (I_{\alpha }\otimes
g_{i}^{\beta \gamma })\Delta _{\alpha ,\beta \gamma }F_{\alpha \beta \gamma
}(v_{1})\otimes \Delta _{\beta ,\gamma }(e_{i}^{\beta \gamma }) \\
\xi _{2}(v_{1}) &=&\sum_{h=1}^{n^{\beta }}\sum_{l=1}^{n^{\gamma }}[\sim
(\sim \otimes I_{K})(I_{\alpha }\otimes (g_{h}^{\beta }\otimes g_{l}^{\gamma
})\Delta _{\beta ,\gamma }]\Delta _{\alpha ,\beta \gamma }F_{\alpha \beta
\gamma }(v_{1})\otimes e_{h}^{\beta }\otimes e_{l}^{\gamma }
\end{eqnarray*}
then $\xi _{1}=\xi _{2}$.
\end{lem}

\begin{proof}
We put \;\;\;\;\; $F_{\alpha \beta \gamma }(v_{1})=h_{\alpha \beta
\gamma }=\sum_{j=1}^{n^{\alpha \beta \gamma }}\lambda
_{j}e_{j}^{\alpha \beta \gamma }$\\
$ \Delta _{\beta ,\gamma }(e_{i}^{\beta \gamma
})=\sum_{r=1}^{n^{\beta }}\sum_{s=1}^{n^{\gamma }}\eta
_{rs}^{i}e_{r}^{\beta }\otimes e_{s}^{\gamma } \;\;\mbox{and}
\;\;\Delta _{\alpha ,\beta \gamma }(e_{i}^{\alpha \beta \gamma
})=\sum_{r=1}^{n^{\alpha }}\sum_{s=1}^{n^{\beta \gamma }}\theta
_{rs}^{i}e_{r}^{\alpha }\otimes e_{s}^{\beta \gamma }$ then we
have
\begin{eqnarray*}
\Delta _{\alpha ,\beta \gamma }F_{\alpha \beta \gamma }(v_{1})
&=&\Delta _{\alpha ,\beta \gamma }(h_{\alpha \beta \gamma })
=\sum_{j=1}^{n^{\gamma }}\lambda _{j}\Delta _{\alpha ,\beta \gamma
}(e_{j}^{\alpha \beta \gamma })=\sum_{j=1}^{n^{\gamma
}}\sum_{r=1}^{n^{\alpha }}\sum_{s=1}^{n^{\beta \gamma }}\lambda
_{j}\theta _{rs}^{j}e_{r}^{\alpha }\otimes e_{s}^{\beta \gamma }
\\
\mbox{imply that}
\end{eqnarray*}
\begin{eqnarray*}
\sum_{i=1}^{n^{\beta \gamma }} &\sim &(I_{\alpha }\otimes
g_{i}^{\beta \gamma })\Delta _{\alpha ,\beta \gamma }F_{\alpha
\beta \gamma
}(v_{1})\otimes \Delta _{\beta ,\gamma }(e_{i}^{\beta \gamma }) \\
&=&\sum_{i=1}^{n^{\beta \gamma }}[\sim (I_{\alpha }\otimes
g_{i}^{\beta \gamma })\sum_{j=1}^{n^{\gamma
}}\sum_{r=1}^{n^{\alpha }}\sum_{s=1}^{n^{\beta \gamma }}\lambda
_{j}\theta _{rs}^{j}e_{r}^{\alpha }\otimes e_{s}^{\beta \gamma
}]\otimes \sum_{h=1}^{n^{\beta }}\sum_{l=1}^{n^{\gamma }}\eta
_{hl}^{i}e_{h}^{\beta }\otimes e_{l}^{\gamma }
\\
&=&\sum_{i=1}^{n^{\beta \gamma }}\sum_{j=1}^{n^{\gamma
}}\sum_{r=1}^{n^{\alpha }}\sum_{s=1}^{n^{\beta \gamma }}\lambda
_{j}\theta _{rs}^{j}g_{i}^{\beta \gamma }(e_{s}^{\beta \gamma
})e_{r}^{\alpha }\otimes \sum_{h=1}^{n^{\beta
}}\sum_{l=1}^{n^{\gamma }}\eta _{hl}^{i}e_{h}^{\beta
}\otimes e_{l}^{\gamma } \\
&=&\sum_{i=1}^{n^{\beta \gamma }}\sum_{j=1}^{n^{\gamma
}}\sum_{r=1}^{n^{\alpha }}\lambda _{j}\theta
_{ri}^{j}e_{r}^{\alpha }\otimes \sum_{h=1}^{n^{\beta
}}\sum_{l=1}^{n^{\gamma }}\eta _{hl}^{i}e_{h}^{\beta
}\otimes e_{l}^{\gamma } \\
&=&\sum_{i=1}^{n^{\beta \gamma }}\sum_{j=1}^{n^{\gamma
}}\sum_{r=1}^{n^{\alpha }}\sum_{h=1}^{n^{\beta
}}\sum_{l=1}^{n^{\gamma }}\eta _{hl}^{i}\lambda _{j}\theta
_{ri}^{j}e_{r}^{\alpha }\otimes e_{h}^{\beta }\otimes
e_{l}^{\gamma }.
\end{eqnarray*}%
\begin{eqnarray*}
\sum_{h=1}^{n^{\beta }}\sum_{l=1}^{n^{\gamma }}[ &\sim &(\sim
\otimes I_{K})(I_{\alpha }\otimes (g_{h}^{\beta }\otimes
g_{l}^{\gamma })\Delta _{\beta ,\gamma }]\Delta _{\alpha ,\beta
\gamma }F_{\alpha \beta \gamma
}(v_{1})\otimes e_{h}^{\beta }\otimes e_{l}^{\gamma } \\
&=&\sum_{h=1}^{n^{\beta }}\sum_{l=1}^{n^{\gamma }}[\sim (\sim
\otimes I_{K})(I_{\alpha }\otimes (g_{h}^{\beta }\otimes
g_{l}^{\gamma })\Delta _{\beta ,\gamma }]\lbrack
\sum_{j=1}^{n^{\gamma }}\sum_{r=1}^{n^{\alpha
}}\sum_{i=1}^{n^{\beta \gamma }}\lambda _{j}\theta
_{ri}^{j}e_{r}^{\alpha }\otimes e_{i}^{\beta \gamma }]\otimes
e_{h}^{\beta }\otimes e_{l}^{\gamma }
\\
&=&\sum_{h=1}^{n^{\beta }}\sum_{l=1}^{n^{\gamma
}}\sum_{j=1}^{n^{\gamma }}\sum_{r=1}^{n^{\alpha
}}\sum_{i=1}^{n^{\beta \gamma }}\lambda _{j}\theta _{ri}^{j}\sim
(\sim \otimes I_{K})[e_{r}^{\alpha }\otimes (g_{h}^{\beta }\otimes
g_{l}^{\gamma })\Delta _{\beta ,\gamma }(e_{i}^{\beta \gamma
})]\otimes e_{h}^{\beta }\otimes e_{l}^{\gamma } \\
&=&\sum_{h=1}^{n^{\beta }}\sum_{l=1}^{n^{\gamma
}}\sum_{j=1}^{n^{\gamma }}\sum_{r=1}^{n^{\alpha
}}\sum_{i=1}^{n^{\beta \gamma }}\sum_{q=1}^{n^{\beta
}}\sum_{s=1}^{n^{\gamma }}\lambda _{j}\theta _{ri}^{j}\eta
_{qs}^{i}g_{h}^{\beta }(e_{q}^{\beta })g_{l}^{\gamma
}(e_{s}^{\gamma
})e_{r}^{\alpha }\otimes e_{h}^{\beta }\otimes e_{l}^{\gamma } \\
&=&\sum_{h=1}^{n^{\beta }}\sum_{l=1}^{n^{\gamma
}}\sum_{j=1}^{n^{\gamma }}\sum_{r=1}^{n^{\alpha
}}\sum_{i=1}^{n^{\beta \gamma }}\eta _{hl}^{i}\lambda _{j}\theta
_{ri}^{j}e_{r}^{\alpha }\otimes e_{h}^{\beta }\otimes
e_{l}^{\gamma }
\end{eqnarray*}%
Therefore, we have \; $\sum_{i=1}^{n^{\beta \gamma }} \sim
(I_{\alpha }\otimes g_{i}^{\beta \gamma })\Delta _{\alpha ,\beta
\gamma }F_{\alpha \beta \gamma }(v_{1})\otimes \Delta _{\beta
,\gamma }(e_{i}^{\beta \gamma }) =\sum_{h=1}^{n^{\beta
}}\sum_{l=1}^{n^{\gamma }}[\sim (\sim \otimes I_{K})(I_{\alpha
}\otimes (g_{h}^{\beta }\otimes g_{l}^{\gamma })\Delta _{\beta
,\gamma }]\Delta _{\alpha ,\beta \gamma }F_{\alpha \beta \gamma
}(v_{1})\otimes e_{h}^{\beta }\otimes e_{l}^{\gamma }$, and hence
$\xi _{1}=\xi _{2}$
\end{proof}

\begin{thm}
$W=\{W_{\alpha }\}_{\alpha \in \pi }$ is right $\pi -$comodule over $H$ by $%
\Omega =\{\Omega _{\alpha ,\beta }\}_{\alpha ,\beta \in \pi }$as
\begin{equation*}
\Omega _{\alpha ,\beta }:W_{\alpha \beta }\rightarrow W_{\alpha }\otimes
H_{\beta }\text{ as }\Omega _{\alpha ,\beta }(F_{\alpha \beta })=\sim
(I_{\alpha }\otimes g^{\beta })\Delta _{\alpha ,\beta }F_{\alpha \beta
}\otimes e^{\beta }.\newline
\end{equation*}
\end{thm}

\begin{proof}
Firstly, we ${^{\prime }}$ll prove that $\sim (I_{\alpha }\otimes
g^{\beta })\Delta _{\alpha ,\beta }F_{\alpha \beta }\in W_{\alpha
}$.
\begin{eqnarray*}
&&L_{1,\alpha }[\sim (I_{\alpha }\otimes g^{\beta })\Delta
_{\alpha ,\beta }F_{\alpha \beta }] \\
&=&\sim (L_{1,\alpha }\otimes g^{\beta })\Delta _{\alpha ,\beta
}F_{\alpha \beta }\\
&=&(I_{C_{1}}\otimes \sim (I_{\alpha }\otimes g^{\beta
}))(L_{1,\alpha }\otimes I_{\beta })\Delta _{\alpha
,\beta }F_{\alpha \beta }\\
&=&(I_{C_{1}}\otimes \sim (I_{\alpha }\otimes g^{\beta
}))(I_{C_{1}}\otimes \Delta _{\alpha ,\beta })L_{1,\alpha \beta
}F_{\alpha \beta }\\
&=&(I_{C_{1}}\otimes \sim (I_{\alpha }\otimes g^{\beta
}))(I_{C_{1}}\otimes
\Delta _{\alpha ,\beta })(I_{C_{1}}\otimes F_{\alpha \beta })\rho _{1,1} \\
&=&(I_{C_{1}}\otimes \sim (I_{\alpha }\otimes g^{\beta })\Delta
_{\alpha ,\beta }F_{\alpha \beta })\rho _{1,1}\text{ \ \ \ \ \ \ \
\ \ \ \ \ \ \ \ \ \ \ \ \ \ \ \ \ \ \ \ \ \ \ \ \ \ \ \ \ \ \ \ \
\ \ \ \ \ \ \ \ \ \ \ \ \ \ \ \ \ \ \ \ \ \ \ \ \ \ \ \ \ \ \ \ \
\ \ \ \ \ \ \ \ \ \ \ \ \ \ \ \ \ \ \ \ \ \ \ \ \ \ \ \ \ \ \ \ \
\ \ \ \ \ \ \ \ \ \ \ \ \ \ \ \ \ \ \ \ \ \ \ }
\end{eqnarray*}%
Now, we will prove that $\Omega $ is coaction on $W,$ i.e., the
two diagrams \newline
\begin{tabular}{lll}
$\ \ \ \ \ W_{\alpha \beta \gamma }$ & $%
\begin{tabular}{l}
$\Omega _{\alpha \beta ,\gamma }$ \\
$\rightarrow $%
\end{tabular}%
$ & $W_{\alpha \beta }\otimes H_{\gamma }$ \\
$\Omega _{\alpha ,\beta \gamma }\downarrow $ &  & $\downarrow
(\Omega
_{\alpha ,\beta }\otimes I_{\gamma })$ \\
$W_{\alpha }\otimes H_{\beta \gamma }$ &
\begin{tabular}{l}
$\longrightarrow $ \\
$I_{W_{\alpha }}\otimes \Delta _{\beta ,\gamma }$%
\end{tabular}
& $W_{\alpha }\otimes H_{\beta }\otimes H_{\gamma }$%
\end{tabular}%
, \ \ \
\begin{tabular}{lll}
$W_{\alpha }$ & $%
\begin{tabular}{l}
$\Omega _{\alpha ,1}$ \\
$\rightarrow $%
\end{tabular}%
$ & $W_{\alpha }\otimes H_{1}$ \\
& $\sim $ $\searrow $ & $\downarrow (I_{W_{\alpha }}\otimes \epsilon )$ \\
&  & $W_{\alpha }\otimes K$%
\end{tabular}
are commutes.
\begin{eqnarray*}
&&(\Omega _{\alpha ,\beta }\otimes I_{\gamma })\Omega _{\alpha
\beta ,\gamma }F_{\alpha \beta \gamma } \\
&=&(\Omega _{\alpha ,\beta }\otimes I_{\gamma })(\sim (I_{\alpha
\beta }\otimes g^{\gamma })\Delta _{\alpha \beta ,\gamma
}F_{\alpha \beta \gamma
}\otimes e^{\gamma }) \\
&=&\Omega _{\alpha ,\beta }(\sim (I_{\alpha \beta }\otimes
g^{\gamma })\Delta _{\alpha \beta ,\gamma }F_{\alpha \beta \gamma
})\otimes e^{\gamma }
\\
&=&\sim (I_{\alpha }\otimes g^{\beta })\Delta _{\alpha ,\beta
}\sim (I_{\alpha \beta }\otimes g^{\gamma })\Delta _{\alpha \beta
,\gamma
}F_{\alpha \beta \gamma }\otimes e^{\beta }\otimes e^{\gamma } \\
&=&\sim (\sim \otimes I_{K})((I_{\alpha }\otimes g^{\beta })\Delta
_{\alpha ,\beta }\otimes g^{\gamma })\Delta _{\alpha \beta ,\gamma
}F_{\alpha \beta
\gamma }\otimes e^{\beta }\otimes e^{\gamma } \\
&=&\sim (\sim \otimes I_{K})(I_{\alpha }\otimes g^{\beta }\otimes
g^{\gamma })(\Delta _{\alpha ,\beta }\otimes I_{\gamma })\Delta
_{\alpha \beta ,\gamma
}F_{\alpha \beta \gamma }\otimes e^{\beta }\otimes e^{\gamma } \\
&=&\sim (\sim \otimes I_{K})(I_{\alpha }\otimes g^{\beta }\otimes
g^{\gamma })(I_{\alpha }\otimes \Delta _{\beta ,\gamma })\Delta
_{\alpha ,\beta \gamma
}F_{\alpha \beta \gamma }\otimes e^{\beta }\otimes e^{\gamma } \\
&=&\sim (\sim \otimes I_{K})(I_{\alpha }\otimes (g^{\beta }\otimes
g^{\gamma })\Delta _{\beta ,\gamma })\Delta _{\alpha ,\beta \gamma
}F_{\alpha \beta \gamma }\otimes e^{\beta }\otimes e^{\gamma }\\
&=&\xi_1 \text{ \ \ \ \ \ \ \ \ \ \ \ \ \ \ \ \ \ \ \ \ \ \ \ \ \
\ \ \ \ \ \ \ \ \ \ \ \ \ \ \ \ \ \ \ \ \ \ \ \ \ \ \ \ \ \ \ \ \
\ \ \ \ \ \ \ \ \ \ \ \ \ \ \ \ \ \ \ \ \ \ \ \ \ \ \ \ \ \ \ \ \
\ \ \ \ \ \ \ \ \ \ \ \ \ \ \ \ \ \ \ \ \ \ \ \ \ \ \ \ \ \ \ \ \
\ \ \ \ \ }
\end{eqnarray*}%
Also
\begin{eqnarray*}
&&(I_{W_{\alpha }}\otimes \Delta _{\beta ,\gamma })\Omega _{\alpha
,\beta \gamma }F_{\alpha \beta \gamma }\\
&=&(I_{W_{\alpha }}\otimes \Delta _{\beta ,\gamma })(\sim
(I_{\alpha }\otimes g^{\beta \gamma })\Delta _{\alpha ,\beta
\gamma }F_{\alpha \beta
\gamma }\otimes e^{\beta \gamma }) \\
&=&\sim (I_{W_{\alpha }}\otimes g^{\beta \gamma })\Delta _{\alpha
,\beta \gamma }F_{\alpha \beta \gamma }\otimes \Delta _{\beta
,\gamma }(e^{\beta \gamma })\\
&=&\xi_2.\text{ \ \ \ \ \ \ \ \ \ \ \ \ \ \ \ \ \ \ \ \ \ \ \ \ \
\ \ \ \ \ \ \ \ \ \ \ \ \ \ \ \ \ \ \ \ \ \ \ \ \ \ \ \ \ \ \ \ \
\ \ \ \ \ \ \ \ \ \ \ \ \ \ \ \ \ \ \ \ \ \ \ \ \ \ \ \ \ \ \ \ \
\ \ \ \ \ \ \ \ \ \ \ \ \ \ \ \ \ \ \ \ \ \ \ \ \ \ \ \ \ \ \ \ \
\ \ \ \ \ }
\end{eqnarray*}%
from lemma 4.2. the first diagram is commute. For the second
diagram
\begin{eqnarray*}
&&(I_{W_{\alpha }}\otimes \epsilon )\Omega _{\alpha ,1}F_{\alpha
}\\
&=&(I_{W_{\alpha }}\otimes \epsilon )[\sim (I_{\alpha }\otimes
g^{1})\Delta _{\alpha ,1}F_{\alpha }\otimes e^{1}] \\
&=&\sim (I_{W_{\alpha }}\otimes g^{1})\Delta _{\alpha ,1}F_{\alpha
}\otimes
\epsilon (e^{1}) \\
&=&\sim (I_{W_{\alpha }}\otimes \epsilon (e^{1})g^{1})\Delta
_{\alpha ,1}F_{\alpha }\otimes 1_{K}\\
&=&\sim (I_{W_{\alpha }}\otimes \epsilon )\Delta _{\alpha
,1}F_{\alpha }\otimes 1_{K} \\
&=&F_{\alpha }\otimes 1_{K}.\text{ \ \ \ \ \ \ \ \ \ \ \ \ \ \ \ \
\ \ \ \ \ \ \ \ \ \ \ \ \ \ \ \ \ \ \ \ \ \ \ \ \ \ \ \ \ \ \ \ \
\ \ \ \ \ \ \ \ \ \ \ \ \ \ \ \ \ \ \ \ \ \ \ \ \ \ \ \ \ \ \ \ \
\ \ \ \ \ \ \ \ \ \ \ \ \ \ \ \ \ \ \ \ \ \ \ \ \ \ \ \ \ \ \ \ \
\ \ \ \ \ \ \ \ \ \ \ \ \ \ }
\end{eqnarray*}
\end{proof}

\begin{rem}
Given a right corepresentation $\rho $ of the $\pi -$coisotropic quantum
subgroup of $(C,\sigma )$ the corresponding corepresentation $(I\otimes
\Delta)$ on $Coind(\rho )$ of $H$ is called coinduced representation from $%
\rho$ on $H$.\\
\end{rem}

Now, we constract an induced and coinduced representation from subHopf $\pi
- $coalgebra.

\begin{thm}
Let $H$ be a finite dimentional Hopf $\pi $-coalgebra and $A=\{A_{\alpha
}\}_{\alpha \in \pi }$ be an isolated subHopf $\pi $-coalgebra of $H.$ If $%
V=\{V_{\alpha }\}_{\alpha \in \pi }$ is left $\pi $-comodule over $A$ by $%
\rho =\{\rho _{\alpha ,\beta }\}_{\alpha ,\beta \in \pi },$ then we can
construct an induced and coinduced representation over $H $.
\end{thm}

\begin{proof}
 Since $A=\{A_{\alpha }\}_{\alpha \in \pi }$ is an isolated
subHopf $\pi $-coalgebra of $H$, there exist a family
$I=\{I_{\alpha }\}_{\alpha \in \pi }$ of Hopf $\pi $-coideal of
$H$ such that $H_{\alpha }=A_{\alpha }\oplus I_{\alpha }$ for all
$\alpha \in \pi .$ We will prove
that $(A,\sigma )$ is left $\pi -$coisotropic quantum subgroup of $H$ where $%
\sigma =\{\sigma _{\alpha }\}_{\alpha \in \pi }$ and $\sigma
_{\alpha }:H_{\alpha }\rightarrow A_{\alpha }$ as $\sigma _{\alpha
}(m_{\alpha }+i_{\alpha })=m_{\alpha }.$ Clear, $A$ is $\pi
$-coalgebra. We ${^{\prime }}$ll prove that $A_{\alpha }$ is left
$H_{\alpha }$-module for all $\alpha \in \pi .$ We define $\Phi
_{\alpha }:H_{\alpha }\otimes A_{\alpha }\rightarrow A_{\alpha }$ as follow $%
\Phi _{\alpha }((m_{\alpha }+i_{\alpha })\otimes a_{\alpha
})=m_{\alpha }a_{\alpha }$. Then $\Phi _{\alpha }(\mu _{\alpha
}\otimes I_{\alpha })((m_{\alpha }+i_{\alpha })\otimes (n_{\alpha
}+j_{\alpha })\otimes a_{\alpha }) =\Phi _{\alpha }((m_{\alpha
}n_{\alpha }+m_{\alpha }j_{\alpha }+i_{\alpha }n_{\alpha
}+i_{\alpha }j_{\alpha })\otimes a_{\alpha })=m_{\alpha }n_{\alpha
}a_{\alpha }$ and $\Phi _{\alpha }(I_{\alpha }\otimes \Phi
_{\alpha })((m_{\alpha }+i_{\alpha })\otimes (n_{\alpha
}+j_{\alpha })\otimes a_{\alpha }) =\Phi _{\alpha }((m_{\alpha
}+i_{\alpha })\otimes n_{\alpha }a_{\alpha })=m_{\alpha }n_{\alpha
}a_{\alpha }$ where $m_{\alpha }j_{\alpha }+i_{\alpha }n_{\alpha
}+i_{\alpha }j_{\alpha }\in I_{\alpha }.$ Also, we have $\Phi
_{\alpha }(\eta _{\alpha }\otimes I_{\alpha })(k\otimes a_{\alpha
}) =\Phi _{\alpha }((\eta _{\alpha }(k)+0)\otimes a_{\alpha })
=\eta _{\alpha }(k)a_{\alpha }=ka_{\alpha }.$ Now, we will prove
that $\sigma _{\alpha }$ is $\pi -$coalgebra map, $(\sigma
_{\alpha }\otimes \sigma _{\beta })\Delta _{\alpha ,\beta
}(h_{\alpha \beta })=(\sigma _{\alpha }\otimes \sigma _{\beta
})\Delta _{\alpha ,\beta }(m_{\alpha \beta }+i_{\alpha \beta })
=(\sigma _{\alpha }\otimes \sigma _{\beta })(\Delta _{\alpha
,\beta }(m_{\alpha \beta })+\Delta _{\alpha ,\beta }(i_{\alpha
\beta })) =\Delta _{\alpha ,\beta }(m_{\alpha \beta }) $ and
$\Delta _{\alpha ,\beta }\sigma _{\alpha \beta }(h_{\alpha \beta
}) =\Delta _{\alpha ,\beta }\sigma _{\alpha \beta }(m_{\alpha
\beta }+i_{\alpha \beta }) =\Delta _{\alpha ,\beta }(m_{\alpha
\beta }). $ We will prove that $\sigma _{\alpha }$ is left module
map, $\Phi _{\alpha }(I_{\alpha }\otimes \sigma _{\alpha
})(h_{\alpha }\otimes k_{\alpha }) =\Phi _{\alpha }(I_{\alpha
}\otimes \sigma _{\alpha })((m_{\alpha }+i_{\alpha })\otimes
(n_{\alpha }+j_{\alpha })) =\Phi _{\alpha }((m_{\alpha }+i_{\alpha
})\otimes n_{\alpha })=m_{\alpha }n_{\alpha } $ and $\sigma
_{\alpha }\mu _{\alpha }(h_{\alpha }\otimes k_{\alpha }) =\sigma
_{\alpha }(m_{\alpha }n_{\alpha }+m_{\alpha }j_{\alpha }+i_{\alpha
}n_{\alpha }+i_{\alpha }j_{\alpha })=m_{\alpha }n_{\alpha } $
Therefore $(A,\sigma )$ is left $\pi -$coisotropic quantum
subgroup of $H.$ Since $V$ is left $\pi -$comodule over $A$, then
we can construct an induced and coinduced representation over $H$,
by theorem 2.6 and theorem 4.3.
\end{proof}

\section{{\protect\large {\ Simplicity of induced representation of Hopf
group coalgebra}}}

In this section we have studied the simplicity of our induced
representations of Hopf $\pi -$coalgebra. Let $H$ be a Hopf $\pi -$%
coalgebra, $(C,\sigma )$ a $\pi -$coisotropic quantum subgroup. If we have
two equivalent right $\pi -$comodules over $C$. What is the relations
between the two induced and the two coinduced representations over $H$ ?
Also, we have the following question: what is the algebraic relations of the
induction with respect to direct sums of representations in point of view of
simplicity. Finally, we have proved that if the induced representation is
simple, then its algebraic induction is constructed from simple
representation of the group coisotropic quantum subgroup. If we have simple
representation of group coisotropic quantum subgroup for Hopf $\pi -$%
coalgebra, then the induced representation is not necessary to be simple.%
\newline

A right $\pi $-comodule $M=\{M_{\alpha }\}_{\alpha \in \pi }$ is said to be
simple if it is non-zero (i.e., $M_{\alpha }\neq 0$ for some $\alpha \in \pi$%
) and if it has no $\pi-$subcomodules other than $0={0}_{\alpha \in \pi}$
and itself.\newline

A right $\pi-$comodules $(V, \rho)$ and $(W, \rho^{\prime})$ over $C$ are
equivalent if there exist a family $F=\{ F_{\alpha}:V_{\alpha}\rightarrow
W_{\alpha} \}_{\alpha \in \pi}$ of vector space isomorphisms such that
\begin{equation*}
\rho^{\prime}_{\alpha,\beta}F_{\alpha\beta}=(F_{\alpha}\otimes
I_{\beta})\rho_{\alpha,\beta} \;\;\;\;\; \text{for all } \alpha, \beta \in
\pi
\end{equation*}

\begin{thm}
Let $(V,\rho )$ and $(W,\rho ^{\prime })$ be a right $\pi -$comodules over a
$\pi -$coisotropic quantum subgroup $(C,\sigma ).$ If there exist a vector
space isomorphism $F_{1}:$ $V_{1}\rightarrow W_{1}$ such that
\begin{equation*}
\text{ \ \ \ \ \ \ \ \ \ \ \ \ \ \ \ \ \ \ \ \ \ \ \ \ \ \ \ \ }\psi
_{1,1}F_{1}=(F_{1}\otimes I_{1})\rho _{1,1}\;\ \ \ \ \ \ \ \ \ \ \ \ \ \ \ \
\ \ \ \ \ \ \ \ \ \ \ \ \ \ \ \ \ \ \ \ \ \ \ \ \ \ \ \ \ \ \ \ \ \ \ \ \ \
\ (\#)\ \ \
\end{equation*}
then the induced representations $Ind(\rho )$ and $Ind(\rho ^{\prime })$ are
equivalent right $\pi -$comodules over $H$.
\end{thm}

\begin{proof} Let us consider the family of vector space
isomorphisms
\begin{equation*}
\overline{F}=\{\overline{F}_{\alpha }=F_{1}\otimes I_{\alpha
}:V_{1}\otimes H_{\alpha }\rightarrow W_{1}\otimes H_{\alpha
}\}_{\alpha \in \pi }.\newline
\end{equation*}

Firstly, we must prove that, for $\alpha \in \pi $,
$\overline{F}_{\alpha }(Ind(\rho )_{\alpha })\subseteq Ind(\rho
^{\prime })_{\alpha }$. Let $f\in Ind(\rho )_{\alpha },$ then
\begin{equation*}
\text{ \ \ \ \ \ \ \ \ \ \ \ \ \ \ \ \ \ \ \ \ }(I_{1}\otimes
L_{1,\alpha })(f)=(\rho _{1,1}\otimes I_{\alpha
})(f)\;\;\;\;\;\;\;\;\;\;\;\;\;\;\ \ \ \ \ \ \ \ \ \ \
\;\;\;\;\;\;\;\;\;\;\;\;\;\;\ \ \;\;(1)
\end{equation*}

We want to prove $\overline{F}_{\alpha }(f)\in Ind(\rho ^{\prime
})_{\alpha },$ i. e., $\ (I_{1}\otimes L_{1,\alpha
})\overline{F}_{\alpha }(f)=(\rho _{1,1}^{\prime }\otimes
I_{\alpha })\overline{F}_{\alpha }(f)$
\begin{eqnarray*}
&&(I_{1}\otimes L_{1,\alpha })\overline{F}_{\alpha }(f) \\
&=&(I_{1}\otimes L_{1,\alpha })(F_{1}\otimes I_{\alpha })(f) \\
&=&(F_{1}\otimes I_{1}\otimes I_{\alpha })(I_{1}\otimes L_{1,\alpha })(f) \\
&=&(F_{1}\otimes I_{1}\otimes I_{\alpha })(\rho _{1,1}\otimes
I_{\alpha
})(f)\;\;\;\;\;\;\;\;\;\;\;\;\;\;\;\;\;\;\;\;\;\;\;\;\;\;\ \ \ \ \
\ \ \ \ \
\ \ \ \ \ \ \ \ \ \ \ \ \ \ \ \ \ \ \ \ \ \ \ \ \ \ \ \ \text{by }(1) \\
&=&[(F_{1}\otimes I_{1})\rho _{1,1}\otimes I_{\alpha }](f) \\
&=&(\rho _{1,1}^{\prime }F_{1}\otimes I_{\alpha
})(f)\;\;\;\;\;\;\;\;\;\;\;\;\;\;\;\;\;\;\;\;\;\;\;\;\;\;\;\;\;\;\;\;\;\;\;%
\;\;\;\ \ \ \ \ \ \ \ \ \ \ \ \ \ \ \ \ \ \ \ \ \ \ \ \ \ \ \ \ \
\ \ \ \ \
\ \ \ \ \ \text{by }(\#) \\
&=&(\rho _{1,1}^{\prime }\otimes I_{\alpha })(F_{1}\otimes I_{\alpha })(f) \\
&=&(\rho _{1,1}^{\prime }\otimes I_{\alpha })\overline{F}_{\alpha
}(f)\text{ \ \ \ \ \ \ \ \ \ \ \ \ \ \ \ \ \ \ \ \ \ \ \ \ \ \ \ \
\ \ \ \ \ \ \ \ \ \ \ \ \ \ \ \ \ \ \ \ \ \ \ \ \ \ \ \ \ \ \ \ \
\ \ \ \ \ \ \ \ \ \ \ \ \ \ \ \ \ \ \ \ \ \ \ \ \ \ \ \ \ \ \ \ \
\ \ \ \ \ \ \ \ \ \ \ \ \ \ \ \ \ \ \ \ \ \ \ \ \ \ \ \ \ \ \ \ \
\ \ \ \ \ \ \ \ \ \ \ \ \ \ \ \ \ \ \ \ \ \ \ \ \ \ \ \ \ \ \ \ \
\ \ \ \ \ \ \ \ \ \ \ \ \ \ \ \ \ \ \ \ \ \ \ \ \ \ \ \ \ \ \ \ \
\ \ \ \ \ \ \ \ \ \ \ \ \ \ \ \ \ \ \ \ \ \ \ \ \ \ \ \ \ \ \ \ \
\ \ \ \ \ \ \ \ \ \ \ \ \ \ \ \ \ \ \ \ \ \ \ \ \ \ \ \ \ \ \ \ \
\ \ \ \ \ \ \ \ \ }
\end{eqnarray*}

Now we shall prove that $\overline{F}$ is $\pi -$comodule map,
which follows from, for $f\in Ind(\rho )_{\alpha \beta }$,
\begin{eqnarray*}
&&(\overline{F}_{\alpha }\otimes I_{\beta })(I_{1}\otimes \Delta
_{\alpha
,\beta })(f) \\
&=&(F_{1}\otimes I_{\alpha }\otimes I_{\beta })(I_{1}\otimes
\Delta _{\alpha
,\beta })(f) \\
&=&(I_{1}\otimes \Delta _{\alpha ,\beta })\overline{F}_{\alpha
}(f).\text{ \ \ \ \ \ \ \ \ \ \ \ \ \ \ \ \ \ \ \ \ \ \ \ \ \ \ \
\ \ \ \ \ \ \ \ \ \ \ \ \ \ \ \ \ \ \ \ \ \ \ \ \ \ \ \ \ \ \ \ \
\ \ \ \ \ \ \ \ \ \ \ \ \ \ \ \ \ \ \ \ \ \ \ \ \ \ \ \ \ \ \ \ \
\ \ \ \ \ \ \ \ \ \ \ \ \ \ \ \ \ \ \ \ \ \ \ \ \ \ \ \ \ \ \ \ \
\ \ \ \ \ \ \ \ \ \ \ \ \ \ \ \ \ \ \ \ \ \ \ \ \ \ \ \ \ \ \ \ \
\ \ \ \ \ \ \ \ \ \ \ \ \ \ \ \ \ \ \ \ \ \ \ \ \ \ \ \ \ \ \ \ \
\ \ \ \ \ \ \ \ \ \ \ \ \ \ \ \ \ \ \ \ \ \ \ \ \ \ \ \ \ \ \ \ \
\ \ \ \ \ \ \ \ \ \ \ \ \ \ \ \ \ \ \ \ \ \ \ \ \ \ \ \ \ \ \ \ \
\ \ \ \ \ \ \ \ \ }
\end{eqnarray*}
\end{proof}

\begin{cor}
If $(V, \rho)$ and $(W, \rho^{\prime})$ are equivalent right $\pi-$comodules
over $C$, then the induced representations $Ind(\rho)$ and $%
Ind(\rho^{\prime})$ are equivalent right $\pi-$comodules over $H$.
\end{cor}

\begin{lem}
Let $(V,\rho )$ and $(W,\psi )$ be two right $\pi -$comodules over $C$, then
$V\oplus W=\{(V\oplus W)_{\alpha }=V_{\alpha }\oplus W_{\alpha }\}_{\alpha
\in \pi }$ is right $\pi -$comodule over $C$ by $\rho \oplus \psi =\{(\rho
\oplus \psi )_{\alpha ,\beta }=\rho _{\alpha ,\beta }\oplus \psi _{\alpha
,\beta }:V_{\alpha \beta }\oplus W_{\alpha \beta }\rightarrow (V_{\alpha
}\oplus W_{\alpha })\otimes C_{\alpha }\}_{\alpha ,\beta \in \pi }$ where
\begin{equation*}
\rho _{\alpha ,\beta }\oplus \psi _{\alpha ,\beta }=(l_{\alpha }^{V}\otimes
I_{\beta })\rho _{\alpha ,\beta }P_{\alpha \beta }^{V}+(l_{\alpha
}^{W}\otimes I_{\beta })\psi _{\alpha ,\beta }P_{\alpha \beta }^{W},\;\;%
\text{where}
\end{equation*}%
\begin{equation*}
P_{\alpha }^{V}:V_{\alpha }\oplus W_{\alpha }\rightarrow V_{\alpha }
\end{equation*}%
\begin{equation*}
P_{\alpha }^{W}:V_{\alpha }\oplus W_{\alpha }\rightarrow W_{\alpha }
\end{equation*}%
are the natural projection maps and
\begin{equation*}
l_{\alpha }^{V}:V_{\alpha }\rightarrow V_{\alpha }\oplus W_{\alpha }
\end{equation*}%
\begin{equation*}
l_{\alpha }^{W}:W_{\alpha }\rightarrow V_{\alpha }\oplus W_{\alpha }
\end{equation*}%
are the natural immersion maps.\newline
\end{lem}

\begin{thm}
Let $(V,\rho )$ and $(W,\psi )$ be two right $\pi -$comodules over $C$, then
\begin{equation*}
Ind(\rho \oplus \psi )\cong Ind(\rho )\oplus Ind(\psi ).
\end{equation*}
\end{thm}

\begin{proof}
We define $\theta :Ind(\rho \oplus \psi )\rightarrow Ind(\rho
)\oplus Ind(\psi )$ as follows: $\theta =\{\theta _{\alpha
}=(P_{1}^{V}\otimes I_{\alpha })+(P_{1}^{W}\otimes I_{\alpha
})\}_{\alpha
\in \pi }$. $\theta _{\alpha }$ is an isomorphism of vector spaces and $%
\theta _{\alpha }:Ind(\rho \oplus \psi )_{\alpha }\rightarrow
Ind(\rho )_{\alpha }\oplus Ind(\psi )_{\alpha }$. We will prove
that, for $f\in Ind(\rho \oplus \psi )_{\alpha }$,

\begin{enumerate}
\item $(P_{1}^{V}\otimes I_{\alpha})(f) \in Ind(\rho)_{\alpha}$

\item $(P_{1}^{W}\otimes I_{\alpha})(f) \in Ind(\psi)_{\alpha}$
\end{enumerate}

and hence $\theta _{\alpha }(f)\in Ind(\rho )_{\alpha }\oplus
Ind(\psi )_{\alpha }$.
\begin{eqnarray*}
&&(I_{1}\otimes L_{1,\alpha })(P_{1}^{V}\otimes I_{\alpha })(f) \\
&=&(P_{1}^{V}\otimes I_{1}\otimes I_{\alpha })(I_{1}\otimes
L_{1,\alpha })(f)
\\
&=&(P_{1}^{V}\otimes I_{1}\otimes I_{\alpha })((\rho \oplus \psi
)_{1,1}\otimes I_{\alpha })(f) \\
&=&(P_{1}^{V}\otimes I_{1}\otimes I_{\alpha })[[(l_{\alpha
}^{V}\otimes I)\rho _{1,1}P_{1}^{V}+(l_{1}^{W}\otimes I)\psi
_{1,1}P_{1}^{W}]\otimes
I_{\alpha }](f) \\
&=&(P_{1}^{V}\otimes I_{1}\otimes I_{\alpha })[(l_{\alpha
}^{V}\otimes I)\rho _{1,1}P_{1}^{V}\otimes I_{\alpha
}](f)+(P_{1}^{V}\otimes I_{1}\otimes I_{\alpha })[(l_{\alpha
}^{W}\otimes I)\psi _{1,1}P_{1}^{W}\otimes I_{\alpha
}](f) \\
&=&(P_{1}^{V}\otimes I_{1}\otimes I_{\alpha })[(l_{\alpha
}^{V}\otimes
I)\rho _{1,1}P_{1}^{V}\otimes I_{\alpha }](f) \\
&=&[(P_{1}^{V}l_{\alpha }^{V}\otimes I)\rho _{1,1}P_{1}^{V}\otimes
I_{\alpha
}](f) \\
&=&(\rho _{1,1}P_{1}^{V}\otimes I_{\alpha })(f) \\
&=&(\rho _{1,1}\otimes I_{\alpha })(P_{1}^{V}\otimes I_{\alpha
})(f)
\end{eqnarray*}%
then $(P_{1}^{V}\otimes I_{\alpha })(f)\in Ind(\rho )_{\alpha
}$.\newline

Similar $(P_{1}^{W}\otimes I_{\alpha })(f)\in Ind(\psi )_{\alpha
}$. Therefore $\theta _{\alpha }(f)\in Ind(\rho )_{\alpha }\oplus
Ind(\psi )_{\alpha }$.

Clear $\theta _{\alpha }$ is an isomorphism. The intertwining
property is
the commutativity of the following diagram%
\begin{equation*}
\begin{array}{ccc}
Ind(\rho \oplus \psi )_{\alpha \beta } &
\begin{array}{c}
\theta _{\alpha \beta } \\
\longrightarrow%
\end{array}
& Ind(\rho )_{\alpha \beta }\oplus Ind(\psi )_{\alpha \beta } \\
I\otimes \Delta _{\alpha ,\beta }\downarrow &  & \downarrow \overline{\rho }%
_{\alpha ,\beta }\oplus \overline{\psi }_{\alpha ,\beta } \\
Ind(\rho \oplus \psi )_{\alpha }\otimes H_{\beta } &
\begin{array}{c}
\longrightarrow \\
\theta _{\alpha }\otimes I_{\beta }%
\end{array}
& (Ind(\rho )_{\alpha }\oplus Ind(\psi )_{\alpha })\otimes H_{\beta }%
\end{array}%
\end{equation*}%
where $Ind(\rho \oplus \psi )$ is right $\pi -$comodule over $H$ by $%
I\otimes \Delta =\{I_{V_{1}\oplus W_{1}}\otimes \Delta _{\alpha
,\beta }\}$. Also, $Ind(\rho )\oplus Ind(\psi )$ is right $\pi
-$comodule over $H$ by
\begin{equation*}
\overline{\rho }\oplus \overline{\psi }=\{\overline{\rho }_{\alpha
,\beta }\oplus \overline{\psi }_{\alpha ,\beta }\}
\end{equation*}%
where
\begin{equation*}
\overline{\rho }_{\alpha ,\beta }\oplus \overline{\psi }_{\alpha
,\beta }=(l_{\alpha }^{Ind(\rho )}\otimes I_{\beta })(I_{1}\otimes
\Delta _{\alpha ,\beta })P_{\alpha \beta }^{Ind(\rho )}+(l_{\alpha
}^{Ind(\psi )}\otimes I_{\beta })(I_{1}\otimes \Delta _{\alpha
,\beta })P_{\alpha \beta }^{Ind(\psi )}
\end{equation*}%
Now,
\begin{eqnarray*}
&&[\overline{\rho }_{\alpha ,\beta }\oplus \overline{\psi
}_{\alpha ,\beta
}]\vartheta _{\alpha \beta }(f) \\
&=&[(l_{\alpha }^{Ind(\rho )}\otimes I_{\beta })(I_{1}\otimes
\Delta _{\alpha ,\beta })P_{\alpha \beta }^{Ind(\rho )}+(l_{\alpha
}^{Ind(\psi )}\otimes I_{\beta })(I_{1}\otimes \Delta _{\alpha
,\beta })P_{\alpha \beta
}^{Ind(\psi )}] \\
&&[(P_{1}^{V}\otimes I_{\alpha \beta })+(P_{1}^{W}\otimes
I_{\alpha \beta
})][(v+w)\otimes h] \\
&=&[(l_{\alpha }^{Ind(\rho )}\otimes I_{\beta })(I_{1}\otimes
\Delta _{\alpha ,\beta })P_{\alpha \beta }^{Ind(\rho )}+(l_{\alpha
}^{Ind(\psi )}\otimes I_{\beta })(I_{1}\otimes \Delta _{\alpha
,\beta })P_{\alpha \beta
}^{Ind(\psi )}][(v\otimes h)+(w\otimes h)] \\
&=&[(l_{\alpha }^{Ind(\rho )}\otimes I_{\beta })(I_{1}\otimes
\Delta _{\alpha ,\beta })P_{\alpha \beta }^{Ind(\rho )}][(v\otimes
h)+(w\otimes h)]+[(l_{\alpha }^{Ind(\psi )}\otimes I_{\beta
})(I_{1}\otimes \Delta
_{\alpha ,\beta })P_{\alpha \beta }^{Ind(\psi )}] \\
&&[(v\otimes h)+(w\otimes h)] \\
&=&[(l_{\alpha }^{Ind(\rho )}\otimes I_{\beta })(I_{1}\otimes
\Delta _{\alpha ,\beta })](v\otimes h)+[(l_{\alpha }^{Ind(\psi
)}\otimes I_{\beta
})(I_{1}\otimes \Delta _{\alpha ,\beta })](w\otimes h) \\
&=&(v\otimes h_{1}+0)\otimes h_{2}+(0+w\otimes h_{1})\otimes h_{2} \\
&=&[(v\otimes h_{1})+(w\otimes h_{1})]\otimes h_{2} \\
&=&[[(P_{1}^{V}\otimes I_{\alpha })+(P_{1}^{W}\otimes I_{\alpha
})]\otimes
I_{\beta }][(v+w)\otimes \Delta _{\alpha ,\beta }(h)] \\
&=&(\theta _{\alpha }\otimes I_{\beta })(I\otimes \Delta _{\alpha
,\beta })(f)
\end{eqnarray*}%
i.e., the above diagram commutes.\newline
\end{proof}

\begin{cor}
Let $H$ be a Hopf $\pi -$coalgebra, $(C,\sigma )$ a $\pi -$coisotropic
quantum subgroup of $H$ and $(V,\rho )$ a right $\pi -$comodules over $C$.
If $Ind(\rho )$ is simple, then $\rho $ is simple.\newline
\end{cor}

The inverse of the above result is not true in general case by the following
example.\newline

\textbf{Example } Let $A$ be the algebra generated by $1$ and the elements $%
x,$ $y$, $x^{-1}$, $y^{-1}$ and $g$ with \
\begin{eqnarray*}
xx^{-1} &=&x^{-1}x\;=\;yy^{-1}=y^{-1}y=1,\text{ \ }xy=yx, \\
xg &=&qgx,\;\;\;yg=qgy,\;\;x^{-1}g=q^{-1}gx^{-1},\;\ \ \
y^{-1}g=q^{-1}gy^{-1}\text{ \ \ \ \ \ \ \ \ \ \ \ \ \ \ \ \ \ \ \ \ \ \ \ \
\ \ \ \ \ \ \ \ \ \ \ \ \ \ \ \ \ \ \ \ \ \ \ \ \ \ \ \ \ \ \ \ \ \ \ \ \ \ }
\end{eqnarray*}%
for some fixed non zero \ $q\in K,$ then $A$ become Hopf algebra with%
\begin{eqnarray*}
\Delta \left(
\begin{array}{l}
g \\
x^{\pm 1} \\
y^{\pm 1}%
\end{array}%
\right) &=&\left(
\begin{array}{l}
g\otimes x+y\otimes g \\
x^{\pm 1}\otimes x^{\pm 1} \\
y^{\pm 1}\otimes y^{\pm 1}%
\end{array}%
\right) \text{ \ \ \ \ \ },\text{ \ \ \ \ \ }\epsilon \left(
\begin{array}{l}
g \\
x^{\pm 1} \\
y^{\pm 1}%
\end{array}%
\right) =\left(
\begin{array}{l}
0 \\
1 \\
1%
\end{array}%
\right) , \\
\text{ \ \ \ }S\left(
\begin{array}{l}
g \\
x^{\pm 1} \\
y^{\pm 1}%
\end{array}%
\right) &=&\left(
\begin{array}{l}
-y^{-1}gx^{-1} \\
x^{\mp 1} \\
y^{\mp 1}%
\end{array}%
\right)
\end{eqnarray*}

Now, we take $\pi $ the set of all group like elements of $A$, then $\pi $
acts on $A$ by Hopf algebra endomorphism as $\lambda :\pi \rightarrow End(A)$
where $\lambda (\alpha )(\beta )=\alpha \beta \alpha ^{-1}$ for all $\alpha
,\beta \in \pi $. Then, by [T], we can construct Hopf $\pi -$coalgebra as $%
H=\{H_{\alpha }\}_{\alpha \in \pi }$ where for each $\alpha \in
\pi $, the algebra $H_{\alpha }$ is a copy of $A$ with
comultiplication $\Delta _{\alpha ,\beta }(i_{\alpha \beta
}(a))=i_{\alpha }(\lambda (\beta )(a_{1}))\otimes i_{\beta
}(a_{2})$, counit $\epsilon _{1}(i_{1}(a))=\epsilon (a)$ and the
antipode $S_{\alpha }(i_{\alpha }(a))=$ $i_{\alpha ^{-1}}(\lambda
(\alpha )(s(a)))$ where $a\in A $, $\Delta (a)=a_{1}\otimes a_{2}$
is the given comultiplication in $A$ written in Sweedler`s sigma
notation.\\

Now take $J=\{J_{\alpha }\}_{\alpha \in \pi }$ where for each $\alpha \in
\pi $, $J_{\alpha }$ is a copy of $J_{R}$ where $J_{R}$ is a right ideal of $%
A$ generated by $x-1,y-1,x^{-1}-1$ and $y^{-1}-1$. $J_{R}$ is a coideal
becouse, for example,
\begin{eqnarray*}
\Delta (x-1) &=&x\otimes x-1\otimes 1=(x-1)\otimes x+1\otimes (x-1)\text{ }
\\
\epsilon (x-1) &=&\epsilon (x)-\epsilon (1)=0-0=0,\text{ \ \ \ \ \ \ \ \ \ \
\ \ \ \ \ \ \ \ \ \ \ \ \ \ \ \ \ \ \ \ \ \ \ \ \ \ \ \ \ \ \ \ \ \ \ \ \ \
\ \ \ \ \ \ \ \ \ \ \ \ \ \ \ \ \ \ \ \ \ \ \ \ \ \ \ \ \ \ \ \ \ \ \ \ \ \
\ \ \ \ \ \ \ \ \ \ \ \ \ }
\end{eqnarray*}%
$J$ is a family of right ideals and also $J$ is $\pi -$ coideal. Then, by
theorem 1.2, $C=H/J$ is $\pi -$ coalgebra such that $\overline{\sigma }=\{%
\overline{\sigma }_{\alpha }=\sigma :H_{\alpha }\rightarrow C_{\alpha
}=H_{\alpha }/J_{\alpha }\}_{\alpha \in \pi }$ is $\pi -$coalgebra
epimorphism map. Furthermore, every $C_{\alpha }$ is $H_{\alpha }-$ module
such that every $\overline{\sigma }_{\alpha }$ is module map.\newline

Now take $V=\{V_{\alpha }\}_{\alpha \in \pi }$ is simple $\pi -$comodule
over $C$ by $\rho =\{\rho _{\alpha ,\beta }\}_{\alpha ,\beta \in \pi }$,
this implies that $V_{1}$ is simple $C_{1}-$comodule by $\rho _{1,1}$. Since
$C_{1}$ is cocommutative, then $V_{1}$ has one dimensional. So, we can
suppose that $\rho _{1,1}(w)=w\otimes 1$. The induced representation in this
case is $Ind(\rho )=\{Ind(\rho )_{\alpha }\}$ where $Ind(\rho )_{\alpha }$
is generated by $x,x^{-1},y$ and $y^{-1}$. Also, $Ind(\rho )_{\alpha
}\subseteq H_{\alpha }$. Now, take $B=\{B_{\alpha }\}_{\alpha \in \pi }$
where $B_{\alpha }$ is a subalgebra of $H_{\alpha }$ generated by $x,$ and $%
x^{-1}$. Clear $B$ is subHopf $\pi -$coalgebra of $H$ and every ${B_{\alpha
}\subseteq Ind(\rho )_{\alpha }}$ and every element in $B_{\alpha }$ is
group like, the $B$ is non trivial $\pi -$ subcomodule of $Ind(\rho )$ and
hence $Ind(\rho )$ is not simple $\pi -$comodule over $H$.\newline

Now, we return to study the equivalent and direct sum of coinduced
representation.\newline

\begin{thm}
Let $(V,\rho )$ and $(W,\rho ^{\prime })$ be a right $\pi -$comodules over a
$\pi -$coisotropic quantum subgroup $(C,\sigma )$. If there exist a vector
space isomorphism $F_{1}:$ $V_{1}\rightarrow W_{1}$ such that
\begin{equation*}
\text{ \ \ \ \ \ \ \ \ \ \ \ \ \ \ \ \ \ \ \ \ \ \ \ \ \ \ \ \ \ \ \ }\psi
_{1,1}F_{1}=(F_{1}\otimes I_{1})\rho _{1,1}\;\text{\ \ \ \ \ \ \ \ \ \ \ \ \
\ \ \ \ \ \ \ \ \ \ \ \ \ \ \ \ \ \ \ \ \ \ (1)}
\end{equation*}%
then the coinduced representations $Coind(\rho )$ and $Coind(\psi )$ are
equivalent right $\pi -$comodules over $H$.\newline
\end{thm}

\begin{proof} Let us consider the vector space isomorphisms $\overline{F}%
=\{\overline{F}_{\alpha }:Coind(\psi )_{\alpha }\rightarrow
Coind(\rho )_{\alpha }\}_{\alpha \in \pi }$, where
\begin{equation*}
\overline{F}_{\alpha }(f)=f\circ F_{1}\text{ \ \ \ for all \ }f\in
Coind(\psi )_{\alpha }.\newline
\end{equation*}

Firstly, we must be prove $\overline{F}_{\alpha }(f)=f\circ
F_{1}\in Coind(\rho )_{\alpha }$, for $\alpha \in \pi $ for
$\alpha \in \pi $, $f\in Coind(\psi )_{\alpha }$.
\begin{eqnarray*}
&&L_{1,\alpha }\overline{F}_{\alpha }(f) \\
&=&L_{1,\alpha }f\circ F_{1} \\
&=&(I\otimes f)\psi _{1,1}F_{1} \\
&=&(I\otimes f)(I\otimes F_{1})\rho _{1,1}\text{ \ \ \ \ \ \ \ \ \
\ \ \ \ \ \ \ \ \ \ \ \ \ \ \ \ \ \ \ \ \ \ \ \ \ \ \ \ \ \ \ \ \
\ \ \ \ \ \ \ \ \ \ \ \ \ \ \ \ \ \ \ \ \ \ \ by (1)} \\
&=&(I\otimes f\circ F_{1})\rho _{1,1} \\
&=&(I\otimes \overline{F}_{\alpha }(f))\rho _{1,1}.\text{ \ \ \ \
\ \ \ \ \ \ \ \ \ \ \ \ \ \ \ \ \ \ \ \ \ \ \ \ \ \ \ \ \ \ \ \ \
\ \ \ \ \ \ \ \ \ \ \ \ \ \ \ \ \ \ \ \ \ \ \ \ \ \ \ \ \ \ \ \ \
\ \ \ \ \ \ \ \ \ \ \ \ \ \ \ \ \ \ \ \ \ \ \ \ \ \ \ \ \ \ \ \ \
\ \ \ \ \ \ \ \ \ \ \ \ \ \ }
\end{eqnarray*}

Now, we will prove the following diagram commutes%
\begin{equation*}
\begin{array}{ccc}
Coind(\rho )_{\alpha \beta } & \overset{\overline{F}_{\alpha \beta }}{%
\longrightarrow } & Coind(\rho )_{\alpha \beta } \\
\Omega _{\alpha ,\beta }\downarrow &  & \downarrow \overline{\Omega }%
_{\alpha ,\beta } \\
Coind(\psi )_{\alpha }\otimes H_{\beta } &
\underset{\overline{F}_{\alpha }\otimes I_{\beta
}}{\longrightarrow } & Coind(\rho )_{\alpha }\otimes
H_{\beta }%
\end{array}%
\end{equation*}%
For $h\in Coind(\psi )_{\alpha \beta }$
\begin{eqnarray*}
&&\overline{\Omega }_{\alpha ,\beta }\overline{F}_{\alpha \beta }(h) \\
&=&\overline{\Omega }_{\alpha ,\beta }(h\circ F_{1}) \\
&=&\sim (I_{\alpha }\otimes g^{\beta })\Delta _{\alpha ,\beta
}(h\circ
F_{1})\otimes e^{\beta } \\
&=&\overline{\Omega }_{\alpha }(\sim (I_{\alpha }\otimes g^{\beta
})\Delta
_{\alpha ,\beta }h)\otimes e^{\beta } \\
&=&(\overline{\Omega }_{\alpha }\otimes I_{\beta })((\sim
(I_{\alpha
}\otimes g^{\beta })\Delta _{\alpha ,\beta }h\otimes e^{\beta }) \\
&=&(\overline{\Omega }_{\alpha }\otimes I_{\beta })\Omega _{\alpha
,\beta }(h).\text{ \ \ \ \ \ \ \ \ \ \ \ \ \ \ \ \ \ \ \ \ \ \ \ \
\ \ \ \ \ \ \ \ \ \ \ \ \ \ \ \ \ \ \ \ \ \ \ \ \ \ \ \ \ \ \ \ \
\ \ \ \ \ \ \ \ \ \ \ \ \ \ \ \ \ \ \ \ \ \ \ \ \ \ \ \ \ \ \ \ \
\ \ \ \ \ \ \ \ \ \ \ \ \ \ \ \ \ }
\end{eqnarray*}%
Now we will prove $\overline{F}_{\alpha }$ is an isomorphism. Let
$t\in Coind(\rho )_{\alpha }$, there exist $t\circ F_{1}^{-1}\in
Coind(\psi
)_{\alpha }$ such that $\overline{F}_{\alpha }(t\circ \overline{F}%
_{1}^{-1})=t\overline{F}_{1}^{-1}F_{1}=t.$ Therefore
$\overline{F}_{\alpha }$ is onto, by finite dimensional
$\overline{F}_{\alpha }$ is an isomorphism.\\

The remaining is prove $t\circ F_{1}^{-1}\in Coind(\psi )_{\alpha }$, i.e., $%
L_{1,\alpha }t\circ F_{1}^{-1}=(I\otimes t\circ F_{1}^{-1})\psi
_{1,1}$ i.e., $(I\otimes t)\rho _{1,1}F_{1}^{-1}=(I\otimes t\circ
F_{1}^{-1})\psi
_{1,1}$. We have $\psi _{1,1}F_{1}=(I\otimes F_{1})\rho _{1,1}$ implies $%
(I\otimes F_{1}^{-1})\psi _{1,1}F_{1}=\rho _{1,1}$ that is means
\begin{equation*}
(I\otimes F_{1}^{-1})\psi _{1,1}F_{1}=\rho _{1,1}F_{1}^{-1}F_{1}
\end{equation*}%
and hence
\begin{equation*}
\lbrack (I\otimes F_{1}^{-1})\psi _{1,1}-\rho
_{1,1}F_{1}^{-1}]F_{1}=0
\end{equation*}%
implies that
\begin{equation*}
\lbrack (I\otimes F_{1}^{-1})\psi _{1,1}-\rho _{1,1}F_{1}^{-1}]F_{1}(v)=0%
\text{ \ \ \ \ \ \ for all \ \ }v\in V_{1}.
\end{equation*}%
Since $F_{1}$ is onto, then
\begin{equation*}
\lbrack (I\otimes F_{1}^{-1})\psi _{1,1}-\rho
_{1,1}F_{1}^{-1}](w)=0\text{ \ \ \ for all \ \ \ \ }w\in W_{1},
\end{equation*}%
that is $(I\otimes F_{1}^{-1})\psi _{1,1}-\rho _{1,1}F_{1}^{-1}=0$ i.e., $%
(I\otimes F_{1}^{-1})\psi _{1,1}=\rho _{1,1}F_{1}^{-1}$ therefore
\begin{equation*}
(I\otimes t)(I\otimes F_{1}^{-1})\psi _{1,1}=(I\otimes t)\rho
_{1,1}F_{1}^{-1},
\end{equation*}%
then
\begin{equation*}
(I\otimes t\circ F_{1}^{-1})\psi _{1,1}=L_{1,\alpha }t\circ
F_{1}^{-1}.
\end{equation*}%
Thus $Coind(\psi )$ and $Coind(\rho )$ are equivalent.
\end{proof}

\begin{cor}
If $(V, \rho)$ and $(W, \psi)$ are equivalent right $\pi-$comodules over $C$%
, then the coinduced representations $Coind(\rho)$ and $Coind(\psi)$ are
equivalent right $\pi-$comodules over $H$.\newline
\end{cor}

\begin{thm}
Suppose $(V,\rho )$ and $(W,\psi )$ are left $\pi -$comodules over $C$, then
\begin{equation*}
Coind(\rho \oplus \psi )\cong Coind(\rho )\oplus Coind(\psi ).\newline
\end{equation*}
\end{thm}

\begin{proof} Since $(V,\rho )$ and $(W,\psi )$ are
left $\pi -$comodules over $C$, then $(Coind(\rho ),\Omega )$ and
$(Coind(\psi ),\Upsilon )$ are right $\pi -$comodules over $H$ and
hence, by lemma 5.3, $(Coind(\rho )\oplus (Coind(\psi ))$ is right
$\pi -$comodule over $H$ by $\Lambda =\{\Lambda _{\alpha ,\beta
}\}$ where
\begin{equation*}
\Lambda _{\alpha ,\beta }=(l_{\alpha }^{Coind(\rho )}\otimes
I_{\beta })\Omega _{\alpha ,\beta }P_{\alpha \beta }^{Coind(\rho
)}+(l_{\alpha }^{Coind(\psi )}\otimes I_{\beta })\Upsilon _{\alpha
,\beta }P_{\alpha \beta }^{Coind(\psi )}.
\end{equation*}

Also, by lemma 5.3, $V\oplus W=\{V_{\alpha }\oplus W_{\beta }\}$ is left $%
\pi -$comodule over $C$ by
\begin{equation*}
\rho \oplus \psi =\{(\rho \oplus \psi )_{\alpha ,\beta }=\rho
_{\alpha ,\beta }\oplus \psi _{\alpha ,\beta }:V_{\alpha \beta
}\oplus W_{\alpha \beta }\rightarrow C_{\alpha }\otimes (V_{\beta
}\oplus W_{\beta })\}_{\alpha ,\beta \in \pi }
\end{equation*}
where
\begin{equation*}
\rho _{\alpha ,\beta }\oplus \psi _{\alpha ,\beta }=(I_{\alpha
}\otimes l_{\beta }^{V})\rho _{\alpha ,\beta }P_{\alpha \beta
}^{V}+(I_{\alpha }\otimes l_{\beta }^{W})\psi _{\alpha ,\beta
}P_{\alpha \beta }^{W}.
\end{equation*}

Then $Coind(\rho\oplus \psi)$ is right $\pi-$comodule over $H$ by $%
\Gamma=\{\Gamma_{\alpha,\beta}\}$.

Now, we define $\chi =\{\chi _{\alpha }:Coind(\rho \oplus \psi
)_{\alpha }\rightarrow Coind(\rho )_{\alpha }\oplus Coind(\psi
)_{\alpha }\}_{\alpha \in \pi }$, where
\begin{equation*}
\chi_{\alpha }(g)=g\circ l_{1}^{V}+g\circ l_{1}^{W}\;\;\;\text{for all}%
\;\;\;g\in Coind(\rho \oplus \psi )_{\alpha }.
\end{equation*}

Firstly, we will prove $g\circ l_{1}^{V}\in Coind(\rho )_{\alpha
}$
\begin{eqnarray*}
&&L_{1,\alpha }g\circ l_{1}^{V}\\
&=&(L_{1,\alpha }g)l_{1}^{V} \\
&=&(I_{1}\otimes g)(\rho _{1,1}\otimes \psi _{1,1})l_{1}^{V} \\
&=&(I_{1}\otimes g)[(I_{1}\otimes l_{1}^{V})\rho
_{1,1}P_{1}^{V}+(I_{1}\otimes l_{1}^{W})\psi _{1,1}P_{1}^{W}]l_{1}^{V} \\
&=&(I_{1}\otimes g)(I_{1}\otimes l_{1}^{V})\rho _{1,1}P_{1}^{V}l_{1}^{V} \\
&=&(I_{1}\otimes g)(I_{1}\otimes l_{1}^{V})\rho _{1,1} \\
&=&(I_{1}\otimes g\circ l_{1}^{V})\rho _{1,1}.\text{ \ \ \ \ \ \ \
\ \ \ \ \ \ \ \ \ \ \ \ \ \ \ \ \ \ \ \ \ \ \ \ \ \ \ \ \ \ \ \ \
\ \ \ \ \ \ \ \ \ \ \ \ \ \ \ \ \ \ \ \ \ \ \ \ \ \ \ \ \ \ \ \ \
\ \ \ \ \ \ \ \ \ \ \ \ \ \ \ \ \ \ \ \ \ \ \ \ \ \ \ \ \ \ \ \ \
\ \ \ \ \ \ \ \ \ \ \ \ \ \ \ \ \ \ \ \ \ \ \ }
\end{eqnarray*}

Similar $g\circ l^{W}_{1} \in Coind(\psi)_{\alpha}$, and hence $%
\chi_{\alpha}(g)\in Coind(\rho)_{\alpha}\oplus Coind(\psi)_{\alpha}$.%
\newline

We will prove $\chi _{\alpha }$ is 1-1 (and hence isomorphism by
finite dimensional). Suppose $\chi _{\alpha }(g)=0$, then $(\chi
_{\alpha
}(g))(v_{1}+w_{1})=0$ for all $v_{1}\in V_{1}$ and $w_{1}\in W_{1}$ implies $%
g(v_{1}+w_{1})=0$ for all $v_{1}\in V_{1}$ and $w_{1}\in W_{1}$,
then $g=0.$

Finally, we will prove that the following diagram commutes%
\begin{equation*}
\begin{array}{ccc}
Coind(\rho \oplus \psi )_{\alpha \beta } & \overset{\chi _{\alpha ,\beta }}%
{\longrightarrow } & Coind(\rho )_{\alpha \beta }\oplus Coind(\psi
)_{\alpha
\beta } \\
\Gamma _{\alpha ,\beta }\downarrow &  & \downarrow \Lambda
_{\alpha ,\beta }
\\
Coind(\rho \oplus \psi )_{\alpha }\otimes H_{\beta } &
\underset{\chi _{\alpha }\otimes I_{\beta }}{\longrightarrow } &
Coind(\rho )_{\alpha
}\oplus Coind(\psi )_{\alpha }\otimes H_{\beta }%
\end{array}%
\end{equation*}
\begin{eqnarray*}
&&(\chi _{\alpha }\otimes I_{\beta })\Gamma _{\alpha ,\beta }(g) \\
&=&(\chi _{\alpha }\otimes I_{\beta })[\sim (I_{\alpha }\otimes
f^{\beta
})\Delta _{\alpha ,\beta }g\otimes e^{\beta }] \\
&=&[(\sim (I_{\alpha }\otimes f^{\beta })\Delta _{\alpha ,\beta
}g)l_{1}^{V}+(\sim (I_{\alpha }\otimes f^{\beta })\Delta _{\alpha
,\beta
}g)l_{1}^{W}]\otimes e^{\beta } \\
&=&[\sim (I_{\alpha }\otimes f^{\beta })\Delta _{\alpha ,\beta
}(g\circ l_{1}^{V}+g\circ l_{1}^{W})]\otimes e^{\beta }\text{ \ \
\ \ \ \ \ \ \ \ \ \ \ \ \ \ \ \ \ \ \ \ \ \ \ \ \ \ \ \ \ \ \ \ \
\ \ \ \ \ \ \ \ \ \ \ \ \ \ \ \ \ \ \ \ \ \ \ \ \ \ \ \ \ \ \ \ \
\ \ \ \ \ \ \ \ \ \ \ \ \ \ \ \ \ \ \ \ \ \ \ \ \ \ \ \ \ \ \ \ \
\ \ \ \ \ \ \ \ \ \ \ \ \ \ \ \ \ \ \ \ \ \ \ \ \ \ \ \ }
\end{eqnarray*}%
On other hand
\begin{eqnarray*}
&&\Lambda _{\alpha ,\beta }\chi _{\alpha ,\beta }(g) \\
&=&[(l_{\alpha }^{Coind(\rho )}\otimes I_{\beta })\Omega _{\alpha
,\beta }P_{\alpha \beta }^{Coind(\rho )}+(l_{\alpha }^{Coind(\psi
)}\otimes I_{\beta })\Upsilon _{\alpha ,\beta }P_{\alpha \beta
}^{Coind(\psi
)}](gl_{1}^{V}+gl_{1}^{W}) \\
&=&(l_{\alpha }^{Coind(\rho )}\otimes I_{\beta })\Omega _{\alpha
,\beta }g\circ l_{1}^{V}+(l_{\alpha }^{Coind(\psi )}\otimes
I_{\beta })\Upsilon
_{\alpha ,\beta }g\circ l_{1}^{W} \\
&=&(l_{\alpha }^{Coind(\rho )}\otimes I_{\beta })[\sim (I_{\alpha
}\otimes f^{\beta })\Delta _{\alpha ,\beta }g\circ
l_{1}^{V}\otimes e^{\beta }]+(l_{\alpha }^{Coind(\psi )}\otimes
I_{\beta })[\sim (I_{\alpha }\otimes
f^{\beta })\Delta _{\alpha ,\beta }g\circ l_{1}^{W}\otimes e^{\beta }] \\
&=&[\sim (I_{\alpha }\otimes f^{\beta })\Delta _{\alpha ,\beta
}g\circ l_{1}^{V}\oplus 0]\otimes e^{\beta }+[0\oplus \sim
(I_{\alpha }\otimes
f^{\beta })\Delta _{\alpha ,\beta }g\circ l_{1}^{W}]\otimes e^{\beta } \\
&=&\sim (I_{\alpha }\otimes f^{\beta })\Delta _{\alpha ,\beta
}(g\circ l_{1}^{V}+g\circ l_{1}^{W})\otimes e^{\beta }.
\end{eqnarray*}
\end{proof}

\begin{cor}
Let $H$ be a Hopf $\pi -$coalgebra, $(C,\sigma )$ a $\pi -$coisotropic
quantum subgroup of $H$ and $(V,\rho )$ a right $\pi -$comodules over $C$.
If $Coind(\rho )$ is simple, then $\rho $ is simple.\newline
\end{cor}

\end{document}